\documentclass[reqno, a4paper,12pt]{amsart} 
\usepackage{latexsym}
\usepackage[english]{babel}
\usepackage{fancyhdr}
\usepackage[mathscr]{eucal}
\usepackage{amsmath}
\usepackage{mathrsfs}
\usepackage{amsthm}
\usepackage{amsfonts}
\usepackage{amssymb}
\usepackage{amscd}
\usepackage{bbm}
\usepackage{graphicx}
\usepackage{graphics}
\usepackage{latexsym}
\usepackage{color}
\usepackage{ascmac}

\theoremstyle{plain}
\newtheorem{theorem}{Theorem}[section]
\newtheorem{remark}[theorem]{Remark}%[section]
\newtheorem{lemma}[theorem]{Lemma}
\newtheorem{corollary}[theorem]{Corollary}
\newtheorem{proposition}[theorem]{Proposition}

\theoremstyle{definition}
\newtheorem{definition}[theorem]{Definition}

\newcommand {\absleq} {{\leq_{|\, \cdot\, |}\, }}

\def\Lg {{\mathcal L}}
\def\Pg {{\mathcal P}}
\def\Wg {{\mathcal W}}
\def\Xg {{\mathcal X}}
\def\Yg {{\mathcal Y}}
\def\tRg {{\tilde{\mathcal R}}}
\def\tT {{\tilde{T}}}
\numberwithin{equation}{section}
\def\HL {{L^2(\R^2)}}

\def\ben{\begin{enumerate}}
\def\een{\end{enumerate}}
\def\bgdf{\begin{definition}}
\def\eddf{\end{definition}}
\def\bglm{\begin{lemma}}
\def\edlm{\end{lemma}}
\def\bgpf{\begin{proof}}
\def\edpf{\end{proof}}
\def\bgth{\begin{theorem}}
\def\edth{\end{theorem}}
\def\bgprop{\begin{proposition}}
\def\edprop{\end{proposition}}
\def\bgrm{\begin{remark}}
\def\edrm{\end{remark}}
\def\bgcor{\begin{corollary}}
\def\edcor{\end{corollary}}

\def\lbeq(#1){\label{eqn:#1}}
\def\refeq(#1){{\rm (\ref{eqn:#1})}}
\def\refeqs(#1,#2){{\rm (\ref{eqn:#1}) and (\ref{eqn:#2})}}
\def\refeqss(#1,#2,#3){{\rm (\ref{eqn:#1}),\ (\ref{eqn:#2}) and (\ref{eqn:#3})}}
\def\lbth(#1){\label{th:#1}}
\def\refth(#1){{\rm Theorem \ref{th:#1}}}
\def\refths(#1,#2){{\rm Theorems \ref{th:#1} and \ref{th:#2}}}
\def\refthb(#1){{\bf Theorem \ref{th:#1}}}
\def\lblm(#1){\label{lm:#1}}
\def\reflm(#1){{\rm Lemma \ref{lm:#1}}}
\def\reflms(#1,#2){{\rm Lemmas \ref{lm:#1} and \ref{lm:#2}}}
\def\reflmss(#1,#2,#3){{\rm Lemmas \ref{lm:#1}, \ref{lm:#2} and \ref{lm:#3}}}
\def\reflmsss(#1,#2,#3,#4){{\rm Lemmas \ref{lm:#1},\, \ref{lm:#2},\, \ref{lm:#3} and \ref{lm:#4}}}
\def\reflmb(#1){{\bf Lemma \ref{lm:#1}}}
\def\lbprop(#1){\label{prp:#1}}
\def\refprop(#1){{\rm Proposition \ref{prp:#1}}}
\def\refprops(#1,#2,#3,#4){{\rm Propositions \ref{prp:#1},\, \ref{prp:#2},
\, \ref{prp:#3} \, and \ref{prp:#4}}}
\def\refpropb(#1){{\bf Proposition \ref{prp:#1}.}}
\def\lbcor(#1){\label{cor:#1}}
\def\refcor(#1){{\rm Corollary \ref{cor:#1}}}
\def\refcors(#1,#2){{\rm Corollaries \ref{cor:#1} and \ref{cor:#2}}}
\def\lbrm(#1){\label{rm:#1}}
\def\refrm(#1){{\rm Remark \ref{rm:#1}}}
\def\refrmss(#1,#2,#3){{\rm Remark \ref{rm:#1},\, \ref{rm:#2}\, and \ref{rm:#3}}}
\def\lbass(#1){\label{ass:#1}}
\def\refass(#1){{\rm Assumption \ref{ass:#1}}}
\def\lbdf(#1){\label{df:#1}}
\def\refdf(#1){{\rm Definition \ref{df:#1}}}
\def\refdfs(#1,#2){{\rm Definitions \ref{def:#1} and \ref{def:#2}}}
\def\lbsec(#1){\label{s:#1}}
\def\refsec(#1){{\rm \S\ref{s:#1}}}
\def\lbsubsec(#1){\label{ss:#1}}
\def\refsubsec(#1){{\rm \S\ref{ss:#1}}}

\def\Ag{{\mathcal A}}
\def\Bg{{\mathcal B}}
\def\Fg{{\mathcal F}}
\def\Gg{{\mathcal G}}
\def\Mg{{\mathcal M}}
\def\Ng{{\mathcal N}}

\def\Og{{\mathcal O}}

\newcommand{\lam}{\lambda}

\def\Bb{{\bf B}}

\def\ph{{\varphi}}
\def\bqn{\begin{equation}}
\def\eqn{\end{equation}}
\def\C{{\mathbb C}}
\def\N{{\mathbb N}}

 \def\Cb{{\overline{\mathbb C}}}

\def\R{{\mathbb R}}
\def\Rg {{\mathcal R}}
\def\a{\alpha}
\def\b{\beta}
\def\c{\gamma}

\def\d{\delta}
\def\Cg{{\mathcal C}}
\def\Dg{{\mathcal D}}
\def\Eg{{\mathcal E}}
\def\Sg{{\mathcal S}}

\def\p{\psi}

\def\ep{\varepsilon}
\def\z{\zeta}

\def\th{\theta}

\def\k{\kappa}
\def\m{\mu}

\def\r{\rho}
\def\s{\sigma}
\def\t{\tau}
\def\w{\omega}
\def\W{\Omega}
\def\Hg {{\mathcal H}}
\def\Zg {{\mathcal Z}}
\def\la{\langle}
\def\ra{\rangle}

\def\lap{\Delta}
\def\ax{{\la x \ra}}
\def\ay{{\la y \ra}}
\def\az{{\la z \ra}}
\def\pa{{\partial}}

\def\bglm{\begin{lemma}} 
\def\edlm{\end{lemma}} 
\def\br{\begin{array}}
\def\er{\end{array}}

\def\Ker{\rm Ker\,}
\def\rank{{\rm rank\,}}
\def\low{{\rm low}}

\begin{document}

\title[$L^p$-boundedness of wave operators in $\R^2$]
{The $L^p$-boundedness of wave operators for 
two dimensional Schr\"odinger operators 
with threshold singularities}

\author[K.~Yajima]{Kenji Yajima}
\address{Department of Mathematics \\ Gakushuin University 
\\ 1-5-1 Mejiro \\ Toshima-ku \\ Tokyo 171-8588 (Japan). \\ 
\footnote{Supported by JSPS grant in aid for scientific research No. 19K03589}}
\email{kenji.yajima@gakushuin.ac.jp}

\begin{abstract} 
We generalize the recent result of Erdo{\u g}an, Goldberg and Green 
on the $L^p$-boundedness of wave operators for two dimensional Schr\"odinger 
operators and prove that they are bounded in $L^p(\R^2)$ for all $1<p<\infty$ 
if and only if the Schr\"odinger operator possesses no $p$-wave threshold 
resonances, viz. Schr\"odinger equation $(-\lap + V(x))u(x)=0$ possesses no 
solutions which satisfy $u(x)= (a_1x_1+a_2 x_2)|x|^{-2}+ o(|x|^{-1})$ 
as $|x|\to \infty$  for an $(a_1, a_2) \in \R^2\setminus \{(0,0)\}$ and, 
otherwise, they are bounded in $L^p(\R^2)$ for $1<p\leq 2$ and 
unbounded for $2<p<\infty$.  We present also a new proof for  
the known part of the result. 
\end{abstract}

\date{}

\maketitle

\section{Introduction and main result}\label{sec:theorems} 

Let $H_0= -\lap$, $D(H_0)=W^{2,2}(\R^2)$ be the free Schr\"odinger operator 
on $\R^d$, $d\geq 1$ and $V(x)$ a real measurable function on $\R^d$. 
Suppose that, for some $\c>1/2$,
$\ax^{\c}|V(x)|^{1/2}(H_0+ 1)^{-\frac12}$, $\ax=(1+|x|^2)^\frac12$, 
is a compact operator 
on $L^{2}(\R^d)$. Then, Schr\"odinger operator 
$H=-\lap + V$ defined via the quadratic form is selfadjoint (\cite{RS2});  
the spectrum $\s(H)$ consists of the absolutely continuous part 
$[0,\infty)$ and the point spectrum which is discrete 
in $\R \setminus \{0\}$ (\cite{Ag}); $L^2(\R^d)$ is the orthogonal sum of 
the absolutely continuous subspace $L^2_{ac}(\R^d)$ for $H$ 
and the space of eigenfunctions of $H$. 
The scattering theory compares the large time behavior of scattering 
solutions $e^{-itH}\ph$, $\ph \in L^2_{ac}(\R^d)$ of the time dependent 
Schr\"odinger equation: 
\bqn \lbeq(Sch) 
i\pa_t u(t) = Hu(t), \quad u(0) = \ph \in L^2_{ac}(\R^d)
\eqn 
with that of free solutions $e^{-itH_0}\ph_0$ and wave operators 
are defined by the strong limits:   
\bqn \lbeq(limit)
W_{\pm}= \lim_{{ t\to\pm \infty}} e^{itH }e^{-itH_0}\,.
\eqn 
It is well known (\cite{Ag,Ku}) that $W_{\pm}$ exist, 
${\rm Image}\ W_\pm =L^2_{ac}(H)$ and,   
all scattering solutions $e^{-itH}P_{ac}(H)\ph$  
become asymptotically free in the remote past and far future:
\[
\lim_{{ t\to\pm \infty}} \|e^{-itH}P_{ac}(H)\ph - e^{-itH_0}\ph_\pm \| =0, 
\quad \ph_\pm = W_{\pm}^\ast \ph,
\]
where $P_{ac}(H)$ is the orthogonal projection onto $L^2_{ac}(H)$. 

Wave operators $W_\pm$ satisfy the intertwining property, viz.  
for Borel functions $f(\lam)$ of $\lam\in \R^1$  
\bqn \lbeq(inter) 
f(H)P_{ac}(H)= W_\pm f(H_0)W_\pm^\ast  
\eqn 
and various properties of $f(H)P_{ac}(H)$ may be derived from those of 
the Fourier multiplier $f(H_0)$ if $W_\pm$ satisfy appropriate properties. 
If $W_\pm$ are bounded in $L^p(\R^d)$ for $p\in I\subset [1,\infty]$, then 
\refeq(inter) produces a set of estimates that 
for $\{p,q\} \in I\times I^\ast$, $I^\ast=\{p/p-1, p \in I\}$,   
\begin{gather*}
\|f(H)P_{ac}(H)\|_{\Bb(L^q, L^p)} \leq C \|f(H_0)\|_{\Bb(L^q, L^p)}, \\
\|f(H_0)\|_{\Bb(L^p, L^q)} 
\leq C^{-1} \|f(H)P_{ac}(H) \|_{\Bb(L^p, L^q)} 
\end{gather*} 
where $C=C_{p,q}$ is independent of $f$. Here $\Bb(X,Y)$ is the Banach space 
of bounded operators from $X$ to $Y$ 
and $\Bb(X)=\Bb(X,X)$. Such estimates are 
very useful and have many applications (e.g. \cite{Schlag-rew}). 

Thus, the problem of whether or not $W_\pm$ are bounded in $L^p(\R^d)$ has 
attracted interest of many authors and various results have been obtained 
under various assumptions. We briefly review here some results under 
simplified assumption that $|V(x)|\leq C\ax^{-\s}$ 
for a sufficiently large $\s>2$. We need some notation: 
Multiplication operator $M_m$ with a function $m(x)$ is often denoted 
simply by $m$. 
$\C^{+}=\{\lam\in \C \colon \Im \lam>0\}$ and 
$\Cb^{+}= \C^{+} \cup \R$; $G_0(\lam)= (H_0-\lam^2)^{-1}$ for $\lam \in \C^{+}$; 
$\mathcal{G}_\lam(x)$ is the 
convolution kernel of $G_0(\lam)$: $G_0(\lam)u(x)= (\Gg_\lam\ast u)(x)$; 
if $d=2$, $\Gg_\lam(x)$ is given by 
\begin{equation}\lbeq(1-1)
\mathcal{G}_\lam(x)
= \frac1{(2\pi)^2}\int_{\R^2}\frac{e^{ix\xi}d\xi}{\xi^2-\lam^2}= 
\frac{i}{4} H_0^{(1)}(\lam |x|)\,,
\end{equation}
where $H_0^{(1)}(z)$ is the Hankel function of the first kind.
\[
\la u, v \ra = \int_{\R^d} u(x) v(x) dx 
\]
without complex conjugation.  
\[
U(x) = \left\{\br{ll} 1 \ &\ \mbox{if}\ \ V(x)\geq 0, \\ -1  
&\ \mbox{if}\ \ V(x)<0, \er \right. 
\quad v(x) = |V(x)|^{1/2}, 
\quad w(x)= U(x) v(x)\, .
\] 
The limiting absorption priniciple (e.g. \cite{Ag,Ku}) implies that 
the holomorphic function 
$\C^{+} \ni \lam \to vG_0({\lam})v \in \Bb_c(L^2(\R^d))$, 
the space of compact operators on $L^2(\R^d)$, has a continuous 
extension to $\Cb^{+}$ 
if $d \geq 3$ and to $\lam \in \Cb^{+}\setminus \{0\}$ if $d=1$ and $d=2$;  
\bqn 
\lbeq(Mdef)
M(\lam) \stackrel{\rm def}{=} U + vG_0({\lam})v, \quad 
\lam \in \Cb^{+}\setminus\{0\}  
\eqn 
is invertible unless $\lam^2$ is an eigenvalue of $H$ and the absence of 
positive eigenvalues (\cite{Kato-e,Koch-Tataru}) impies that 
$M(\lam)^{-1}\in \Bb(L^2(\R^d))$ exists for any $\lam\in \R \setminus\{0\}$. 

For $d \geq 3$, we say $H$ is regular at zero if $M(0)^{-1}\in \Bb(L^2)$ 
exists and is singular at zero otherwise. 
It is known that $H$ is singular if and only if for 
some $1/2<\c<\s-1/2$
\bqn \lbeq(Ngdef)
\Ng\stackrel{\rm def}{=} 
\{u\in \ax^{\c}L^2(\R^d) \colon (-\lap + V(x))u(x)=0\}\not=\{0\}; 
\eqn 
$\Ng$ is independent of $\c$ and $u\in \Ng$ is called 
threshold resonance; $u(x)=O(|x|^{2-d})$ as $|x|\to \infty$ 
and $u(x) = O(|x|^{1-d-j})$ if it satisfies $\la x^\a V, u\ra=0$ for 
$|\a|\leq j$, $j=0,1,\dots$;  $u \in \Ng$ is an eigenfunction of $H$ 
with eigenvalue $0$ if $d\geq 5$. 
When $d=1$ or $2$, we say $H$ is 
regular at zero if $\Ng_\infty =\{0\}$ where 
\bqn \lbeq(Nginfty)
\Ng_\infty =\{u\in L^\infty(\R^d) \colon (-\lap + V(x))u(x)=0\} 
\eqn 
and $H$ is singular at zero if otherwise. 
If $d=2$, it is known (see \reflm(resonances)) that $u \in \Ng_\infty$ 
satisfies for constants $c,b_1, b_2$ and $\ep>0$ that 
\bqn 
u(x) = c + \frac{b_1x_1+ b_2x_2}{|x|^2} 
+ O(|x|^{-1-\ep}),
\quad (|x|\to\infty) 
\lbeq(reso-intro)
\eqn 
 and $u\in \Ng_\infty\setminus\{0\}$ 
is called $s$-wave resonance if $c\not=0$, $p$-wave resonance if $c=0$ 
but $(b_1, b_2)\not=(0,0)$ and it is a zero energy eigenfunction of $H$ 
if $c=b_1=b_2=0$. 

We list some known results. 
If $d=1$, $W_\pm$ are bounded in $L^p(\R^1)$ for 
$1<p<\infty$ and are in general unbounded in $L^p(\R^1)$ 
for $p=1, \infty$  
(\cite{RW,GY,DF}). 

If $H$ is regular at zero, we have 
almost complete results: If $d=2$, $W_\pm\in \Bb(L^p(\R^2))$ 
for $1<p<\infty$ 
(\cite{Y-2dim,JY-2} but no results for $p=1$ and $p=\infty$); 
if $d \geq 3$, $W_\pm\in \Bb(L^p(\R^d))$ for 
$1\leq p \leq \infty$ (\cite{Y-3dim,Y-even,BS}).

If $H$ is singular at zero, a rather complete result is known for 
$d\not=2,4$. If $d\geq 5$, 
$W_\pm\in \Bb(L^p(\R^d))$ for $1\leq p<d/2$, 
for $1\leq p<d$ if and only if $\la V, u \ra=0$ for all $u\in \Ng$ 
and, for $1\leq p<\infty$ if and only if 
$\la x^\a V, u \ra=0$ for all $|\a|\leq 1$ 
(\cite{GG,KY,FY,Y-odd-sing}). If $d=3$, $W_\pm\in \Bb(L^p(\R^3))$ 
for $1<p<3$; for $p=1$ if and only if all 
$u\in \Ng$ satisfy $\la V, u\ra=0$; for 
$1\leq p<\infty$ if and only if $\la x^\a V, u\ra=0$ for all $|\a|\leq 1$; 
for $p=\infty$ if $\la x^\a V, u\ra=0$ for all $|\a|\leq 2$.
(\cite{Y-3d-sing}).

However, only a partial result is known when $d=2$ or $d=4$: 
For $d=2$, $ W_\pm\in \Bb(L^p(\R^2))$ for $1<p<\infty$ 
if $\Ng_{\infty}$ consists only of zero energy eigenfunctions 
or only of $s$-wave resonances (\cite{EGG}). If $d=4$ and if all $u\in \Ng$ 
satisfy $\la V, u\ra=0$, then $ W_\pm\in \Bb(L^p(\R^4))$ for $1\leq p \leq 4$ 
and, for $1\leq p < \infty$ 
if $\la x^\a V, u\ra=0$ for $|\a|\leq 1$.(\cite{GG-4, JY-4}). 

The purpose of this paper is to prove the following theorem 
for the case $d=2$ which fills the missing part in the results of 
\cite{Y-2dim, JY-2} and \cite{EGG}. 

\bgth \lbth(th-main) Suppose 
$\ax^2 V \in L^{\frac43}(\R^2)$ and $\ax^{\c}|V(x)|\in L^1(\R^2)$ for 
a constant $\c>8$. Then, $W_\pm $ are bounded in $L^p(\R^2)$  for  
$1<p<\infty$ if and only if $H$ has no $p$-wave resonances. 
$W_\pm$ are otherwise bounded in $L^p(\R^2)$ for $1<p\leq 2$ and 
unbounded for $2<p<\infty$. 
\edth 

We shall also give a new proof for the known parts 
(\cite{Y-2dim, JY-2, EGG}) under slightly weaker assumptions.  
The problem for $p=1$ and $p=\infty$, however, is left open. 
A similar result has recently been obtained for Schr\"odinger 
operators with point interactions on $\R^2$
(\cite{CMY}, \cite{CMYE} and \cite{Ya-point}) and the main idea of 
the proof is borrowed from \cite{CMY} and \cite{Ya-point}.

We briefly explain here the basic strategy for the proof of 
\refth(th-main), introducing 
some more notation and displaying the plan of the paper. 
We prove it only for $W_{+}$. The result for $W_{-}$ then follows 
via the complex conjugation ${\mathcal C} u(x) = \overline{u(x)}$: 
$W_{-}= {\mathcal C} W_{+}{\mathcal C}^{-1}$. 
\[
\hat{u}(\xi)=\Fg u(\xi) = \frac1{2\pi}\int_{\R^2} e^{-ix\xi} u(x) dx.\quad 
\check{u}(\xi)=(\Fg^{-1}u)(\xi)
\]
is the Fourier transform and its inverse; 
$\Sg(\R^2)$ is the Schwartz space;
\bqn 
\Dg_\ast= \{u \in \Sg(\R^2) \ | \ \hat u \in C_0^\infty(\R^2 \setminus\{0\})\};
\eqn 
$\Dg_\ast$ is dense in $L^p(\R^2)$, $1 \leq p<\infty$;
$\|u\|_p=\|u\|_{L^p(\R^2)}$; 
\[
\t_y u(x)=u(x-y), \quad y \in \R^2.
\]
For a Borel function $f(\lam)$ of $\lam>0$, $f(|D|)$ is the Fourier 
multiplier 
\[
f(|D|)u(x)= \frac1{2\pi}\int_{\R^2} e^{ix\xi} f(|\xi|) \hat{u}(\xi)d\xi.
\]
We take and fix $\chi \in C_0^\infty(\R)$ such that 
\bqn  
\chi(\lam)=1\ \mbox{for}\ |\lam|\leq 1/2 \ \mbox{and}\  
\chi(\lam)=0 \ \mbox{for}\ |\lam|\geq 1 
\eqn 
and define for $a>0$ 
\bqn 
\chi_{\leq a}(\lam) = \chi(\lam/a), \quad   
\chi_{> a}(\lam) = 1- \chi_{\leq a}(\lam).
\eqn 
We decompose $W_{+}$ into the high and the low energy parts: 
\bqn \lbeq(W-d)
W_{+}= W_{+}\chi_{> 2a}(|D|)+ W_{+}\chi_{\leq 2a}(|D|) 
\eqn 
and study them separately. Here $a>0$ is arbitrary.  

The proof is based on the stationary formula for $W_{+}$ 
(e.g. \cite{Ku}): 
\begin{gather} \lbeq(sta-00)
W_{+}u(x) = u(x) - \W_{+}u(x), \\
\W_{+}u(x) = 
\int_0^\infty 
(G_0(-\lam)vM(\lam)^{-1}v \Pi({\lam})u)(x) \lam d\lam\,, \lbeq(sta-0)
\end{gather} 
where $\Pi({\lam})u(x)=(i\pi)^{-1}(G_0(\lam)-G_0(-\lam))u(x)$ and for $\lam>0$    
\bqn 
\Pi({\lam})u(x)= 
\frac1{2\pi}\int_{{\mathbb S}^1}e^{i\lam{\w}x} \hat{u} (\lam{\w})d\w\, 
= \frac1{2\pi}\int_{{\mathbb S}^1} (\Fg \tau_{-x} {u})(\lam{\w})d\w\,.
\lbeq(Pi-express)
\eqn  
We define $\Pi({\lam})u(x)$ also for $\lam=0$ by \refeq(Pi-express). 
Then, for $f\in C([0,\infty))$,    
\bqn \lbeq(mult)
f(\lam) \Pi(\lam ) u = \Pi(\lam) f(|D|) u, \quad \lam\geq 0\,.
\eqn 
We have $W_{+}\chi_{> 2a}(|D|)u =\chi_{> 2a}(|D|) u + \W_{{\rm high},2a}u$ and 
$W_{+}\chi_{\leq 2a}(|D|) u =\chi_{\leq 2a}(|D|) u + \W_{{\rm low},2a}u$ where   
\begin{gather} \lbeq(sta-0-high)
\W_{{\rm high},2a}u = \int_0^\infty 
G_0(-\lam)vM(\lam)^{-1}v \Pi({\lam})u \chi_{> 2a}(\lam)\lam d\lam\,,
\\
\W_{{\rm low},2a}u = \int_0^\infty 
G_0(-\lam)vM(\lam)^{-1}v \Pi({\lam})u \chi_{\leq 2a}(\lam)\lam d\lam\, .
\lbeq(sta-0-low)
\end{gather}
and we have only to study $\W_{{\rm high},2a}$ and $\W_{{\rm low},2a}$.

In section 2, we state and prove some estimates related to the Hankel 
function and recall some properties of 
the integral operator $K$ introduced in \cite{CMY}:  
\bqn \lbeq(DFEK)
Ku(x)=\frac1{2\pi}\int_0^{+\infty} \mathcal{G}_{-\lambda}(x) \lam 
\left(\int_{{\mathbb S}^1} ({\mathcal F}{u})(\lam \w) d\w \right) d\lambda\,;
\eqn 
$K$ is bounded in $L^p(\R^2)$ for all $1<p<\infty$ and is closely connected 
to $\W_{+}$:   
\bqn \lbeq(how-0)
(\t_y K \t_{-z}u)(x)= \int_0^\infty \Gg_{-\lam}(x-y) \Pi(\lam)u(z) \lam d\lam\,. 
\eqn 
We shall give a simpler proof of the 
$L^p(\R^2)$-boundedness of $K$ in section 2. 
 
In section 3 we prove the following \refprop(funda) which will be 
repeatedly used for estimating the operators produced by \refeq(sta-0) 
by the less singular part of $M(\lam)^{-1}$ as $\lam\to 0$. 
$\Hg_2$ is the Hilbert space of 
Hilbert-Schmidt operators on $\HL$. 
$\Lg_1$ is the Banach space of integral operators $T$ with 
$T(x,y) \in L^1(\R^2 \times \R^2)$ with the norm 
$\|T\|_{\Lg_1}= \|T\|_{L^1(\R^2 \times \R^2)}$. 

\bgdf \lbdf(term) We use the following terminology. We say:
\begin{itemize}
\item[{\rm (1)}]  $X$ is a {\it good operator} 
if $X \in \Bb(L^p(\R^2))$ for all $1<p<\infty$. 

\item[{\rm (2)}] $f\in C^2((0,\infty))$ 
is a {\it good multiplier} if it satisfies 
$|f^{(j)}(\lam)|\leq C_j \lam^{-j}$ for $j=0,1,2 $. 
$f(|D|)$ is then a good operator by virtue of Mikhlin's theorem 
(cf. \cite{Stein}). $f$ is a {\it good multiplier for small $\lam>0$} if  
$\chi_{\leq a}(\lam) f(\lam)$ is a good multiplier for a $a>0$. 
$\Mg(\R^2)$ is the space of good multipliers. We denote 
$\|f\|_{\Mg,p}=\|f(|D|)\|_{\Bb(L^p)}$, $1<p<\infty$.

\item[{\rm (3)}] 
Let $k=0, 1, \dots$ and $h(\lam)>0$. 
$T(\lam)\in \Og^{(k)}_{\Lg_1}(h)$ as $\lam\to 0$ 
(or $\lam\to \infty$) if $T(\lam,x,y)$ is a function of $\lam>0$ 
of $C^k$-class for a.e. $(x,y)$ 
and simultaneously as an $\Lg_1$-valued function and 
$\|\pa_\lam^j T(\lam) \|_{\Lg_1} 
\leq C |h(\lam)| \lam^{-j}$ for $0\leq j \leq k$ 
as $\lam\to 0$ (or $\lam\to \infty$).

\item[{\rm (4)}] $T(\lam)\in \Og_{2}(h)$ if $T(\lam,x,y)$ 
satisfies the properties of {\rm (3)} with $k=2$ and 
$L^2(\R^2 \times \R^2)$ and $\Hg_2$  
replacing  $L^1(\R^2 \times \R^2)$ and $\Lg_1$ respectively. 
If  $T(\lam)\in \Og_{2}(h)$ and $v, w \in \HL$, then 
$v T(\lam) v \in   \Og^{(2)}_{\Lg_1}(h)$. 
\end{itemize}
We often denote by $\Og^{(k)}_{\Lg_1}(h)$ etc. an operator in 
the class $\Og^{(k)}_{\Lg_1}(h)$ etc.  
\eddf 

\bgdf  The operators defined by \refeq(sta-0) with $T$ or $T(\lam)$ 
replacing $vM(\lam)^{-1}v$ are denoted respectively by $W(T)$ and 
$\Wg(T(\lam))$: 
\begin{align}
& W(T)u (x)= \int_0^\infty 
(G_0(-\lam)T \Pi({\lam})u)(x) \lam d\lam, \ u \in \Dg_{\ast}\,, \lbeq(WTdef) \\
& \Wg(T(\lam))u (x)= \int_0^\infty 
(G_0(-\lam)T(\lam) \Pi({\lam})u)(x) \lam d\lam, \ u \in \Dg_{\ast}\,.  
\lbeq(WgTdef)
\end{align}
We say $T$ or $T(\lam)$ is a {\it good producer} if $W(T)$ 
(resp.$\Wg(T(\lam))$) extends to a good operator and it is a 
{\it good producer for small $\lam>0$ or large $\lam$} if 
$\chi_{\leq 2a}(\lam) T$ or $\chi_{>2a}(\lam) T$ 
(resp. $\chi_{\leq 2a}(\lam)T(\lam)$ or $\chi_{>2a}(\lam)T(\lam)$) 
is a good producer for some $a>0$. 
\eddf

The one dimensional operator $f \mapsto v(x) \int_{\R^2} w(y)f(y) dy$ 
{\it without complex conjugate} will be indiscriminately denoted 
by $v \otimes w$ or $|v\ra \la w|$.  

\bgprop\lbprop(funda) Let $f\in \Mg(\R^2)$, $F \in L^1(\R^2)$ and 
$T\in \Lg_1$. Suppose $a>0$, $\ep>0$ and $k\in \N$. 
Then, we have the following statements for all $1<p<\infty$: 
%The followings are good producers:   
\ben 
\item[{\rm (1)}] $\|W(M_F)u\|_p \leq C\|F\|_1 \|u\|_p$\,. 
% If Multiplication operator $M_F$ with  $F(x)\in L^1(\R^2)$.
\item[{\rm (2)}] $\|\Wg(f(\lam)M_F)u\|_p 
\leq C\|f\|_{\Mg,p}\|F\|_1 \|u\|_p$\,.

\item[{\rm (3)}] $\|W(T)u\|_p \leq C\|T\|_{\Lg_1} \|u\|_p$\,.

\item[{\rm (4)}] $\|\Wg(f(\lam)T)u\|_p 
\leq C\|f\|_{\Mg,p}\|T\|_{\Lg_1} \|u\|_p$\,. 

\item[{\rm (5)}] $\chi_{\leq 2a}(\lam)T(\lam)$ is a good producer if  
$T(\lam) \in \Og_{\Lg_1}^{(2)}(\lam^{1+\ep})$ as $\lam \to 0$. 

\item[{\rm (6)}] $\chi_{>2a}(\lam)T(\lam)$ is a good producer if 
$T(\lam)\in \Og_{\Lg_1}^{(2)}(\lam^{-\ep})$ as $\lam\to \infty$. 

\item[{\rm (7)}] Let 
$I_{k,a}^{\p,\ph}(\lam) 
=f(\lam)\chi_{\leq 2a}(\lam)(\log \lam)^k (\p\otimes \ph)$ for 
$\ph,\p\in L^1(\R^2)$\,.  Assume further that  
either $\ax\ph\in L^1(\R^2)$ and $\int_{\R^2}\ph dx=0$ or 
$\ax\p \in L^1(\R^2)$ and $\int_{\R^2}\p dx=0$. Then, 
$I^{\p,\ph}_{k,a}(\lam)$ 
is a good producer and we respectively have 
\begin{gather}
\|\Wg(I_{k,a}^{\p,\ph}(\lam))u\|_p \leq C_p \|f\|_{\Mg,p}
\|\ax\ph\|_1 \|\p\|_1 \|u\|_p, \lbeq(a) \\
\|\Wg(I_{k,a}^{\p,\ph}(\lam))u\|_p \leq C_p \|f\|_{\Mg,p}
\|\ph\|_1 \|\ax\p\|_1 \|u\|_p\,. \lbeq(b)
\end{gather}
\een
\edprop

In section 4 we prove  
that $W_{+}\chi_{> 2a}(|D|)$ is a good operator if $a>0$ 
under the condition that $\ax^2 V \in L^{\frac43}(\R^2)$. 
If we expand $M(\lam)^{-1}$ in \refeq(sta-0-high) as 
\bqn \lbeq(mlam)
M(\lam)^{-1}= \sum_{j=0}^4 (-1)^j U(vG_0(\lam)w)^j 
- U(vG_0(\lam)w)^5(1+ vG_0(\lam)w)^{-1} ,
\eqn
$\W_{{\rm high},2a}$ becomes the sum $\sum_{j=0}^5 (-1)^j \W_{h,j}$ 
and we prove $\W_{h,j}$, $0\leq j \leq 5$ are good operators separately. 
$\W_{h,0}=Z_V$ is a good operator by virtue of \refprop(funda) (1). 
For proving the same for $\W_{h,1}$ we represent it as 
\bqn \lbeq(Wh1) 
\W_{h,1}u= \int_{\R^2} W(M_{V^{(2)}_y}) (\Hg(|y||D|)\tau_y \chi_{> 2a}(|D|)u) dy. 
\eqn 
where we used the short hand notation 
\[
\Hg(\lam) = (i/4) H_0^{(1)}(\lam) \ \mbox{and} \ V^{(2)}_y(x)= V(x)V(x-y).
\]
Since $\Hg(\lam) = e^{i\lam}\w(\lam)$ with symbol $\w(\lam)$ of order $-1/2$ 
(see \refeq(large-lam)), the theory of the spatially homogeous Fourier 
integral operators (\cite{Peral,Tao}) implies  
\[
\|\Hg(|y||D|)\chi_{> 2a}(|D|)\|_{\Bb(L^p(\R^2))}\leq C_p(1+|\log |y||)
\]
(see \reflm(H02y)). Then, \refprop(funda) (1) implies 
\[
\|\W_{h,1}\|_{\Bb(L^p)} 
\leq C_p \int_{\R^2}|V(x)V(x-y)|(1+|\log |y||)dxdy <\infty, 
\]
and $\W_{h,1}$ is a good operator. For obtaning the expression \refeq(Wh1), 
it is important to observe that $VG_0(\lam)V$ is the superposition of 
product of mulitiplication, multiplier and translation: 
\bqn 
VG_0(\lam)Vu(x)= 
\int_{\R^2} M_{V^{(2)}_y}\Hg(|y|\lam) (\tau_y u)(x) dy\,. \lbeq(wG0w) 
\eqn 
Then, \refeq(Wh1) follows by substituting \refeq(wG0w) for 
$vM(\lam)^{-1}v$ in \refeq(sta-0-high) and by applying \refeq(mult) 
which implies $\Hg(|y|\lam)\Pi(\lam)= \Pi(\lam)\Hg(|y||D|)$.  
We show in \refprop(n-th-estimate) that $\W_{h,j}$, $j=2,3,4$ 
are represented likewise (see \refeq(Wn)) and that they are also 
good operators. We prove that 
$\W_{h,5}$ is a good operator by showing that 
$\|\pa_{\lam}^j vG_0(\lam)w\|_{\Hg_2}\leq C \lam^{-1/2}$ as 
$\lam \to \infty$ for $j=0,1,2$ (see \reflm(AB-dec)) 
and, if sandwitched by $v$, the final term on the right of \refeq(mlam)  
is of class $\Og_{\Lg_1}^{(2)}(\lam^{-\frac12})$.

We study $\W_{{\rm low},2a}$ in section 5, which will be divided into 
seven subsections. We need study $M(\lam)^{-1}$ near $\lam=0$ 
in detail. Define three integral operators $N_0$, $G_1$ and $G_2$ by 
\begin{gather}
N_0u(x) =-\frac1{2\pi} \int_{\R^2}\log |x-y|u(y) dy\,, \\
G_1 u(x) = \frac1{4}\int_{\R^2} |x-y|^2 u(y) dy\,, \\
G_2 u(x) = \frac1{8\pi} \int_{\R^2} |x-y|^2 
\log\left(\frac{e}{|x-y|}\right) u(y) dy \,.
\end{gather}
Then, the expansion \refeq(hankel) of the Hankel function $\Hg(\lam)$ 
with the definition \refeq(g-def) of $g(\lam)$ 
implies that, as $\lam \to 0$, $M(\lam)$ is equal to  
\bqn  \lbeq(m-exp)
 U + g(\lam)v \otimes v + v N_0 v %\\ 
+ \lam^2 g(\lam)vG_1 v + \lam^2 v G_2 v + \Og_2(g\lam^4).
\eqn  
We define projections $P$, $Q$ and the operator $T_0$ by  
\[
P=({v}/\|v\|_2)\otimes (v/\|v\|_2), \ \ Q=1-P, \ \ 
T_0 = U + v N_0 v.
\]
Following definition is due to \cite{JN} (see also \cite{Schlag,EG, EGG}): 

\bgdf \lbdf(singularities) We say $H$ is regular at 
zero if $QT_ 0Q\vert_{Q\HL} $ is invertible. Otherwise, 
$H$ is singular at zero. If $H$ is singular at zero, let  
$S_1$ be the projection in $Q\HL$ onto $\Ker_{Q\HL}QT_ 0Q$. 
We say:  
\ben 
\item[{\rm (1)}] $H$ has singularities of the first kind at zero  
if $T_1= S_1 QT_0 P T_0 QS_1\vert_{S_1 L^2}$ is non-singular. 
If $T_1$ is singular, let $S_2$  
be the projection in $S_1 \HL$  onto $\Ker T_1 $.
\item[{\rm (2)}]  
$H$ has singularities of the second kind at zero  
if $T_2= S_2(v G_1 v)S_2\vert_{S_2 L^2}$ is non-singular. 
\item[{\rm (3)}] $H$  has singularities of the third kind at zero 
if $T_2$ is singular. Let $S_3$ be the projection in 
$S_2\HL$ onto $\Ker T_2 $.
\een 
\eddf 

In subsection 5.1, we recall from \cite{JN} 
the relation between the kind of singularities of $H$ at zero 
and the existence/absence of specific kinds of resonances. 
We shall give a brief proof for readers' convenience and add 
some remarks. In subsection 5.2 we recall without proof 
the Feshbach formula, Jensen-Nenciu's lemma (\cite{JN}) and  
two estimates from \cite{EG,EGG} on 
\[
M_0(\lam) = M(\lam) - g(\lam)v\otimes v - T_0, \ 
M_1(\lam) = M_0(\lam)- v\lam^2 (g G_1 + G_2)v .
\]

In subsection 5.3 we prove that $\W_{{\rm low},2a}$ is a good operator 
if $H$ is regular at zero by showing that 
\bqn 
M(\lam)^{-1}= (g(\lam) + c)^{-1}L+ \Bg + \Og_2(g \lam^2) \quad (\lam \to 0)
\eqn 
and applying \refprop(funda), where $c$ is a constant and $\Bg$ is the sum of 
the multiplication by a bounded function and a Hilbert-Schmidt operator. 
This is a considerably 
simpler proof than the one in \cite{Y-2dim}.

In subsection 5.4 we present a few identities which are necessary 
for studying $M(\lam)^{-1}$ near $\lam=0$. 
It will be important to observe that operators which appear as the coefficients 
of the singularities of $M(\lam)^{-1}$ at $\lam=0$ are those on 
the finite dimensional subspace $S_1\HL$ and that all $\z\in S_1\HL$ satisfy 
\bqn \lbeq(moment)
\int_{\R^2} v(x) \z(x)=0, 
\eqn 
which will cancel some singularities of $M(\lam)^{-1}$ at $\lam=0$.  

In section 5.5 we give a new proof of Erdo{\u g}an, Goldberg and Green's 
result (\cite{EGG}) that $\W_{{\rm low},2a}$ is a good operator 
if $H$ has singularities of the first kind at zero
under the condition $\ax^{4+\ep}V \in L^1(\R^2)$ for an $\ep>0$. 
In this case $\rank S_1= 1$ 
and $S_1 = \z\otimes \z$ for a normalized $\z\in S_1\HL$. 
We then show that 
$vM(\lam)^{-1}v\equiv - c_1\log \lam (v\z\otimes \z{v})$ 
modulo a good producer for small $\lam>0$ and hence, is itself 
a good producer for small $\lam>0$.  

In subsection 5.6 and 5.7 we asume $\ax^\c V \in L^1(\R^2)$ for $\c>8$. 
We study the case that $H$ has singularities of the second kind in 
subsection 5.6. 
Then, $1\leq \rank T_2= \rank S_2 \leq 2$ and we assume $\rank S_2=2$ 
as the case  $\rank S_2=1$ is easier. 
Let $\{\z_1, \z_2\}$ be the orthonormal basis of $S_2L^2(\R^2)$ 
such that $T_2\z_j = -\k_j^2 \z_j$, $\k_j>0$, $j=1,2$ and 
$\tRg_1(\lam)= S_2 v(G_1+ g(\lam)^{-1}G_2)v S_2$. Then, 
the representation matrix $C(\xi)=(c_{jk}(\lam))$ 
for $\tRg_1(\lam)$ on $S_2 \HL$ with respect to the basis 
$\{\z_1, \z_2\}$ satisfies 
\[
c_{jk}(\lam) = -\k_j^2 \d_{jk}+ \Og_2(g(\lam)^{-1}), \quad j,k=1,2.
\] 
Write $D(\lam)=(d_{jk}(\lam))$ for $C(\lam)^{-1}$. $d_{jk}(\lam)$ are 
good multiplers 
for small $\lam>0$. Via the threshold analysis, we first prove that modulo 
a good operator $\W_{{\rm low},2a}u(x)$ is equal to 
\bqn 
\lbeq(tWlow-a)
- \sum_{j,k=1}^2\int_0^\infty g(\lam)^{-1}\lam^{-2} d_{jk}(\lam)
(G_0(-\lam)v\z_j)(x) \la \z_k v, \Pi(\lam)u\ra {\lam}\chi_{\leq 2a}(\lam) d\lam.
\eqn 
Recalling \refeq(moment) that $\la \z_k v, 1\ra=0$, we replace 
$\Pi(\lam)u(z)$ by 
\[
\Pi(\lam)u(z)-\Pi(\lam) u(0)= \frac1{2\pi} \int_{{\mathbb S}^1} 
(e^{i\lam z\w{x}}-1) \hat{u}(\lam\w)d\w 
\]
which we decompose into the sum of 
{\it good part} $\tilde{g}(\lam,z)$ and {\it bad part} $\tilde{b}(\lam,z)$ 
as follows by Taylor expanding $e^{i\lam z\w{x}}$ upto the second order:
\begin{align}
\tilde{g}(\lam,z)& = \frac{-\lam^2}{2\pi}\int_{{\mathbb S}^1}
\left( \int_0^1 (1-\th) ({z}\w)^2 
e^{i\lam{z}\w\th}d\th\right)\hat{u}(\lam\w)d\w \lbeq(good-def-intro)  \\
& = \sum_{j,k=1}^2 z_j z_k 
\lam^2 \int_0^1 (1-\th) \left(\frac{-1}{2\pi}
\int_{{\mathbb S}^1} \Fg(\tau_{-\th{z}}R_j R_k u)(\lam\w)d\w \right) 
d\th. \notag  \\
\tilde{b}(\lam,z)& = \frac{i\lam}{2\pi}\int_{{\mathbb S}^1} 
({z}\w)\hat{u}(\lam{\w})d\w 
=\frac{i\lam}{2\pi}\sum_{l=1}^2 z_l 
\int_{{\mathbb S}^1} \Fg(R_l u)(\lam{\w})d\w.
\lbeq(bad-def-intro)
\end{align} 
Let $\tilde{\W}_{(g)}$ and $\tilde{\W}_{(b)}$ 
be respectively defined by \refeq(tWlow-a) by replacing $\Pi(\lam)u(z)$ 
by $\tilde{g}(\lam,z)$ and $\tilde{b}(\lam,z)$. 
Then, $\tilde{\W}_{(g)}$ is a good operator because 
the singularity $\lam^{-2}$ is cancelled by $\tilde{g}(\lam,z)$,
$\m_{jk}(\lam){=} 
\tilde{g}(\lam)^{-1}d_{jk}(\lam)\chi_{\leq 2a}(\lam) \in \Mg(\R^2)$ and 
$\|\ay^2 v\z\|_1 \leq C \|\ay^{4+\ep}V\|_1$ for $p$-wave resonances. 
The bad part $\tilde{b}(\lam,z)$ has only the factor 
$\lam$ and we show $\tilde{\W}_{(b)}$ is bounded in $L^p(\R^2)$ 
only for $1<p\leq 2$ and is unbounded for $2<p<\infty$. 
We avoid outlining the proof of the boundedness part for not 
making the introduction too long. 

For proving that 
$\tilde{\W}_{(b)}$ is unbounded in $L^p(\R^2)$ for $2<p<\infty$ 
it suffices to prove the same for $\chi_{> 4a}(|D|)\tilde{\W}_{(b)}u(x)$ 
which we represent as 
\[
- \sum_{j,k=1}^2\int_0^\infty \lam^{-2} \m_{jk}(\lam)
(\chi_{> 4a}(|D|)\Gg_{-\lam}\ast (v\z_j))(x) 
\la \z_kv. \tilde{b}(\lam,\cdot)\ra {\lam} d\lam\,.
\] 
Let $\m(\xi)=\chi_{>4a}(|\xi|)|\xi|^{-2}$. Then  
$\chi_{> 4a}(|D|)\Gg_{-\lam}=(2\pi)^{-1}\hat{\m}
+ \lam^2 \m(|D|)\Gg_{-\lam}$ and $\lam^2 \m(|D|)\Gg_{-\lam}$ 
produces a good operator cancelling the singularity $\lam^{-2}$. 
We show that $(2\pi)^{-1}\hat{\m}(x)$ produces an operator which 
is unbounded in $L^p(\R^2)$ for $2<p<\infty$. 
Thanks to the fact that $(2\pi)^{-1}\hat{\m}(x)$  is $\lam$-independent, 
the operator in question becomes    
\bqn \lbeq(quest)
- \sum_{j,k=1}^2  a_{j}(x) \int_0^\infty \lam^{-2} \m_{jk}(\lam)
\la \z_k v, \tilde{b}(\lam,\cdot) \ra {\lam} d\lam\,
= \sum_{j=1}^2 a_j(x) \ell_j(u)
\eqn 
where, for $j=1,2$, $a_j(x) {=}(2\pi)^{-1}(\hat{\m}\ast v\z_j)(x)\in L^p(\R^2)$ 
for $1<p<\infty$ and Parseval's identity implies that 
the linear functional $\ell_j(u)$ is equal to   
\begin{align}
\ell_j(u)& =\frac{i}{2\pi}
\sum \la z_l v, \z_k \ra 
\left(\int_0^\infty \int_{{\mathbb S}^1} \m_{jk}(\lam) 
\Fg(R_l u)(\lam\w) d\w d\lam \right) \\ \notag 
& = \frac{i}{2\pi} \sum \la z_l v, \z_k \ra  \int_{\R^2} u(x)
\Fg(\m_{jk}(|\xi|)\xi_l|\xi|^{-2})(x)dx, \lbeq(lu)
\end{align}
where the sum is taken  over $k,l=1,2$. We show  
$a_1, a_2 \in L^p(\R^2)$ are linearly independent 
if $a>0$ is sufficiently small. It follows 
by the Hahn-Banach theorem that, if $\sum_{j=1}^2 a_j(x) \ell_j(u)$ 
were bounded in $L^p(\R^2)$ for such an $a>0$, $\ell_1$ and $\ell_2$ 
must be bounded on $L^p(\R^2)$ hence, 
by virtue of the Riesz theorem, it must be that for $q=p/(p-1)$ 
\[
\sum_{k,l=1}^2 \la z_l v | \z_k\ra
\Fg(\m_{jk}(|\xi|)\xi_l|\xi|^{-2}) \in L^q(\R^2) , \quad j=1,2.
\] 
Then, since $1<q<2$ for $2<p<\infty$, by virtue of 
Hausdorff-Young's inequality, 
\[
\sum_{k=1}^2  d_{jk}(|\xi|) 
\sum_{l=1}^2 \la z_l v | \z_k \ra \chi_{\leq 2a}(\xi)\xi_l 
g(|\xi|)^{-1}|\xi|^{-2} \in L^p(\R^2) 
\] 
and, since $C(|\xi|)=D(|\xi|)^{-1}$ is 
bounded for $|\xi|\leq 2a$, 
\[ 
\frac{\chi_{\leq 2a}(|\xi|)}
{g(|\xi|)|\xi|^{2}}\begin{pmatrix}
\la z_1 v | \z_1 \xi_1\ra + \la z_1 v | \z_2 \ra \xi_2 \\ 
\la z_2 v | \z_1 \xi_1 \ra + \la z_2 v | \z_2 \ra \xi_2 
\end{pmatrix}\in L^p(\R^2, \C^2). 
\]
But, this can happen only when $\la z_l v | \z_k\ra=0$ 
for $1\leq j,k \leq 2$, which is impossible if $T_2$ is non-singular, 
see \refeq(G1-inner) and $\sum_{j=1}^2 a_j(x) \ell_j(u)$ is unbounded 
in $L^p(\R^2)$ for $2<p<\infty$. 

In the final subsection 5.7 we assume that $H$ has singularities of 
the third kind at zero. Then, $T_3= S_3 G_2 S_3$ is necessarily 
non-singular and the Feshbach formula implies 
$(S_2\tRg_1S_2)^{-1}= g(\lam) S_3 T_3^{-1}S_3 + L_4(\lam)$, 
$L_4(\lam)\in \Mg(\R^2)$ being an operator in 
$S_2\HL$; $\z \in S_3\HL$ 
satisfies the extra moment conditions 
\bqn \lbeq(cancellation-2a)
\int_{\R^2} x_1 \z(x)v(x)dx  = \int_{\R^2} x_2 \z(x)v(x)dx  = 0. 
\eqn 
Modulo a good operator $\W_{{\rm low},2a}$ is still given by 
\refeq(tWlow-a) and, if $(d_{jk}(\lam))$ is replaced by 
the matrix for $g(\lam) S_3 T_3^{-1}S_3$ then, it produces 
a good operator. This is because  
we have  $\la \z {v}, \Pi(\lam) u \ra= \la \z {v}, \tilde{g}(\lam, \cdot)\ra$ 
for $\z \in S_3 \HL$ by virtue of \refeq(cancellation-2a) and 
the factor $\lam^2$ of $\tilde{g}(\lam,x)$ cancels 
the singularity $\lam^{-2} g(\lam)^{-1}$ leaving $g(\lam)^{-1}$ which 
then cancels $g(\lam)$ in the front of $g(\lam) S_3 T_3^{-1}S_3$.  
If $S_3=S_2$, in which case $p$-wave resonances are absent, 
then $L_4(\lam)=0$ and $\W_{{\rm low},2a}$ becomes a good operator. 
If $S_2\not=S_3$, then $L_4(\lam) \not=0$ but $L_4(\lam)$ has the structure 
similar to that of $\tRg_1(\lam)^{-1}$ of subsection 5.6, 
and we prove by slightly modifying the argument in subsection 5.6 
that it produces an operator which is bounded in $L^p$ for $1<p\leq 2$ 
and unbounded for $2<p<\infty$. 
In the rest of the paper we shall give the details of the proof. 

\section{Preparatory estimates}

We often write $\lam$ for $z\in \Cb^{+}$ when we want to emphasize that 
$z$ can also be real. $a\absleq b$ means that $|a|\leq |b|$. 
Recall that $G_0(\lam)u(x)= (\Gg_{\lam}\ast u)(x)$, 
$\Gg_{\lam}(x) = \Hg(\lam|x|)$ and $\Hg(\lam)=(i/4)H_0^{(1)}(\lam)$. 
The Hankel function $\Hg(\lam)$ has the well known 
series expansion (\cite{AS}, p.358, (9.1.12) and (9.1.13)) and the 
integral representation (\cite{AS},p. 360, (9.1.23)):
\begin{align} 
& \Hg(\lam )= g(\lam) 
\sum_{k=0}^\infty \frac{(-1)^k}{(k!)^2}\Big(\frac{\lam ^2}{4}\Big)^k 
 \notag \\
& +\frac1{2\pi}\left(\frac{\frac14 \lam ^2}{(1!)^2}
-\Big(1+\frac12\Big)\frac{(\frac14 \lam ^2)^2}{(2!)^2} 
+ \Big(1+\frac12+\frac13\Big)
\frac{(\frac14 \lam ^2)^3}{(3!)^2} - \cdots \right) \lbeq(hankel) \\
& = \frac{e^{i\lam}}{2^{\frac32}\pi} 
\int_0^\infty e^{-t}t^{-\frac12}\left(\frac{t}2-i\lam\right)^{-\frac12}dt 
\lbeq(Hank) 
\end{align}
where $z^\frac12$ is positive for positive $z$ and, 
with Euler's constant $\c$ and with principal branch of $\log z$ 
\bqn \lbeq(g-def)
g(\lam)= -\frac1{2\pi}\log\Big(\frac{\lam }{2}\Big)+ \frac{i}4 -
\frac{\gamma}{2\pi}\,.
\end{equation}

The following lemma is obvious 
from \refeq(hankel) and \refeq(Hank).

\bglm \lblm(hankel)
\ben
\item[{\rm (1)}] $\Hg(\lam)$ and $\Gg_\lam(x)$ satisfy for any $\d>0$ 
\begin{align} 
& \Hg(\lam)  = g(\lam)  + \lam^2\Big(-\frac{g(\lam)}{4} 
+ \frac1{8\pi}\Big) + O(g(\lam) \lam^4), \lam \to 0 , \lbeq(small-lam) \\
& \Gg_\lam(x)= g(\lam) + N_0(x) + O((\lam|x|)^{2-\d}), \ \ \lam|x| \to 0. 
 \lbeq(hankel-est-separate) 
\end{align}
\item[{\rm (2)}] For $\lam\geq 1$, 
$\Hg(\lam) = e^{i\lam}\w(\lam)$ with  symbol $\w(\lam)$ of order $-1/2$:   
\bqn 
|\w^{(j)}(\lam)|\leq C_j \lam^{-\frac12-j}, \quad \lam \geq 1, \ 
j=0,1,\dots.  \lbeq(large-lam) 
\eqn 
\item[{\rm (3)}] There exists a constant $C>0$ such that for 
\bqn 
|\Gg_\lam(x)|\leq  C\left\{\br{ll} \la \log \lam |x| \ra,  & \quad 
|\lam |x| |\leq 1, \\ 
\la \lam|x| \ra^{-1/2}, & \quad |\lam |x| |\geq 1 \,. \er
\right.  
\lbeq(hankel-est) 
\eqn 
\een
\edlm 

\bglm  \lblm(hankel-a)
Let $0<\a<\b<\infty$. Then, for a constant $C_{\a,\b}$,
\bqn \lbeq(hankel-int)
\int_\a^\b |\Gg_\lam(x)|d\lam \leq C_{\a,\b}(|x|^{-\frac12}+ |x|^{-1}), 
\quad x \in \R^2\,.
\eqn 
\edlm 
\bgpf Applying \refeq(hankel-est), we estimate the integral by 
\bqn \lbeq(hankel-int-1)
C \int_{\a <\lam <|x|^{-1}} (|\log \lam|x||+ 1)d\lam + 
C \int_{|x|^{-1}<\lam <\b} (\lam |x|)^{-1/2} d\lam .
\eqn 
For $|x|<\b^{-1}$ the second integral vanishes and via change of variable 
\[
\int_\a^\b |\Gg_\lam(x)|d\lam 
\leq \frac{C}{|x|}\int_{\a|x|}^1 (-\log s + 1)ds
\leq \frac{C}{|x|}\int_{0}^1 (-\log s + 1)ds =\frac{C}{|x|}.
\]
When $|x|>\a^{-1}$ the first integral vanishes and  
\[
\int_\a^\b |\Gg_\lam(x)|d\lam  \leq 
\frac{2C}{|x|^{\frac12}}(\b^\frac12 - |x|^{-\frac12})
\leq \frac{2C\b^\frac12}{|x|^\frac12}\, .
\]
Since the left side is continuous for $\b^{-1}\leq |x|\leq \a^{-1}$, 
\refeq(hankel-int) follows. 
\edpf 

\bglm \lblm(first-step) Let $u \in \Dg_\ast$. Then, 
$\Pi(\lam)u(x)\in C^\infty(\R_\lam \times \R^2_x)$ and 
there exist constants $0<\a<\b<\infty$ such that 
\begin{gather} 
\lbeq(support)
{\rm supp}_\lam \Pi(\lam)u(x) \subset (\a,\b).  \\
|\Pi(\lam)u(x)|\leq \min(\|u\|_1, \ C_u \ax^{-1/2}). \lbeq(Pi-est)
\end{gather} 
\edlm
\bgpf The definition of $\Dg_\ast$ trivially implies \refeq(support).
The estimate \refeq(Pi-est) is well known, see e.g. \cite{Stein}, p.348.
\edpf

The following lemma is proved in \cite{CMY}. 
$K$ is defined by \refeq(DFEK). 

\begin{lemma} 
\lblm(1) 
$K$ satisfies the identity 
\refeq(how-0) for $x,y,z\in \R^2$. 
Moreover:  
\ben 
\item[{\rm (1)}] 
$Ku(x)$ is a rotationary invariant and  
\bqn \lbeq(K-sing)
Ku(x) = \lim_{\ep \downarrow 0} 
\frac{-1}{(2\pi)^2}\int_{\R^2} \frac{u(y)dy}{x^2- y^2-i\ep}.  
\eqn 
\item[{\rm (2)}]
$K$ is bounded in $L^p(\R^2)$ for any $1<p<\infty$:
\bqn \lbeq(K-est)
\| Ku \|_p \leq C_p \|u\|_p\,. 
\eqn 
\een
\end{lemma}
\bgpf We give a simpler proof of \refeq(K-sing). 
Since $d\eta= \lam d\lam d\w$ in the polar coordinates 
$\eta= \lam \w$, we have for $u, v \in \Dg_\ast$ that 
\begin{align*}
& \la Ku, v\ra = \frac{1}{2\pi}
\int_0^{+\infty} \lam \left(
\int_{{\mathbb S}^1} \hat{u}(\lam \w) d\w 
\right)
\la \mathcal{G}_{-\lambda}, {v}\ra d\lambda\, \notag \\
&= \lim_{\ep \downarrow 0}\frac{1}{(2\pi)^2}
\int_0^{+\infty} 
\lam \left(\int_{{\mathbb S}^1} \hat{u}(\lam \w) d\w \right) 
\left(\int_{\R^2}\frac{\hat{v}(\xi)}{\xi^2-\lam^2+i\ep} d\xi\right) 
d\lambda\, .\notag \\
& =  
\lim_{\ep \downarrow 0}\frac{1}{(2\pi)^2}
\iint_{\R^4} 
\frac{\hat{u}(\eta)\hat{v}(\xi)}{\xi^2-\eta^2+i\ep} 
d\xi d\eta \notag \\
& = \lim_{\ep \downarrow 0} \frac{-i}{2}\int_0^\infty 
e^{-t\ep} 
\left(\frac1{2\pi}\int_{\R^2}    
e^{-it\eta^2/2}\hat{u}(\eta)d\eta\right)\left(\frac1{2\pi}\int_{\R^2}
e^{it\xi^2/2}\hat{v}(\xi)d\xi \right) dt\, .
\end{align*}
Since $u,v \in \Dg_\ast$, integration by parts shows that both 
\[
\frac1{2\pi}\int_{\R^2}    
e^{-it\eta^2/2}\hat{u}(\eta)d\eta \ \ \mbox{and}\ \ 
\frac1{2\pi}\int_{\R^2}
e^{it\xi^2/2}\hat{v}(\xi)d\xi
\]
are rapidly decreasing as $t \to \infty$  and 
they are equal respectively to 
\[
\frac{1}{2t\pi{i}}\int_{\R^2}
e^{{i}x^2/2t}u(x) dx \ \ \mbox{and}\ \ \frac{i}{2t\pi}\int_{\R^2}
e^{-{i}x^2/2t}u(x) dx, \quad t>0 .
\]
Thus, the dominated convergence theorem implies that 
\begin{align*}
\la Ku, v\ra& =\lim_{\ep \downarrow 0} \frac{-i}{2(2\pi)^2}\int_0^\infty 
e^{-\ep/2t}\left(\iint_{\R^4} e^{i(x^2-y^2)/2t}u(x)v(y) dxdy \right)
\frac{dt}{t^2} \\
& = \lim_{\ep \downarrow 0} \frac{-i}{(2\pi)^2}\int_0^\infty 
\left(\iint_{\R^4} e^{is(x^2-y^2+i\ep)}u(x)v(y) dxdy \right){ds} \\
& = \lim_{\ep \downarrow 0} \frac{1}{(2\pi)^2}
\iint_{\R^4} \frac{u(x)v(y)}{x^2-y^2+i\ep} dxdy .
\end{align*}
This proves \refeq(K-sing). Denote $Ku(|x|)=Ku(x)$. Then, 
\[
Ku(\sqrt{r}) = \lim_{\ep \downarrow 0} 
\frac{-1}{4\pi}
\int_{0}^\infty 
\frac{1}{r- \r-i\ep}\left(\frac{1}{2\pi}\int_{{\mathbb S}^1}
u(\sqrt{\r}\w)d\w \right) d\r. 
\]
This is essentially the Hilbert transform and  
\begin{align*}
\|Ku\|_p^p & = \pi \int_0^\infty |Ku(\sqrt{r})|^p dr 
\leq C \int_0^\infty 
\left|\frac{1}{2\pi}\int_{{\mathbb S}^1}
u(\sqrt{\r}\w)d\w \right|^p d\r \\
& \leq \frac{C}{2\pi} \int_0^\infty 
\int_{{\mathbb S}^1}
|u(\sqrt{\r}\w)|^p d\w  d\r = \frac{C}{\pi} \|u\|_p^p,  
\end{align*}
where H\"older's inequality is used in the third step. 
\edpf

\section{Integral estimates}

In this section, we prove \refprop(funda). We often identify integral operators 
with their integral kernels and say e.g. the integral operator $T(x,y)$. 
We always suppose $u \in \Dg_\ast$ which will not be mentioned anymore. 
In this section we assume $V \in L^1(\R^2)$ if otherwise stated explicitly

\begin{lemma} For $a, b \in \R^2$, define    
\bqn \lbeq(Wab)
\W_{a,b}u(x)= \int_0^\infty \Gg_{-\lam}(x-a) (\Pi(\lam)u)(b)
\lam d\lam.
\eqn 
Then, for a constant $C_p$ independent of $a,b$ and $u$, 
\bqn \lbeq(a-b) 
\|\W_{a,b} u \|_p \leq C_p \|u\|_p.
\eqn 
\edlm 
\bgpf By virtue of \reflms(hankel-a,first-step), the integral is 
absolutely convergent unless $x=a$ and, \refeq(how-0) implies 
$\W_{a,b} u(x) =(\tau_a K \tau_{-b}{u})(x)$. 
Since translations are isometries of $L^p$,  
\refeq(K-est) implies \refeq(a-b).
\edpf 

\begin{lemma} \lblm(multiplier) 
Let $F \in L^1(\R^2)$ and $f \in \Mg(\R^2)$. Then: 
\bqn \lbeq(prod-fM)
\|\Wg(f(\lam) M_F)u\|_p \leq C \|F\|_1 \|f\|_{\Mg,p} \|u\|_p, \quad 1<p<\infty. 
\eqn 
In particular $M_F$ is good producer.
\edlm
\bgpf We may assume $f(\lam)=1$ since 
$f(\lam) \Pi(\lam)u= \Pi(\lam)f(|D|)u$ by virtue of \refeq(mult) 
and $\|f(|D|)u\|_p \leq \|f\|_{\Mg,p} \|u\|_p$. We have 
\bqn \lbeq(first-1)
W(M_F)u(x) 
= \int^{\infty}_0 \left(
\int_{\R^2} \Gg_{-\lam}(x-y)F(y)\Pi(\lam)u(y) dy \right)\lam d\lam.
\eqn 
If we integrate with respect to $\lam$ first, \refeq(Wab) implies  
\bqn \lbeq(first-3)
W(M_F)u(x) = 
\int_{\R^2} F(y)(\W_{y,y}u)(x) dy, \quad {\rm a.e.}\ x \in \R^2
\eqn 
and Minkowski's inequality and \refeq(a-b) yield the desired estimate: 
\[
\|W(M_F)u\|_p 
\leq \int_{\R^2} |F(y)|\|\W_{y,y}u(x)\|_p dy 
\leq C_p \|F\|_1 \|u\|_p. 
\]
To change the order of integrations we argue as follows. 
The argument is standard and will be repeatedly used in what follows. 
Let $B_R=\{x\colon |x|\leq R\}$. Since $u \in \Dg_\ast$, 
$\Pi(\lam)u(z)=0$ for $\lam \not\in (\a,\b)$ 
and $|\Pi(\lam)u(z)|\leq C $ by virtue of \reflm(first-step). 
It follows by using \refeq(hankel-int) that 
\begin{multline} \lbeq(1-23)
\int_{B_R} \left(\int_0^\infty 
|\Gg_{-\lam}(x-y)F(y)(\Pi(\lam)u)(y)|\lam{d\lam}\right) dx  
\\
\leq C |F(y)|\int_{B_R} (|x-y|^{-1/2}+ |x-y|^{-1})dx 
\leq C |F(y)|\,.
\end{multline}
and the following $5$-dimensional integral is (absolutely) convergent: 
\[
\int_{\R^2 \times B_R \times [0,\infty)} 
|\Gg_{-\lam}(x-y)F(y)(\Pi({\lam})u)(y) \lam | d\lam dx dy <\infty.
\]
Then, Fubini's theorem implies that for a.e. $x\in \R^2$  
\refeq(first-1) is integrable on $[0,\infty) \times \R^2$ 
and the order of the integrals can be changed. 
This concludes the proof. 
\edpf

\begin{lemma} \lblm(rank1)  
Let $T \in \Lg_1$ and $f \in \Mg(\R^2)$. Then 
for a constant $C_p>0$ 
\bqn 
\|\Wg(f(\lam) T)u\|_p 
\leq C_p \|T\|_1 \|f\|_{\Mg,p} \|u\|_p\,, \quad 1<p<\infty.  \lbeq(rank-1)
\eqn  
In particular $W(T)$ is a good operator 
\edlm 
\bgpf It suffices to show \refeq(rank-1) for $W(T)u$ as previously. 
If the integral on the right side is integrable for a.e. 
$x\in \R^2$, we may integrate with respect to $\lam$ first. 
Then, \refeq(how-0) and \refeq(v-w) imply
\begin{align} \lbeq(v-w)
& W(T)u(x)= \int_0^\infty 
\iint_{\R^2\times \R^2}\Gg_{-\lam}(x-y)T(y,z)\Pi(\lam)u(z) dz dy d\lam \\ 
& \qquad = \int_{\R^2\times \R^2 } T(y,z) 
\left\{\int_0^\infty \Gg_{-\lam}(x-y) \Pi(\lam)u(z) \lam d\lam \right\} 
dy dz \notag \\
& \qquad = \int_{\R^2\times \R^2 } T(y,z)(\t_y K \t_{-z}{u})(x) dy dz\,.
\lbeq(wGgl-a)
\end{align}
Then, Minkowski's inequality and \reflm(1) imply 
\[
\|W(T)u\|_p \leq \int_{\R^2\times \R^2 } |T(y,z)|
\|K\|_{\Bb(L^p)}\|u\|_p dy dz \leq C \|u\|_p. 
\]
To see that \refeq(v-w) is integrable for a.e. $x\in \R^2$, we repeat 
the argument of the proof of the previous lemma. 
Applying \refeq(support), 
\refeq(Pi-est) and \refeq(hankel-int) and using the notation of 
there, we estimate as in \refeq(1-23) 
\begin{align}
& \int_{\R^2\times \R^2}
|T(y,z)| 
\left(\int_{B_R \times [0,\infty)} 
\lam |\Gg_{-\lam}(x-y)| |\Pi(\lam)u(z)|dx d\lam\right)
dy dz   \notag \\
& \leq 
C \int_{\R^2\times \R^2} \az^{-\frac12} |T(y,z)|dy dz 
\leq C \|T\|_1 <\infty.   \lbeq(Wuv-est)
\end{align}
This completes the proof. \edpf 

\bglm \lblm(W-T-abs-0) Let $\s>1$ and $a>0$. 
Suppose $T(\lam)\in  \Og_{\Lg_1}^{(2)}(\lam^\s)$ as 
$\lam \to 0$.  Then, for a constant $C_p>0$ 
\bqn  \lbeq(W-T-abs) 
\|\Wg(\chi_{\leq a}(\lam)T(\lam)) u \|_p 
\leq C_p \|u\|_p, \quad 1<p<\infty.
\eqn  
\edlm 
\bgpf Since $T(0,x,y)=T'(0,x,y)=0$ for a.e. $(x,y)\in \R^2 \times \R^2$, 
integration by parts shows as a Riemann integral in $L^1(\R^4)$ that 
\bqn \lbeq(tlam)
T(\lam,x,y)= \int^\infty_0 ({\lam}-{\mu})_{+} T''(\m,x,y)d\m, 
\eqn 
and also pointwise ${\rm a.e.}$ and 
$\Wg(\chi_{\leq a}(\lam)T(\lam))u(x)$ is equal to 
\begin{multline} \lbeq(hDast)
\int_0^\infty\lam \chi_{\leq a} (\lam)  
\left[
\int_{\R^2} \Gg_{-\lam}(x-y)\left\{
\int_0^{\infty}(\lam-\mu)_{+}  \right. \right. \\
\left. \left. \times 
\left(
\int_{\R^2} T''(\m, y,z) \Pi(\lam)u(z)dz
\right)
\chi_{\leq 2a}(\m) d\mu
\right\} 
dy 
\right]
d\lam.
\end{multline}
We have inserted ${\chi_{\leq 2a}}(\m)$ which is equal to $1$ 
if $(\lam-\mu)_{+}\chi_{\leq a}(\lam)\not=0$. 
Since $\chi_{\leq a}(\lam)\lam(\lam-\mu)_{+}\leq a^2 \chi_{\leq a}(\lam)$, 
$\Pi(\lam)u(z)=0$ for $\lam\not\in (\a,\b)$ and $|\Pi(\lam)u(z)|\leq C$, 
the argument similar to the one used 
in previous lemmas implies that \refeq(hDast) is integrable 
for a.e. $x\in \R^2$ and we can change   
the order of integration with respect to $d\lam$ and $d\mu$. We split as 
\bqn \lbeq(break-up)
({\lam}-{\mu})_{+}= {\lam}-{\mu} + (\mu-{\lam})_{+}
\eqn 
and decompose accordingly  
$\Wg(\chi_{\leq a}(\lam)T(\lam)) u= \sum_{j=1}^3 \Wg_{(j)}u$:  
\begin{align*}
\Wg_{(1)} u & =
\int_0^\infty \left(
\int_0^{\infty}G_0(-\lam)T''(\mu) 
\Pi(\lam)|D|\chi_{\leq a}(|D|) {u} \lam d\lam
\right) {\chi}_{\leq 2a}(\m) d\mu \\
\Wg_{(2)} u & =
- \int_0^\infty \left(\int_0^{\infty}G_0(-\lam)T''(\m)  
\Pi(\lam)\chi_{\leq a}(|D|)u \lam d\lam\right)\mu {\chi}_{\leq 2a} (\m) d\mu. \\
\Wg_{(3)} u & =\int_0^\infty 
\left(\int_0^{\infty}G_0(-\lam)T''(\m) 
\Pi(\lam) (1-|D|/\mu)_{+}
\tilde{u}\lam d\lam\right)\mu{\chi}_{\leq 2a}(\m) d\mu,
\end{align*}
where we defined $\tilde{u}= \chi_{\leq a}(|D|)u$ in the last formula.  
We apply \reflm(rank1). Then, since $\lam\chi_{\leq a}(\lam)$ and 
$\chi_{\leq a}(\lam)\in \Mg(\R^2)$, Minkowski's inequality and 
$T(\lam)\in  \Og_{\Lg_1}^{(2)}(\lam^\s)$ for $\s>1$ as $\lam \to 0$ 
%$\|T''_\m\|_{L^1(\R^4)} \leq C\m^{\s-2}$ with $\s>1$ as $\mu \to 0$ 
jointly imply that for $j=1,\, 2$  
\bqn \lbeq(Wg12)
\|{\Wg}_{(j)} u\|_p \leq C  \|u\|_p 
\int_0^{2a}\|T''(\mu)\|_{\Lg_1} \mu^{j-1} d\mu \leq C  \|u\|_p. 
\eqn 
We have $\sup_{\m>0}\|(1-|D|/\mu)_{+}\|_{\Bb(L^p)}\leq C$ for any 
$1\leq p \leq \infty$ since the Fourier transform of $(1-|\xi|/\mu)_{+}$ 
is integrable with $\mu$-indendent $L^1(\R^2)$-norm 
(see p. 426 of \cite{Stein}). 
It follows as previously that 
\bqn \lbeq(Wg3)
\|\Wg_{(3)} u\|_p  
\leq C \|u\|_p 
\int_0^{2a} \mu \|T''_(\mu)\|_{\Lg_1} d\mu 
\leq C \|u\|_p. 
\eqn 
Combining \refeq(Wg12) and \refeq(Wg3), we obtain the lemma. 
\edpf 

\bglm \lblm(W-T-abs-infi) Let 
$T(\lam) \in \Og^{(2)}_{\Lg_1}(\lam^{-\s})$ as 
$\lam \to \infty$ for a $\s>0$ and let $a>0$.  
Then, for a constant $C_p>0$ independent of $u$, 
\bqn  \lbeq(W-T-abs-geq) 
\|\Wg(\chi_{> a}(\lam)T(\lam)) u \|_p 
\leq C_p \|u\|_p, \quad 1<p<\infty.
\eqn   
\edlm 
\bgpf By integration by parts we have as previously 
\bqn \lbeq(tlamint)
T(\lam,x,y) = \int_{\lam}^\infty (\mu-\lam)_{+}T''(\m,x,y)d\m, 
\quad {\rm a.e.}\ (x,y).
\eqn 
If we use \refeq(tlamint) in place of \refeq(tlam), then   
$\Wg(\chi_{> a}(\lam)T(\lam)) u(x)$ may be expressed 
as in \refeq(hDast) 
with $(\m-\lam)_{+}$, $\chi_{> a}(\lam)$ and $\chi_{> a/2}(\mu)$ 
in place of $(\lam-\m)_{+}$, 
$\chi_{\leq a}(\lam)$ and $\chi_{\leq 2a}(\mu)$ respectively. 
Then, we repeat the argument in the proof of \reflm(W-T-abs-0). 
The argument actually is simpler here because we do not have to use 
the splitting as in \refeq(break-up). We should safely be able to 
omit the repetitious details. 
\end{proof}

For proving statement (7) of \refprop(funda), we need the following lemma. 
The lemma must be well known, however, we present a proof for readers' 
convenience. 

\bglm \lblm(logk)
Let $\k\in C_0^\infty(\R^2)$. Then for $k=0,1, \dots$. 
\bqn \lbeq(logk-Fourier)
\Fg(g(|\xi|)^{k+1} \k(\xi))(x) \absleq C_k \la \log |x| \ra^k \ax^{-2}
\eqn 
\edlm 
\bgpf It suffices to show \refeq(logk-Fourier) with $\log|\xi|$ 
replacing $g(|\xi|)$. Since $\lap \log|\xi|= -2\pi \d(\xi)$, 
$\Fg(\log|\xi|)(x) = |x|^{-2}$ in $\R^2\setminus \{0\}$. Thus, if we define 
the regularization $S(x)\in \Sg'(\R^2)$ of $|\xi|^{-2}$ by 
\[
\la S, u \ra = \int_{\R^2} \frac{u(x)-\chi(x)u(0)}{|x|^2}dx 
\]
then, $\Fg(\log|\xi|)(x)-S(x) $ is a finite sum of $C_\a D_x^\a \d(x)$ and 
\[
\Fg(\log|\xi|\k(\xi))(x)= (2\pi)(S \ast \hat\k)(x)+ \sum C_\a D_x^\a \hat{\k}(x).
\]
$D_x^\a \hat{\k}(x)$ is clearly rapidly decreasing and 
\[% begin{multline}
(S \ast \hat\k)(x)%= \int_{\R^2} \frac{\hat\k(x-y)-\chi(y)\hat\k(x)}{|y|^2}dy \\
= \left(\int_{|y|<1}+ 
\int_{|y|\geq 1}\right) \frac{\hat\k(x-y)-\chi(y)\hat\k(x)}{|y|^2}dy
= I_1(x) + I_2(x). 
\]%\end{multline*}
For $|y|\leq 1$, we have 
$\hat\k(x-y)-\chi(y)\hat\k(x)= (1-\chi(y))\hat\k(x) + 
y\nabla\hat\k(x-\th{y})$ for a $0<\th<1$ and  
\[
|I_1(x)| \leq  C |\hat\k(x)| + 
\int_{|y|<1} \frac{|\nabla \hat\k(x-\th{y})|}{|y|}dy 
\]
is rapidly decreasing and, $I_2(x)-\frac1{|x|^2}\int_{\R^2}\hat{\k}(y)dy$ 
is equal to 
\[
\frac1{|x|^2}\int_{|x-y|\leq 1}\hat{\k}(y)dy + 
\int_{|x-y|\geq 1}\hat{\k}(y)\left(\frac{1}{|x-y|^2}-\frac1{|x|^2}\right)dy 
\]
and is bounded by $C|x|^{-3}$ for large $|x|$. 
\edpf 

The idea of the proof of the following two lemmas is borrowed from 
\cite{Ya-point}. We first show that the statement (7) is satisfied 
more generally if $\ph$ on the right satisfies the vanishing moment 
condition $\la 1, \ph\ra=0$. Define for $u \in \Dg_\ast$ 
\bqn \lbeq(Wpphmu-def)
\W(\p,\ph,\mu)u(x) = 
\int_0^\infty G_0(-\lam)|\p \ra \la \ph, \Pi(\lam)u \ra 
\mu(\lam)\lam d\lam.
\eqn
Note that $\Wg(I_{k,a}^{\p,\ph})=\W(\p,\ph,g(\lam)^k \chi_{\leq 2a}(\lam))$.

\begin{lemma} \lblm(swave-r) 
Suppose $\p(x), \ax\ph(x) \in L^1(\R^2)$, 
$\int_{\R^2} \ph(x) dx=0$ and $\r(\lam)$ be such that 
$\r(\lam)\lam \in \Mg(\R^2)$. Then:
\bqn \lbeq(Wpphmu-r)
\|\W(\p,\ph,\rho)u\|_p \leq C_p \|\ax\ph\|_1 \|\p\|_1 \|u\|_p, 
\quad 1<p<\infty. 
\eqn 
In particular, \refeq(a) is satisfied for $\Wg(I_{k,a}^{\p,\ph})$. 
\edlm 
\bgpf Since $\int_{\R^2} \ph(x) dx=0$, 
$\la \ph(y), \Pi(\lam)u(y) \ra= \la \ph(y), \Pi(\lam)u(y) -\Pi(\lam)u(0)\ra$ 
and via Taylor's formula 
\begin{align*}
\Pi(\lam)u(y) -\Pi(\lam)u(0)& =\frac{1 }{2\pi}\int_0^1 \left(\int_{{\mathbb S}^1}
i\lam{y\w}e^{i\lam{y}\w\th} (\Fg u) (\lam{\w})d\w\right)d\th  
\lbeq(first-expansion-0)
\\
& =i\sum_{j=1}^2 \frac{\lam}{2\pi}\int_0^1 
\left(\int_{{\mathbb S}^1 }y_j \Fg(\t_{{-\th}y} R_j u)(\lam\w) d\w\right)
d\th,  %\lbeq(first-expansion)
\end{align*}
where $R_j= D_j/|D|$, $j=1,2$ are Riesz transforms. It follows by virtue of 
\refeq(mult) that with $\tilde{\r}(\lam) = \lam \rho(\lam) \in \Mg(\R^2)$  
\bqn \lbeq(canvce)
\rho(\lam)\la \ph, \Pi(\lam)u \ra= 
\frac{i}{2\pi} \sum_{j=1}^2 \int_{\R^2} y_j \ph(y) 
\Big( \int_0^1 \int_{{\mathbb S}^1} 
(\Fg R_j \tau_{-{\th}y}\tilde{\r}(|D|)u)(\lam\w)d\w d\th\Big)dy.   
\eqn  
Multiply \refeq(canvce) by 
$G_0(-\lam)\p(x)=  \int_{\R^2} \Gg_{-\lam}(x-z)\p(z) dz$ from the left 
and integrate with respect to $\lam d\lam$ first, 
which is regitimated as in the proof of \reflm(multiplier). 
Then, $\W(\p,\ph,\rho)u(x)$ becomes 
\[
\sum_{j=1}^2 
\int_0^1 
\left( \iint_{\R^4} y_j \ph(y)\p(z)(\tau_z K \t_{-\th{y}} 
\Rg_j \tilde{\rho}(|D|)u)(x) dy dz \right ) d\th .
\]
Then, Minkowki's inequality, Riesz theorem,  
Mikhlin's multiplier theorem and  \reflm(1)  jointly imply 
\[
\|\W(\p,\ph,\rho)u\|_p \leq C
\sum_{j=1}^2 \int_0^1 
\left( \iint_{\R^4} |y_j \ph(y)\p(z)| 
\|K\|_{\Bb(L^p)}\|u\|_p dy dz \right )d\th 
\]
and \refeq(Wpphmu-r) follows.  
\edpf 

\bglm \lblm(swave-l) 
Suppose $\ax\p(x), \ph(x) \in L^1(\R^2)$ and 
$\int_{\R^2} \p(x) dx=0$. Then, \refeq(b) is satisfied for $k=0,1, \dots$:
\bqn \lbeq(Wpphmu-l)
\|\Wg(I_{k,a}^{\p,\ph})u\|_p \leq C_{k,p} \|\ax\p\|_1 \|\ph\|_1 \|u\|_p, 
\quad 1<p<\infty. 
\eqn 
\edlm 
\bgpf We split $\Wg(I_{k,a}^{\p,\ph})u$ by multiplying it by  
$\chi_{>4a}(|D|)+ \chi_{\leq 4a}(|D|)=1$ from the left.
We deal with $\chi_{>4a}(|D|)\Wg(I_{k,a}^{\p,\ph})u$ first. 
Define $\m(\xi)= \chi_{>4a}(\xi)|\xi|^{-2}$. Then, for $\lam<2a$, 
functions of both sides of   
\[
\chi_{>4a}(|\xi|)(\xi^2-\lam^2)^{-1}
= \mu(\xi)+ \lam^2\mu(\xi)(\xi^2-\lam^2)^{-1}
\]
are smooth and the Fourier transform yields   
\[
\chi_{>4a}(|D|)\Gg_{-\lam}(x) 
= \frac{1}{2\pi} \hat{\mu}(x)+ \lam^2\mu(|D|)\Gg_{-\lam}(x).
\] 
It follows that  
\bqn \lbeq(6-2)
\chi_{>4a}(|D|)G_0(-\lam)\p(x)= \frac1{2\pi}(\hat{\mu}\ast \p)(x) 
+ \lam^2 \mu(|D|)G_0(-\lam)\p(x) 
\eqn 
and $\chi_{>4a}(|D|)\Wg(I^{\p,\ph}_{k,a})u(x)$ splits into 
the sum $T_1 u(x)+ T_2 u(x)$:  
\begin{gather}
T_1 u(x)= \frac1{2\pi}(\hat{\mu}\ast \p)(x)\cdot 
\int_0^\infty \la \ph, \Pi(\lam)u\ra g(\lam)^k 
\chi_{\leq 2a}(\lam) \lam d\lam \lbeq(T-1)
\\ 
T_2 u(x) = 
\mu(|D|)\int_0^\infty G_0(-\lam)|\ph \ra \la \p| \Pi(\lam)u \ra 
\r_k(\lam)\lam d\lam\,. \lbeq(T-2)
\end{gather}
where $\r_k(\lam) = \lam^2 g(\lam)^k \chi_{\leq 2a}(\lam)\in \Mg(\R^2)$.
The second factor on the right of \refeq(T-1) 
which we denote by $\ell(u)$  is a inear functional and   
\[  
\ell(u) = \frac1{2\pi} \int_{\R^2} (\Fg^{-1}\ph)(\xi)\hat{u}(\xi)
 g(|\xi|)^k \chi_{\leq 2a}(|\xi|) d\xi
=\la u, \ph \ast \Fg(g^k \cdot \chi_{\leq 2a})\ra
\]
Since $\hat{\m}\ast \p \in L^p(\R^2)$ for 
$1\leq p <\infty$ and 
$\Fg(g^k \cdot \chi_{\leq 2a})\in L^q(\R^2)$ for any $1<q<\infty$ 
by virtue of \reflm(logk), H\"older's inequality implies 
\[
\|T_1 u\|_p \leq C_p \|\hat{\m}\|_p  \|\ph\|_1 \|\p\|_1 \|u\|_p .
\]
We have $\|T_2 u\|_p \leq C\|\ph\|_1 \|\p\|_1 \|u\|_p$ 
by virtue of \reflm(rank1). Combining these two estimates 
we see that  
$\|\chi_{>4a}(|D|)\Wg(I^{\p,\ph}_{k,a})u\|_p$ is bounded by the right side of 
\refeq(Wpphmu-l). 

For $\chi_{\leq 4a}(|D|)\Wg(I^{\p,\ph}_{k,a})u$ we write 
\bqn \lbeq(6-20-a)
\chi_{\leq 4a}(|D|)G_0(-\lam)\p(x)
 = \frac1{2\pi} \int_{\R^2} 
\frac{e^{ix\xi}\chi_{\leq 4a}(|\xi|)\hat{\p}(\xi)}
{\xi^2 -(-\lam +i0)^2}d\xi  
\eqn 
where $+i0$ means that we first change $i0$ by $i\s$ with $\s>0$ and let 
$\s\to 0$ after the integration. Since $\int_{\R^2}\p(x)dx=0$,  
\bqn \lbeq(exp-Fp)
\hat{\p}(\xi)= \frac1{2\pi} \int_0^1 
\left(\int_{\R^2}(-ix\xi)e^{-i{\th}x\xi}\p(x)dx\right)d\th,
\eqn 
and \refeq(6-20-a) becomes $\sum_{m=1}^2 $ of 
\begin{multline} \lbeq(6-21-a)
\frac{-i}{(2\pi)^2}\int_{\R^2}\left(\int_{\R^2}\int_0^1 
\frac{\xi_m\chi_{\leq 4a}(|\xi|)e^{i(x-{\th}z)\xi}z_m \p(z)}
{\xi^2 -(-\lam +i0)^2}d\th dz \right) d\xi \\
= \frac{-i}{2\pi}\int_0^1 \int_{\R^2}z_m \p(z) \tau_{{\th}z}
\left(\frac1{2\pi}\int_{\R^2}e^{ix\xi}\frac{\xi_m\chi_{\leq 4a}(|\xi|)}
{\xi^2 -(-\lam +i0)^2} d\xi \right)d\th dz. 
\end{multline}
Here the change of order of integral is trivially justified 
if $i0$ is replaced by $i\s$;  the change of the order 
of the limit $\s\to  0$ and the integral with respect to $dz d\th$ 
can also be justified by observing that $z_m \ph(z) \in L^1(\R^2)$ 
and that, for $\a\leq \lam\leq 2a$ outside of which the integrand 
vanishes if $\Pi(\lam){u}=$ for $|\lam|< \a$ (recall $u \in \Dg_\ast$), 
\[
\frac1{2\pi}\int_{\R^2}e^{ix\xi}\frac{\xi_m\chi_{\leq 4a}(|\xi|)}
{\xi^2 -(-\lam +i\s)} d\xi 
= G_0(-\lam +i\s)\Fg(\xi_m\chi_{\leq 4a})(x) 
\]
is bounded in $\R^2$ uniformly for $0<\s<1$ by virtue of \refeq(hankel-est) 
and it converges uniformly as $\s \to +0$. Elementary estimate involving 
Cauchy's principal value implies that 
\[
\lim_{\s\to 0+}\int_{\R^2}\frac{f(\xi)(|\xi|-\lam)}{|\xi|^2-(-\lam+i\s)^2}d\xi 
= \int_{\R^2}\frac{f(\xi)}{|\xi|+ \lam}d\xi, \quad f \in \Sg(\R^2),
\]
and we may set in \refeq(6-21-a) as   
\[
\frac{\xi_m}{\xi^2-(-\lam+i0)^2}= 
\frac{\lam\w_m}{\xi^2-(-\lam+i0)^2} + \frac{\w_m }{|\xi|+ \lam} , 
\quad \w_m=\frac{\xi_m}{|\xi|}.
\] 
Then the inner integral in the right of \refeq(6-21-a) becomes 
\bqn \lbeq(rmchi-a)
2\pi R_m \lam \chi_{\leq 4a}(|D|)\Gg_{-\lam}(x) + R_m\left(\frac1{2\pi} 
\int_{\R^2} e^{ix\xi} \frac{\chi_{\leq 4a}(|\xi|)}{|\xi|+\lam}d\xi\right).
\eqn
The contribution  
of the first term of \refeq(rmchi-a) for   
$\chi_{\leq 4a}(|D|)\Wg_{k,a}^{\p,\ph}u$ is given by 
$\sum_{m=1}^2 \int_0^1 d\th$ of 
\begin{align}
& W_{m,\th}^{(1)}u(x)\stackrel{\rm def}{=}  
-i\int_{\R^2}z_m \p(z) \ph(y) \tau_{{\th}z} 
R_m \chi_{\leq 4a }(|D|) \notag \\
& \hspace{1cm} \times 
\left\{
\int_0^\infty \Gg_{-\lam}(x) \left(
\int_{{\mathbb S}^1} (\Fg \t_{-y}u )(\lam\w)d\w 
\right) \lam g(\lam)^k \chi_{\leq 2a}(\lam)\lam {d\lam} 
\right\}dzdy  \notag \\
& = 
-2\pi{i}
\int_{\R^2}z_m \p(z) \ph(y) 
(R_m \chi_{\leq 4a}(|D|)K(\t_{-y}\k(|D|)u)(x-\th{z})dz dy  \lbeq(40)
\end{align}
where $\k(\lam)= \lam g(\lam)^k \chi_{\leq a}(\lam)$ is a 
good multiplier. Then, Minkowski's inequality, 
\refeq(K-est) and multiplier theory imply 
\bqn 
\Big\|\sum_{m=1}^2 \int_0^1 W_{m,\th}^{(1)} d\th\Big\|_p
\leq C\|\ax \p\|_1 \|\ph\|_1 \|u\|_p, \quad 1<p<\infty 
\eqn 
We define $F(x)$ by 
\[ 
F(x)= \int_0^\infty 
\left(\frac1{2\pi}
\int_{\R^2} e^{ix\xi} \frac{\chi_{\leq 4a}(|\xi|)}{|\xi|+\lam}d\xi
\right)
\la \ph, \Pi(\lam)u \ra g(\lam)^k \chi_{\leq 2a}(\lam)\lam d\lam .
\] 
Then, the contribution 
of the second term of \refeq(rmchi-a) to 
$\chi_{\leq 4a}(|D|)\Wg_{k,a}^{\p,\ph}u$
is given by 
\bqn \lbeq(W-2-def)
W^{(2)}u(x)= 
\frac{-i}{2\pi} 
\sum_{m=1}^2 \int_0^1 \int_{\R^2} z_m \p(z)(R_m F)(x-\th z)dz d\th  
\eqn 
and Minkowski's inequality implies 
\bqn \lbeq(W-2-mid)
\|W^{(2)}u\|_p \leq C \|\az\p\|_1 \|F\|_p, \quad 1<p<\infty.
\eqn 
Thus, it suffices to show $\|F\|_p \leq C \|\ph\|_1 \|u\|_p$. We have 
\begin{align}
& F(x)= \int_{\R^2}
\left(\int_{\R^2} e^{ix\xi} \frac{\chi_{\leq 4a}(|\xi|)}{|\xi|+|\eta|}d\xi\right)
\check{\ph}(\eta)\hat{u}(\eta)g(|\eta|)^k \chi_{\leq 2a}(|\eta|)d\eta \notag 
\\
& = \int_{\R^2}\left( 
\int_{\R^4} e^{ix\xi-iy\eta} 
\frac{\chi_{\leq 4a}(|\xi|)g(|\eta|)^k\chi_{\leq 2a}(|\eta|)\check{\ph}(\eta)}
{|\xi|+|\eta|}d\xi d\eta \right)u(y) dy. \lbeq(tildeL-a)
\end{align}
Denote by $Q(x,y)$ the integral inside the parenthesis in \refeq(tildeL-a) 
and define 
\begin{gather}
\tilde{L}(x,y)= \int_{\R^4} e^{ix\xi-iy\eta} 
\frac{\chi_{\leq 4a}(|\xi|)\chi_{\leq 4a}(|\eta|)}
{|\xi|+|\eta|}d\xi d\eta, \\
{L}(x,y)= \int_{\R^2} \tilde {L}(x,y-y')\Fg(g^k\chi_{\leq 2a})(y') dy'
\end{gather}
so that 
\bqn 
Q(x,y) = \int_{\R^2}L(x,y-z)\ph(z) dz 
\eqn 
We show $\|\|Q(x,y)\|_{L^p(\R^2_x)}\|_{L^q(\R^2_y)} \leq C \|\ph\|_1$ 
for $1<p<\infty$ and its dual exponent $q=p/(p-1)$ which 
will via Minkowski's and H\"older's inequalities 
yield the desired estimate:
\[
\|F\|_p \leq \int_{\R^2} \|Q(x,y)\|_{L^p(\R^2_x)}|u(y)| dy 
\leq \|\|Q(x,y)\|_{L^p(\R^2_x)}\|_{L^q(\R^2_y)}\|u\|_p 
= C \|u\|_p.
\]
Since $\ph \in L^1(\R^2)$, it suffices to show  
\bqn \lbeq(target)
\|L(x,y)\|_{L^p(\R^2_x)}\in L^q(\R^2_y).   
\eqn 
It is shown that 
$\tilde{L}(x,y)\absleq C \ax^{-1}\ay^{-1}(\ax+ \ay)^{-1}$ 
in Lemma B in Appendix of \cite{Ya-point}  and elemetray estimates 
produce 
\[
\|\tilde{L}(x,y)\|_{L^p_x}
\leq C \left\{\br{ll} \ay^{-2}\, ,  & \ \ p >2, \\
\ay^{-2}\la \log y \ra^\frac12\, , & \ \ p=2,    \\
\ay^{\frac{2}{p}-3}\, , & \ \ 1<p<2\,. 
\er \right.
\]
It follows by Minkowski's inequality that 
$\|L(x,y)\|_{L^p(\R^2_x)}$ is bounded by  
\[
\int_{\R^2} \frac{\la \log |y'| \ra^k dy'}{\la y-y'\ra^2 \la y'\ra^2}, \  
\int_{\R^2} \frac{\la \log |y-y'| \ra \la \log |y'| \ra^k}
{\la y-y'\ra^2 \la y'\ra^2}dy', \  
\int_{\R^2} \frac{\la \log |y'| \ra^k dy'}
{\la y-y'\ra^{3-\frac{2}{p}} \la y'\ra^2} 
\] 
for $p>2$, $p=2$ and $1<p<2$ respectively. Then, Young's inequality implies 
that 
$\|L(x,y)\|_{L^p(\R^2_x)}\in L^q(\R^2_y)$ as desired. 
\edpf 

\section{High energy estimate} 

In this section,  we prove that $\W_{{\rm high},2a}$ defined by 
\refeq(sta-0-high):
\bqn \lbeq(high)
\W_{{\rm high},2a}u = \int_0^\infty 
G_0(-\lam)vM(\lam)^{-1}v \Pi({\lam})u \chi_{> 2a}(\lam)\lam d\lam\,.
\eqn 
is a good operator for any $a>0$. This has been 
known for some years under a slightly strongler assumption that 
for an $\ep>0$, $\ax^{\frac{7}{2}+\ep}|V(x)|\leq C $ 
(cf. \cite{Y-2dim}, Lemma 3.2). We give a new proof which replaces 
the argument via integration parts in \cite{Y-2dim} by the one via the 
Fourier multipliers and the singular 
integral operator $K$ of \refeq(DFEK). 

\bgth \lbth(high) Assume $\ax^2 V \in L^{\frac43}(\R^2)$. 
For any $a>0$, $\W_{{\rm high},2a}$ is a good operator. 
\edth 

We shall often omit the index $2a$ in $\W_{{\rm high},2a}$ and etc. in 
what follows when no confusion is feared.  
Expanding $M(\lam)^{-1}=U(1+ vG_0(\lam)w)^{-1} $ 
as in \refeq(mlam) shows that 
$\W_{\rm high}u$ is the sum of six operators: 
\begin{align} 
& \sum_{j=0}^5 \W_{{\rm high}}^{(j)} u   = \sum_{j=0}^4 
\int_0^\infty G_0({\lam})w (-vG_0(\lam)w)^j v 
\Pi({\lam})u \lam \chi_{>2a}(\lam) d\lam  \notag \\
& + \int_0^\infty G_0({\lam})w (vG_0(\lam)w)^5 
(1+ vG_0(\lam)w)^{-1}v \Pi({\lam})u \lam \chi_{>2a}(\lam) d\lam.
\lbeq(W12345)
\end{align}
and we show that $\W_{{\rm high}}^{(j)}$, $0\leq j\leq 5$ 
are good operators 
separately.  

\subsection{Estimate of $\W_{{\rm high},0}$} 
\bgprop \lbprop(Whigh0) Under the asumption of \refth(high) 
$\W_{{\rm high},0}$ is a good operator. 
\edprop
\bgpf Since $\W_{{\rm high},0}= W(M_V)$, \refprop(funda) (1) implies   
\bqn \lbeq(Wh-0)
\|\W_{{\rm high}}^{(0)} u\|_p \leq C_p \|V\|_1 \|u\|_p , \quad 1<p<\infty.
\eqn
Since $V\in L^1(\R^2)$ if $\ax^2 V \in L^{\frac43}(\R^2)$, 
the proposition follows. 
\edpf

\subsection{Estimate of $\W_{{\rm high}}^{(1)}$} 
In this subsection we prove 

\bgprop \lbprop(Whigh2-estimate) For any $1<p<\infty$ 
\bqn \lbeq(Wh-1)
\|\W_{{\rm high}}^{(1)}u\|_p \leq 
C_p \|u\|_p 
\iint_{\R^4}|V(x)||V(y)|(1+ |\log |x-y||)dx dy.
\eqn 
\edprop

Define $V_y^{(2)}(x)  = V(x)V(x-y)$. 
For a.e. $x\in \R^2$, $V_y^{(2)}\in L^1(\R^2)$. 
\bglm \lblm(Whigh2-2) We have the following expression 
for a.e. $x\in \R^2$:    
\bqn \lbeq(whigh2a)
\W_{{\rm high}}^{(1)}u(x) = - \int_{\R^2}(W(M_{V_y^{(2)}})
\Hg(|y||D|) \chi_{>2a}(|D|)\tau_y u)(x) dy. %\ \ u \in \Dg_\ast. 
\eqn 
\edlm 
\bgpf For $\lam\not=0$, $\Gg_{-\lam}(x-y)V(y) \in L^1(\R^2_y)$ 
and for $u \in \Dg_\ast$  
\begin{align}
VG_0({\lam})V u(x) 
& = \int_{\R^2} V(x)\Hg(\lam|y|)V(x-y) u(x-y)dy  \notag \\  
&= \int_{\R^2} V^{(2)}_y(x)\Hg(|y|\lam)(\tau_y u) (x) dy,  
\ \ \mbox{a.e.}\ x\in \R^2.
\lbeq(v1Gg)
\end{align} 
It follows that  
\begin{multline*} \lbeq(W0-exp)
- \W_{{\rm high}}^{(1)}(x) = \int_0^\infty (G_0(-{\lam})VG_0(\lam)V 
\Pi(\lam)u)(x) \lam \chi_{>2a}(\lam) d\lam \\
=
\int_0^\infty \left(\int_{\R^4} \Gg_{-{\lam}}(x-y) 
V^{(2)}_z(y) \Hg(\lam |z|) (\Pi({\lam}) \tau_z u)(y) 
dzdy \right)\lam \chi_{>2a}(\lam) d\lam. 
\end{multline*} 
We show that the last integral converges absolutely for a.e. 
$x\in \R^2$ by repeating 
the argument used for \refeq(hDast): We let $B_R \subset \R^2_x$. 
\reflm(first-step) implies $\Pi({\lam})(\tau_z u)(y)$ vanishes 
for $\lam \not\in (\a,\b)$ and is bounded; \refeq(hankel-est) that  
$|\Hg(|z|\lam)|\leq (|\log |z||+ C)$ for $\lam \in (\a,\b)$. 
Then, \refeq(hankel-int) implies that 
\begin{align*}
& \int_0^\infty \int_{B_R \times \R^4}|\Gg_{-\lam}(x-y) 
V^{(2)}_z(y) \Hg(|z|\lam)\Pi({\lam})(\tau_z u)(y)|
\lam \chi_{>2a}(\lam) dxdzdyd\lam \\ 
& \leq C \int_\a^\b \int_{B_R \times \R^4}
|\Gg_{-\lam}(x-y)||V^{(2)}_z(y)|(|\log |z||+ C) dxdzdyd\lam \\
& \leq C\int_{\R^4}|V(y) V(z)|(1+|\log |y-z||)dzdy <\infty.
\end{align*}
Thus, we can integrate with respect to $d\lam$ first. Then, 
applying \refeq(mult) for $f(\lam)=\chi_{>2a}(\lam)\Hg(|z|\lam)$ 
we obtain 
\[
\W_{{\rm high}}^{(1)}u(x)=
- \int_{\R^2} 
\Big( \int_0^\infty (G_0(-{\lam})V^{(2)}_y \Pi({\lam}) 
\Hg(|y||D|)\chi_{>2a}(|D|) \tau_y u)(x) \lam  d\lam \Big) dy  
\]
which is nothing but \refeq(whigh2a). This proves the lemma. 
\edpf 

We define for an $a>0$ 
\[
\Hg (\lam)= \Hg_{\rm low}(\lam) + \Hg_{\rm high}(\lam)
\stackrel{\rm def}{=} 
\chi_{\leq a}(\lam)\Hg(\lam)+ \chi_{> a}(\lam)\Hg(\lam). 
\] 
\bglm \lblm(H02y) Let $a>0$ and $y \in \R^2$.  
$\Hg_{{\rm high}}(|y||D|)$, 
$\Hg_{{\rm low}}(|y||D|)\chi_{>2a}(|D|)$ and 
$\Hg_{{\rm low}}(|y||D|)\chi_{>2a}(|D|)$ are all good operator and 
\begin{align} 
& \|\Hg_{{\rm high}}(|y||D|)\|_{\Bb(L^p)} \leq C_p. \lbeq(H02y-add)
\\
& \|\Hg_{{\rm low}}(|y||D|)\chi_{>2a}(|D|)\|_{\Bb(L^p)} 
\leq C_{a,p} (1+ |\log|y||). 
\lbeq(H02y) \\
& \|\Hg(|y||D|)\chi_{>2a}(|D|)\|_{\Bb(L^p)} 
\leq C_{a,p} (1+ |\log|y||)\,. \lbeq(H02y-a)
\end{align}
\edlm 
\bgpf (1) Since $\Hg_{\rm high}(\lam)$ satisfies \refeq(large-lam), 
the theory of  spatially homogenous Fourier integral operators 
(p.138 of \cite{Peral}, see also \cite{SSS,Tao}) implies 
that $\Hg_{\rm high}(|D|)$ is a good operator. Then, the scaling argument 
implies the same for $\Hg_{{\rm high}}(|y||D|)$ and estimate 
\refeq(H02y-add). \\[3pt]
(2) \refeq(small-lam) implies 
$(\Hg(\lam)-g(\lam))\chi_{\leq a}(\lam)\in \Mg(\R^2)$. 
It follows that $\Hg_{\rm low}(|y||D|)- g(|y||D|)\chi_{\leq a}(|y||D|)$ 
is a good operator with $y$-independent $\Bb(L^p)$-norm. 
Thus, it suffices to prove \refeq(H02y) for 
$g(|y||D|)\chi_{\leq a}(|y||D|)\chi_{>2a}(|D|)$ which is equal to%t or for  
\[
(g(|y|)- (2\pi)^{-1}{\log |D|})\chi_{\leq a}(|y||D|)\chi_{>2a}(|D|).
\]
However, $\|g(|y|)\chi_{\leq a}(|y||D|)\chi_{>2a}(|D|)\|_{\Bb(L^p)}
\leq C(1+ |\log|y||)$ is evident and we need consider only 
$-(1/2\pi)(\log |D|) \chi_{\leq a} (|y||D|)\chi_{>2a}(|D|)$. 
Define $F(\lam)= (1/2\pi)(\log \lam) \chi_{\leq a} (|y|\lam) 
\chi_{>2a}(\lam)$. Then, we have 
\bqn \lbeq(Fj)
|F^{(j)}(\lam)|\leq C\lam^{-j}(1 + |\log |y||). \quad j=0,1,2.
\eqn 
Indeed $F(\lam)\not=0$ only if $|y|<1$ and $a <\lam <a/|y|$ 
and,  
\[
|F(\lam) |\leq (2\pi)^{-1}\max(|\log a|, |\log a/|y||) 
\leq (2\pi)^{-1}(|\log |y|| + |\log a|),  
\]
which implies \refeq(Fj) for $j=0$. The proof for $j=1$ and $j=2$ is 
similar. Then, Mikhlin's theorem implies that 
$\|F(|D|)\|_{\Bb(L^p)}\leq C(1+ |\log |y||)$ and 
\refeq(H02y) follows. 
\refeq(H02y-add) and \refeq(H02y) implies \refeq(H02y-a). 
\edpf 

\noindent 
{\bf Proof of \refpropb(Whigh2-estimate)} We apply 
Minkowski's inequality to \refeq(whigh2a). Then, \refprop(funda) (1)  
and \refeq(H02y-a) implies that  
\[
\|\W_{{\rm high}}^{(1)}u\|_p \leq 
C_p \int_{\R^2}\|V^{(2)}_y\|_1 (1+|\log|y|) \|u\|_p dy 
\] 
which is equivalent to \refeq(Wh-1)
\qed 

\subsection{Estimate of $\W_{{\rm high}}^{(n)}$} 
The following proposition implies that 
$\W_{{\rm high}}^{(n)}$, $n=2,3,4$ are good operators. 
Note that $vG_0(-\lam)w\in \Hg_2$ for $\lam>0$. 
We define for $n=0,1,\dots$ 
\begin{gather} \lbeq(W-n)
\W^{(n)}u= \int_0^\infty G_0({\lam})w (-vG_0(\lam)w)^j v 
\Pi({\lam})u \lam \chi_{>2a}(\lam) d\lam, \\
C^{(n)}= 
 \int_{\R^{2(n+1)}} \left(\prod_{i=1}^{n+1} |V(y_i)|\right) 
\left(\prod_{i=1}^{n}(1+ |\log |y_i-y_{i+1}||)\right) 
dy_1 \cdots dy_{n+1}. \notag 
\end{gather} 

\bgprop \lbprop(n-th-estimate) For $1<p<\infty$, 
there exists a constant $C_p$ such that 
\bqn  \lbeq(Wn-est)
\|\W^{(n)}u\|_p \leq C_p C^{(n)}\|u\|_p, \quad u \in \Dg_\ast.
\eqn 
For any $\ep>0$, $C^{(n)}\leq (C_\ep \|\ax^\ep V\|_1)^{n+1}$. 
\edprop 
\bgpf We let $n \geq 2$. Then, repeating the argument which is used 
for obtaining \refeq(v1Gg) implies that for $u \in \Dg_\ast$   
\begin{align*}  
& (v U(-vG_0(\lam)w)^n v u)(x)= (-1)^n 
(V G_0(\lam)V \cdots VG_0(\lam)Vu)(x) 
\notag \\
& = (-1)^n \iint_{\R^{2n}}
V^{(n+1)}_{y_1,\dots, y_n}(x)
\left(\prod_{j=1}^{n} \Hg(\lam|y_j|)\right)\tau_{y_1+\cdots+y_n}u(x) \, 
dy_1 \dots dy_n 
\end{align*}
where 
$V^{(n+1)}_{y_1,\dots,y_n}(x)= V(x) V(x-y_1) \cdots V(x-y_1-\cdots-y_n)$. 
This yields 
\begin{multline}
\W^{(n)}u(x)= (-1)^n \iint_{\R^{2n}}\int_0^\infty \Gg_{-\lam}(x-y) 
V^{(n+1)}_{z_1,\dots, z_n}(y) \\
\times \left(\prod_{j=1}^{n} \Hg(|y_j|\lam)\right)
\Pi(\lam)\tau_{y_1+\cdots+y_n}u(y) \lam \chi_{>2a}(\lam)
dy_1 \cdots dy_n dy d\lam\,.  \lbeq(Wn)
\end{multline} 
The argument similar to the one which is 
used in the proof of \reflm(Whigh2-2) 
implies that \refeq(Wn) is absolutely integrable for a.e. 
$x\in \R^2$, which we avoid to repeat, and 
after applying \refeq(mult) we may integrate \refeq(Wn) by $d\lam$ first. 
This implies that $\W^{(n)}u$ is equal to  
\[
\iint_{\R^{2n}} W(M_{V^{(n+1)}_{y_1,\dots,y_n}})
 \chi_{>2a}(|D|)\prod_{j=1}^{n}
\Hg(|y_j||D|) \t_{y_1+\cdots+ y_n} {u} dy_1 \dots dy_n.   
\]
Minkowski's inequality, \refprop(funda) (1) and \reflm(H02y) then imply  
\[
\|\W^{(n)}u\|_p \leq C_p 
\int_{\R^{2n}}  
\|V^{(n+1)}_{y_1,\dots,y_n}\|_{L^1(\R^2_x)} 
\prod_{j=1}^n (1+ |\log |y_j||)\|u\|_p dy_1 \dots dy_n.
\]
The estimate $C^{(n)}\leq (C\|\ax^\ep V\|_1)^{n+1}$ is evident 
for $(1+ |\log|x-y||) \leq C_\ep \ax^{\ep/2}\ay^{\ep/2}$ 
for any $\ep>0$.  
\edpf 

\subsection{Estimate of $\W_{{\rm high}}^{(5)}$.} 
Estimate \refeq(Wn-est) implies that, if $\|\ax^4 V\|_{\frac43}$ is 
sufficiently small, the expansion 
$\W_{\rm high}= \sum_{j=0}^\infty \W^{(n)}$ 
converges in $\Bb(L^p)$ for any $1<p<\infty$ and 
$\W_{\rm high}$ becomes a good operator. However, 
we do not want the smallness assumption which makes $H$ 
automatically regular at zero. We shall instead exploit the decay 
property of $\|vG_0(\lam)w\|_{\Hg_2}$: 

\bglm \lblm(AB-dec)
{\rm (1)} Let $v, w \in L^\frac{8}{3}(\R^2)$. Then, 
$v G_0(\lam) w\in \Hg_2$ for any $\lam>0$ and for $\ep>0$ there exists a 
$C_\ep>0$ such that  
\bqn \lbeq(ABestima)
\|v G_0(\lam) w\|_{\Hg_2} \leq C_\ep \lam^{-1/2}\|v\|_{8/3}\|w\|_{8/3}\,, 
\quad  \lam \geq \ep.
\eqn
{\rm (2)} Let $j=1,2$. Suppose that 
$\ax^j{v}, {\ax}^j w \in L^\frac{8}{3}(\R^2)$. 
Then, $v G_0(\lam) w$ is an $\Hg_2$-valued $C^j$ function 
of $\lam\in (0,\infty)$. For any $\ep>0$, there exists a constant 
$C_\ep$ such that 
\bqn \lbeq(ABestimb)
\|(d^j/d\lam^j)v G_0(\lam) w \|_{\Hg_2} \leq C_\ep \lam^{-1/2}
\|\ax^j {v}\|_{8/3}
\|\ax^j {w}\|_{8/3}\,, \quad \lam \geq \ep .
\eqn
\edlm 
\bgpf (1) By virtue of \refeq(large-lam) 
and \refeq(small-lam), $\|v G_0(\lam) w\|_{\Hg_2}^2 \leq C(F_1+ F_2)$ where   
\begin{align*}
F_1&= \int_{|x-y|>\ep/\lam} \frac{|v(x)|^2|w(y)|^2 }{\lam |x-y|} dxdy, \\ 
F_2&=\int_{|x-y|<\ep/\lam} (|\log |\lam |x-y||+ C)^2 |v(x)|^2|w(y)|^2 dxdy. 
\end{align*}
The generalized Young's inequality (e.g. \cite{LS}) implies 
\begin{align*}
F_1 & \leq C \lam^{-1} \|v\|_{8/3}^2 \|w\|_{8/3}^2, 
\\
F_2 & \leq C \|v\|_{8/3}^2 \|w\|_{8/3}^2
\left(\int_{|z|<\ep/\lam} (|\log |(\lam z)||+ C)^2 dz \right)^2 \\
& \leq C_\ep \lam^{-1} \|\|v\|_{8/3}^2 \|w\|_{8/3}^2 .
\end{align*}
\refeq(ABestima) follows. \\[3pt]
(2) It is well known (see \cite{AS}, pp. 360-364) that 
$\Gg_{\lam}^{(1)}(x)=(d/d\lam)\Gg_{\lam}(x)$ and 
$\Gg_{\lam}^{(2)}(x)=(d^2/d\lam^2)\Gg_\lam(x)$ are 
respectively given by $i/4$ times 
\bqn \lbeq(first-derivative)
|x|\Big(\frac{d}{dz} H_0^{(1)}\Big)(\lam |x|)
= - |x|H_{-1}^{(1)}(\lam |x|)
\absleq C \left\{ 
\br{ll} \lam^{-\frac12}|x|^{\frac12}, \ & |x|\lam \geq 1\,. \\
\lam^{-1}, \ & |x|\lam<1\,.
\er
\right.
\eqn 
\begin{align}
|x|^2 \Big(\frac{d^2}{dz^2} H_0^{(1)}\Big)(\lam |x|)
& = |x|^2 
\left(H_{-2}^{(1)}(\lam |x|)+ 
\frac{1}{\lam |x|}H_{-1}^{(1)}(\lam |x|)\right) 
\notag \\
& \absleq C 
\left\{ 
\br{ll} \lam^{-\frac12}|x|^{\frac32}, \ & |x|\lam \geq 1. \\
\lam^{-2}, \ & |x|\lam<1. 
\er
\right.  \lbeq(2nd-derivative)
\end{align} 
Let $G_0^{(j)}(\lam)u(x)= (\Gg_{\lam}^{(j)}\ast u)(x)$ for $j=1,2$. 
We have for the integral kernels for a.e. $(x,y)\in \R^2 \times \R^2$ that 
\[
(v G_0^{(j-1)}(\lam) w)(x,y) - (v G_0^{(j-1)}(\mu) w)(x,y) 
= \int^\lam_{\mu} (v G^{(j)}_{(0)}(\r) w)(x,y) d\r\,.
\]
Then, if 
${\ax}^j v(x), {\ax}^j w(x) \in L^{\frac83}(\R^2)$, 
by virtue of \refeqs(first-derivative,2nd-derivative), 
the proof of (1) together with the dominated convergence theorem 
implies that $v G_0^{(j)}(\lam) w \in \Hg_2$ satisfies \refeq(ABestimb) 
with obvious modifications and that 
$\lam \mapsto v G_0^{(j)}(\lam) w \in \Hg_2$ is continuous 
for $0<\lam<\infty$. It follows hat 
$v G_0(\lam) w$ is an $\Hg_2$-valued $C^j$ function of 
$\lam\in (0,\infty)$ and $(d/d\lam)^j vG_0(\lam)w=v G^{(j)}_0(\lam) w$. 
The lemma follows. 
\edpf

\bgprop \lbprop(Wh5) Assume $\ax^2 V \in L^{\frac43}(\R^2)$. Then,  
$\W_{{\rm high}}^{(5)}$ is a good operator. 
\edprop
\bgpf Define $T(\lam) = - w(vG_0(\lam)w)^5 (1+ vG_0(\lam)w)^{-1}v$. 
It suffices by virtue of \refprop(funda) (5) to show that 
$T(\lam)\in \Og_{\Lg_1}^{(2)}(\lam^{-\ep})$ for some $\ep>0$ as 
$\lam\to \infty$. 
We first remark that for any decomposition $V(x) = v(x) w(x)$ 
such that $v, w \in L^{8/3}(\R^2)$ 
\bqn \lbeq(v-w-change)
v(1+ wG_0(\lam)v)^{-1} w =w(1+ vG_0(\lam)w)^{-1}v
\eqn
and they do not depend on how $V$ is decomposed.  
Indeed, \refeq(v-w-change) holds for large $|\lam|$ 
in the cone 
$\Gamma= \{\lam\in \Cb^{+} \colon 0\leq \arg \z <\pi/4\}$
because $\|vG_0(\lam)w\|_{\Hg_2}$ converges to zero as 
$\lam \to \infty$ in 
$\Gamma$ by the proof of \reflm(AB-dec) (1) and 
both sides become $V + V G_0(\lam) V+ \cdots$ by the Neumann expansion. 
Then, 
the analyticity in $\C^{+}\setminus \Eg$ and the continuity 
in $\Cb^{+}\setminus (\Eg \cup\{0\})$ imply that the same holds 
for all $\lam>0$, where $\Eg=\{i\k \colon \k>0, \k^2\in \s_p(H)\}$. 

Let first $V= v w$ with $v(x) = |V(x)|^{\frac12}$ and $w(x) = U(x) v(x)$ 
so that $\ax v, \ax w \in L^{\frac83}(\R^2)$. Then, \reflm(AB-dec) obviously 
implies that $T(\lam)\in \Og_{\Lg_1}^{(1)}(\lam ^{-\frac32})$ as 
$\lam \to \infty$ and that 
\begin{multline} \lbeq(d-1)
T'(\lam)= -\sum_{j=0}^4 w(vG_0(\lam)w)^{j}vG_0'(\lam)w
(vG_0(\lam)w)^{4-j} v \\
 - w(vG_0(\lam)w)^5 
(1+ vG_0(\lam)w)^{-1}vG_0'(\lam)w (1+ vG_0(\lam)w)^{-1}v. %\lbeq(d-1-a)
\end{multline}  
To see $T(\lam)\in \Og_{\Lg_1}^{(2)}(\lam ^{-\frac12})$, we further 
differentiate \refeq(d-1). The result is the sum of terms 
with contain two first derivatives $G_0'(\lam)$ and the ones with 
single $G_0''(\lam)$. Letting $V= v w$ with $v(x) = |V(x)|^{\frac12}$ 
and $w(x) = U(x) v(x)$, we see the former terms are $\Lg_1$-valued 
continuous functions of $\lam>0$ with $\Lg_1$-noms bounded by 
$C\lam^{-\frac52}$ as $\lam \to \infty$. 
To see the same holds for the latter terms,  
we decompose $V = vw$ in such a way that $v, \ax^2 w \in L^{\frac83}(\R^2)$ 
and sandwitch $G''(\lam)$ by $w$ like 
$w(vG_0(\lam)w)^{j-1}vG_0(\lam)v wG_0''(\lam)w(vG_0(\lam)w)^{4-j}v$. 
Thus, $T(\lam) \in \Og_{\Lg_1}^{(2)}(\lam^{-\frac12})$ and the proof is 
completed. 
\edpf 

\section{Low energy estimate}

We now study the low energy part $\W_{{\rm low},2a}$ defined by 
\refeq(sta-0-low):   
\bqn \lbeq(w-low)
\W_{{\rm low},2a} u= \int_0^\infty G_0(-{\lam})vM(\lam)^{-1}v 
\Pi(\lam)u \lam \chi_{\leq 2a}(\lam) d\lam, \ u \in \Dg_\ast.  
\eqn 
In this section. {\it we shall often omit the phrase ``for small $\lam>0$''}
as we shall exclusively work for small $\lam>0$, 
  
\subsection{Resonances} 
The following lemma which has mostly been proved by \cite{JN} 
under a slightly different assumptions gives the relation 
between resonances and the singularities of $H$. 
We assume $V(x)\not\equiv 0$ define $v=|V|^{\frac12}$ and $w=U v$ 
as in the introduction.   

\bglm \lblm(resonances) 
Suppose $\ax^{1+\ep} V \in L^1(\R^2)$ for an $\ep>0$. Then, 
$\Ng_\infty\not=\{0\}$ if and only if $H$ is singular at zero. 
 In this case 
\begin{itemize}
\item $u \in \Ng_\infty$ satisfies $\la w u,v\ra=0$ 
and $S_1 \HL = \{w u \colon u \in \Ng_\infty\}$. 
\item $\Ng_\infty \ni u \mapsto \z=w u\in S_1 \HL$ 
is an isomorphism and the inverse map is given by 
$u= N_0 v\z-\|v\|^{-2} \la PT_0 S_1\z, v\ra$. 

\item $u\in \Ng_\infty$ satisfies as $|x|\to \infty$ that 
\begin{gather} 
u(x) = c + \sum_{j=1}^2 \frac{x_j}{|x|^2} 
\left( \frac{1}{2\pi} \int_{\R^21} y_j V(y) u(y) dy \right)  
+ O(|x|^{-1-\ep}), \lbeq(b-sol) \\
c= \|v\|_2^{-2} \la PT_0 S_1 wu, v \ra.  \lbeq(b-sol-ext) 
\end{gather} 
\end{itemize}
\ben 
\item[{\rm (1)}] $H$ has singularities of the first kind at zero 
if and only if $\Ng_\infty$ consists only of $s$-wave resonances.  
In this case $\rank T_1= \dim\Ng_\infty=1$. 
\item[{\rm (2)}] $H$ has singularities of the second kind at zero 
then $\Ng_\infty$ consists of $s$-wave and $p$-wave resonances 
but no zero energy eigenfunctions.  
In this case  $1 \leq \rank S_2 = \rank T_2 \leq 2$.
$\Ng_\infty$ consists only of $p$-wave resonances if and only if $T_1=0$. 
If $T_1\not=0$, then $u \in \Ng_\infty$ is a $p$-wave resonance 
if $wu \in S_2\HL$ and an $s$-wave resonance otherwise. 
\item[{\rm (3)}] 
$H$ has singularities of the third kind at zero  if and only if 
zero energy eigenfunctions exist. In this case 
$u\in \Ng_\infty$ is eigenfunction if $ wu \in S_3\HL$, 
$p$-wave resonance if $wu \in S_2\HL \setminus S_3\HL$ and 
$s$-wave resonance if $wu \in S_1\HL \setminus S_2\HL$.
\een
\edlm 
\bgpf We give a proof for readers' convenience. 
Suppose $u \in \Ng_\infty\setminus \{0\}$. Then, 
$(-\lap)(u+ N_0 Vu )=0$ and $u+ N_0 Vu$ is a harmonic polynomial. 
But, $u \in L^\infty$ implies $u+ N_0 Vu=O(\log |x|)$ as $|x|\to \infty$ 
and $u+ N_0 Vu=c$ for a constant. 
Hence, $N_0 Vu(x)\in L^\infty$ and it must be that $\int V u dx =0$ or $P(wu)=0$. 
Thus, $cv=(v+ vN_0V)u= T_0 Q (wu)$. It follows $QT_0 Q(wu) = 0$, 
viz $wu \in S_1\HL$ and 
$wu\not=0$ because $wu=0$ would imply $u=c$ and $w=0$ and hence $V=0$, 
which is a contradiction. Hence $H$ is singular at zero. 
Moreover, $u+ N_0 Vu=c$ and $\int_{\R^2}V(x)u(x) dx =0$  imply  \refeq(b-sol);  
from $T_0 Q (wu) =c v$ we obtain \refeq(b-sol-ext). 

Assume conversely that $QT_0Q\z = 0$ for a $\z\in Q\HL\setminus\{0\}$.  Then, 
$(U+ vN_0 v)\z= cv$ for a constant $c$ and $u \stackrel{\rm def}{=}N_0 v\z-c$ 
satisfies both 
$U\z+ vu=0$ and $-\lap u = v\z$. It follows $(-\lap + V)u=0$ 
and $\la v,\z\ra=0$ implies $u\in L^\infty(\R^2)$. 
If $u=0$, then $U\z= - vu=0$ and $\z=0$. It follows that 
$u \in \Ng_\infty \setminus\{0\}$,  
$\z = -w u$, viz. $\z\in \{wu \colon u \in \Ng_\infty\}$ 
and that $\Ng_\infty \ni u \to w u \in S_1 L^2$ is an isomorphism. 
The first part of the lemma follows. 

(1) \refeq(b-sol-ext) implies $c\not=0$ for all $u \in \Ng_\infty$ 
if and only if $PT_0 S_1Q \not=0$ on $S_1 \HL$, or $T_1$ is non-singular 
on $S_1 \HL$. 

(2) If $\z \in S_2\HL \setminus \{0\}$, then $\la v, \z\ra=0$ and, since 
$T_2= S_2 vG_1 v S_2$ is non-singular
\bqn \lbeq(G1-inner)
\la vG_1 v \z, \z\ra = - \sum_{j=1}^2 
\left| \int_{\R^2} x_j v(x)\z(x) dx \right|^2 >0 . 
\eqn 
Moreover, $PT_0 S_1 \z=0$ implies $c=0$. Thus, the corresponding 
$u= N_0 v\z -c \Ng_\infty$ is a $p$-wave resonance. 
\refeq(G1-inner) implies $1\leq \rank S_2=\rank T_2\leq 2$.  

(3) In this case, $u = -w\z$, $\z \in S_3 \HL$ is an eigenfunction. 
The rest of the statement is obvious from the proof of (1) and (2). 
\edpf  

\subsection{Preliminaries} 
For studying the behavior of $M(\lam)^{-1}$ as $\lam\to 0$, 
we repeatedly use Feshbach formula and the lemma due to 
Jensen and Nencie (\cite{JN}) which we recall here. 
Let $A$ be the operator matrix
\[
A = 
\begin{pmatrix} a_{11}  & a_{12} \\ a_{21}  & a_{22} \end{pmatrix} 
\] 
on the direct sum of Banach spaces $\Yg= \Yg_1 \oplus \Yg_2$. 

\bglm \lblm(FS) Suppose $a_{11}$, $a_{22}$ are closed 
and $a_{12}$, $a_{21}$ are bounded operators. 
Suppose that $a_{22}^{-1}$ exists. Then $A^{-1}$ exists if and only if 
$d= (a_{11}- a_{12}a_{22}^{-1}a_{21})^{-1}$ exists. In this case we have 
\bqn \lbeq(FS-formula)
A^{-1} = \begin{pmatrix}  d & -d a_{12} a_{22}^{-1} \\
-a_{22}^{-1}a_{21} d & a_{22}^{-1}a_{21}d a_{12} a_{22}^{-1} + a_{22}^{-1}
\end{pmatrix}.
\eqn 
\edlm 

\begin{lemma}[\cite{JN}] \lblm(JN) 
Let $A$ be a closed operator in a Hilbert space $\Xg$ 
and $S$ a projection. Suppose $A+S$ has a bounded inverse. 
Then, $A$ has a bounded inverse if and only if 
\[
B= S - S(A+S)^{-1}S 
\]
has a bounded inverse in $S\Xg$. In this case, 
\bqn \lbeq(JN-1)
A^{-1}= (A+S)^{-1}+ (A+S)^{-1}SB^{-1}S (A+S)^{-1}.
\eqn 
\end{lemma}
 
Recall that $g(\lam|x|)= N_0(x)+ g(\lam)$. 
The expansion \refeq(hankel) implies  
\begin{align} 
& \pa_\lam^{j} (\Gg_{\lam}(x)-g(\lam|x|))
\absleq C \pa_\lam^j (\la g(\lam)\ra \lam^2)|x|^2 \la\log|x|\ra, 
\lbeq(hankel-add-1) \\
& \pa_\lam^j \big(
\Gg_{\lam}(x)-g(\lam|x|)(1 -\tfrac14(\lam|x|)^2)
+ \tfrac1{8\pi}(\lam|x|)^2\big) \notag \\
& \hspace{4cm} \absleq 
C \pa_{\lam}^j (\la g(\lam)\ra \lam^4)|x|^4 \la\log|x|\ra\,. 
\lbeq(hankel-add-2)
\end{align} 
The following lemma has been proved in \cite{EG} (see also \cite{EGG}) 
under slightly different assumptions. We set 
\bqn \lbeq(g1) 
g_1(\lam)\stackrel{\rm def}{=} g(\lam)\|V\|_1.
\eqn 

\bglm \lblm(step-0) 
{\rm (1)} Suppose $\ax^\c V \in L^{1}(\R^2)$, $\c>4$. 
Then, as $\lam \to 0$,  
\bqn \lbeq(M-1)
M(\lam) = g_1(\lam)P + T_0 
+ M_0(\lam), \ \ M_0(\lam)= \Og_2(g(\lam)\lam^{2}).
\eqn 
{\rm (2)} Suppose $\ax^\s V \in L^1(\R^2)$  for some $\s>8$. Then, as 
$\lam \to  0$, 
\bqn \lbeq(M-2)
M_0(\lam)= -g(\lam)\lam^2 v G_1 v -\lam^2 vG_2 v + 
\Og_2(\lam^{4}\la \log\lam\ra) \,.
\eqn 
\edlm 
\bgpf We may assume $0<\lam<1$. 
\refeq(hankel-add-1) implies that for $j=0,1,2$  
\[
\left(\int_{|x-y|\lam<1}
 |(d/d\lam)^j M_0(\lam,x,y)|^2dx dy\right)^\frac12 \leq 
C |g(\lam)\lam^{2-j}| \|\ax^\c V\|_1.
\]
For $\lam |x-y|\geq 1$, \refeqs(first-derivative,2nd-derivative) imply 
\[
\frac{d^j}{d\lam^j}(M_0(\lam)+v(x)g(\lam|x-y|)v(y))
=O(\lam^{-1/2}|x-y|^{j-1/2})v(x)v(y)
\]
for $j=0,1,2$. Since  $\ax^{-1} \ay^{-1} \leq C \lam$, we have  
\begin{gather*}
\int_{|x-y|>\lam^{-1}} |g(\lam)|^2|v(x)v(y)|^2 dxdy 
\leq C \lam^{4}|g(\lam)|^2 \|\ax^4 V\|_1^2, \\
\int_{|x-y|>\lam^{-1}} |\log |x-y||^2 |v(x)v(y)|^2 dxdy 
\leq C \lam^{4}|g(\lam)|^2 \|\ax^\c V\|_1^2 ,\\
\int_{|x-y|>\lam^{-1}} \lam^{-1} |x-y|^{2j-1}|v(x)v(y)|^2 dxdy 
\leq C \lam^{\c-2j}\|\ax^{\c} V\|_1^2.
\end{gather*}
These estimates yield the first statement of the lemma. 
Proof of (2) is similar and we omit the repetitious details.  
\edpf 

We often use the following trivial but important identity:
\bqn \lbeq(X+1)
(1+ X)^{-1} = 1 - X(1+X)^{-1}= 1 - X + X(1+X)^{-1}X. 
\eqn 

\subsection{The case $H$ is regular at zero} 
In this section we assume $QT_0Q$ is non-singular on $Q\HL$ and 
prove the following theorem, which together with \refth(high) implies 
that $W_{\pm}$ are bounded in $L^p(\R^2)$ for all $1<p<\infty$ 
if $V \in \ax^{-2}L^{\frac43}(\R^2) \cap \ax^{-\c}L^1(\R^2)$ for a $\c>4$ 
and if $H$ is regular at zero. 
This slightly improves the result of  \cite{Y-2dim,JY-2} 
which assumes $|V(x)|\leq C\ax^{-6-\ep}$, $\ep>0$.  

\bgth \lbth(regular-case)
Suppose that $\ax^\c V \in L^1(\R^2)$ for a $\c>4$ 
and that $H$ is regular at zero. Then, $\W_{{\rm low},2a}$ is a 
good operator for any $a>0$. 
\edth 

\refth(regular-case) immediately follows from the next lemma and 
\refprop(funda). We define $v_{\ast}=v/\|v\|_2$,  
\bqn \lbeq(h-def)
h(\lam){=}  (g_1(\lam) + c_1)^{-1}, \quad c_1 =
\la v_{\ast}|T_0- T_0Q(QT_0Q)^{-1}QT_0|v_{\ast}\ra.
\eqn 
Recall that $\Bg=\{M_m+ T \colon m \in L^\infty(\R^2), T \in \Hg_2\}$.

\bglm \lblm(FS-regular) 
We have $Q(QT_0Q)^{-1}Q\in \Bg$ and there exists an 
operator $L$ of rank at most two such that 
\bqn \lbeq(M-reg-case)
M(\lam)^{-1} = h(\lam)L + Q(QT_0Q)^{-1}Q + \Og_2(g\lam^2) \quad \lam \to 0. 
\eqn  
\edlm 
\bgpf  It has been proved (\cite{Schlag}) that $Q(QT_0Q)^{-1}Q\in \Bg$.
In view of \refeq(M-1), we first show 
$g_1(\lam)P+ T_0$ is invertible.  
In the decomposition $L^2(\R^2) = P L^2(\R^2) \oplus Q L^2(\R^2)$, 
\bqn \lbeq(Decom-FS)
g_1(\lam)P + T_0 = 
\begin{pmatrix} g_1(\lam)  + PT_0P & PT_0Q \\ QT_0P & Q T_0 Q \end{pmatrix}.
\eqn
Here $a_{22}\stackrel{\rm def}{=} QT_0Q$ is invertible by 
the assumption; for small $\lam>0$, 
\[
g_1(\lam)P + PT_0 P - PT_0Q a_{22}^{-1}QT_0P = (g_1(\lam)+ c_1)P 
\]
is invertible in $PL^2(\R^2)$ and  
$(g_1(\lam)P + PT_0 P - PT_0Qa_{22}^{-1}QT_0P )^{-1}=h(\lam)P $. Then, 
\reflm(FS) implies that $(g_1(\lam)P + T_0)^{-1}$ exists and it is equal to  
\begin{align}
& h(\lam) \begin{pmatrix} P & -PT_0Qa_{22}^{-1} \\
-a_{22}^{-1}QT_0P  & a_{22}^{-1}QT_0PT_0Q a_{22}^{-1} 
\end{pmatrix} + Qa_{22}^{-1}Q\notag  \\
& \stackrel{\rm def}{=}h(\lam)L + Q(QT_0Q)^{-1}Q.    \lbeq(gPT0-inv)
\end{align}
It is obvious that $\rank L \leq 2$.  It follows from \refeq(M-1) that       
\[
M(\lam)= (1+ M_0(\lam) (g_1(\lam)P+T_0)^{-1})(g_1(\lam)P+T_0)
\]
and $M_0(\lam)(g_1(\lam)P+T_0)^{-1}= \Og_2(g(\lam)\lam^{2})$. 
Hence, for small $\lam>0$, the series converges in $\Hg_2$  and 
\begin{align*} 
M(\lam)^{-1}& = (g_1(\lam)P+T_0)^{-1}(1+ M_0(\lam) (g_1(\lam)P+T_0)^{-1})^{-1} 
\notag \\
& = (g_1(\lam)P+T_0)^{-1}+ 
\sum_{j=1}^\infty (g_1(\lam)P+T_0)^{-1} (M_0(\lam) (g_1(\lam)P+T_0)^{-1})^j, 
\end{align*}
which implies \refeq(M-reg-case) by virtue of \refeq(gPT0-inv). 
\edpf 

\subsection{The case $H$ is singular at zero. Threshold analysis 1} 
In the rest of the paper we assume that $H$ is singular at zero: 
$QT_0 Q$ is singular in $Q\HL$ and 
$S_1$ is the orthogonal projection in $QL^2(\R^2)$ onto 
$\Ker_{QL^2}QT_0 Q$ and that $\ax^\c V \in L^1(\R^2)$ for $\c>4$, 
however, for the results which we need in the following subsections 
we assume $\ax^{\c} V(\R^2) \in L^1$ for a $\c>8$. 

\bglm[\cite{JN,EGG}] \lblm(preli)
${\rm Spec}(QT_0 Q\vert_{QL^2(\R^2)})$ is discrete outside $\{-1,1\}$. 

\noindent 
{\rm (1)} The projection $S_1$ is of finite rank. Let $n= \rank S_1$. 

\noindent
{\rm (2)} $D_0=((QT_0 Q+ S_1)\vert_{QL^2(\R^2)})^{-1}$ is of class $\Bg$. 
\edlm
\bgpf (1) is proved in \cite{JN}.   
Schlag(\cite{Schlag}) proves $D_0\in \Bg$.
\edpf 

\bgdf 
We take and fix an orthonormal basis $\{\z_1, \dots, \z_n\}$ 
of $S_1L^2$. We regard the 
projection $S_1 = \z_1 \otimes \z_1 + \cdots+ \z_n \otimes \z_n$ 
also as an orthogonal projection in $\HL$, viz. we identify 
$QS_1Q$ and $S_1$.
\eddf

We study $M(\lam)^{-1}$ as $\lam \to 0$ by applying 
\reflm(JN) to the pair $(M(\lam),S_1)$. We first study   
\bqn \lbeq(M-1-add)
M(\lam)+ S_1 = g_1(\lam) P+T_0+ S_1 + M_0(\lam).
\eqn 
Define  
\bqn \lbeq(def-h1)
h_1(\lam)=(g_1(\lam) + c_2 )^{-1}, \quad 
c_2 {=} \la v_{\ast}|T_0- T_0Q D_0 QT_0|v_{\ast}\ra. 
\eqn  
Evidently $h_1 \in \Mg(\R^2)$. Define the operator matrix 
\bqn \lbeq(def-L1)
L_1 = \begin{pmatrix} P & -PT_0Q D_0 \\
-D_0 QT_0P  &  D_0 QT_0PT_0Q D_0 
\end{pmatrix}.
\eqn 
in the decomposition $L^2(\R^2) = P L^2(\R^2) \oplus Q L^2(\R^2)$.  
$L_1$ is $\lam$-independent and ${\rm rank}\, L_1 \leq 2$. 

\bglm \lblm(s-wave-f1) 
$g_1(\lam) P+T_0+ S_1$ is invertible in $L^2(\R^2)$ and 
\bqn 
(g_1(\lam) P+T_0+ S_1)^{-1}= h_1(\lam) L_1 + Q D_0 Q. 
\lbeq(gPT0)   
\eqn 
We denote 
$N(\lam)=(g_1(\lam) P+T_0+ S_1)^{-1}$. $vN(\lam)v$ is a good producer. 
\edlm 
\bgpf We use \reflm(FS). 
In the decomposition $L^2(\R^2) = P L^2 \oplus Q L^2$, 
\bqn \lbeq(sing-f-0)
g_1(\lam)P + T_0 + S_1= 
\begin{pmatrix} g_1(\lam) P + 
PT_0P & PT_0Q \\ QT_0P & Q T_0 Q + S_1 \end{pmatrix}.
\eqn 
By virtue of \reflm(preli), $D_0=(Q T_0 Q + S_1)^{-1}$ exists and 
$D_0 \in \Bg$. It is obvious that for small $\lam>0$ 
\[
g_1(\lam)P + PT_0 P - PT_0Q D_0 QT_0P = (g_1(\lam) + c_2)P 
\]
has the inverse $h_1(\lam)P$. It follows that \refeq(sing-f-0) 
is invertible for small $\lam>0$ and \refeq(FS-formula) implies \refeq(gPT0). 
\edpf  

We define 
\bqn 
\Rg_1(\lam)\stackrel{\rm def}{=} v(G_1+ g(\lam)^{-1} G_2)v .  \lbeq(T1-Rg1-def)
\eqn  

\bglm \lblm(A+S-inv) 
$M(\lam)+S_1$ is invertible for small $\lam>0$ and  
\bqn 
(M(\lam)+ S_1)^{-1} =  N(\lam) + \Og_2(g(\lam) \lam^2), \quad \lam \to 0.
\lbeq(A+S1)
\eqn 
We denote $\Ag_0(\lam)= (M(\lam)+ S_1)^{-1}$. 
$v\Ag_0(\lam)v$ is a good producer.  
If $\ax^{\c}V \in L^1(\R^2)$ for a $\c>8$, then \refeq(A+S1) improves 
and as $\lam \to 0$  
\bqn 
\Ag_0(\lam)=  N(\lam) + 
\lam^2 g(\lam) N(\lam) \Rg_1(\lam) N(\lam) + \Og_2(\lam^4 g(\lam)^2).
\lbeq(A+S1-a)
\eqn 
\edlm 
\bgpf We have 
$M(\lam)+S_1 = (1+ M_0(\lam)N(\lam))(g_1(\lam)P + T_0 + S_1)$ 
by virtue of \reflm(s-wave-f1). \refeq(M-1) and \refeq(gPT0) 
then imply that it is invertible for small $\lam>0$ and 
\begin{multline} \lbeq(6-25a)
\Ag_0(\lam)=(M(\lam)+S_1)^{-1}= N(\lam) - N(\lam) M_0(\lam) 
N(\lam)  \\
+ N(\lam) M_0(\lam) N(\lam) 
(1+ M_0(\lam) N(\lam))^{-1}M_0(\lam) N (\lam).
\end{multline} 
This implies \refeq(A+S1) and the 
second line of \refeq(6-25a) is of $\Og_2(\lam^4 g(\lam)^2)$. 
In particular, $v\Ag_0(\lam)v$ is a good producer 
by virtiue of \refprop(funda). If $\ax^{\c}V \in L^1(\R^2)$ for a $\c>8$, 
then \refeq(M-2) implies 
$M_0(\lam)=-g(\lam) \lam^2 v(G_1+ g(\lam)^{-1}G_2)v + \Og_2(g(\lam) \lam^4)$ 
and \refeq(A+S1-a) follows from \refeq(6-25a). 
\edpf 

We define $B_1(\lam)$ on $S_1L^2(\R^2)$ by 
\bqn 
\lbeq(B1-def)
B_1(\lam) = S_1 - S_1 \Ag_0(\lam)S_1, \quad T_1=  S_1Q T_0 P T_0 Q S_1 .
\eqn 

\bglm \lblm(B1-expression) On $S_1\HL$, we have as $\lam \to 0$ 
\bqn 
B_1(\lam)= - h_1(\lam)(T_1 - \lam^2 X(\lam)),  
\quad X(\lam) \in  \Og_2(g(\lam)^2).  \lbeq(B1-s)
\eqn 
If $\ax^{\c}V \in L^1(\R^2)$ for a $\c>8$, then in the right of \refeq(B1-s)
\begin{multline}  
X(\lam)=
 -S_1 h_1(\lam)^{-1}g(\lam)\big(\Rg_1(\lam) 
+ h_1(\lam) (L_1\Rg_1(\lam)+ \Rg_1(\lam)L_1) 
\\
+ h_1(\lam)^2 L_1 \Rg_1(\lam) L_1 + \Og_2 (g(\lam)\lam^2)\big) S_1.  
\lbeq(X-def)
\end{multline} 
\edlm 
\bgpf We have $S_1 N(\lam) S_1 = S_1 + T_1$ since 
$S_1QD_0 QS_1= S_1$ and   
\[
S_1 {L}_1 S_1= \begin{pmatrix} 0 & 0 \\ 0  &  S_1 \end{pmatrix} 
\begin{pmatrix} 
P & -PT_0Q D_0 \\
-D_0 QT_0P  &  D_0 QT_0PT_0Q D_0  
\end{pmatrix}
\begin{pmatrix} 0 & 0 \\ 0  &  S_1 \end{pmatrix} = T_1 .
\]
Then, \refeq(gPT0) and \refeq(A+S1) imply \refeq(B1-s). 
If $\ax^{\c}V \in L^1(\R^2)$, we obtain \refeq(X-def) by using 
also \refeq(A+S1-a) .
\edpf

If $B_1(\lam)^{-1}$ exists in $S_1\HL$, then \reflm(JN) implies  
that   
\bqn \lbeq(1-s)
M(\lam)^{-1}= \Ag_0(\lam)+  \Ag_1(\lam), \ 
\Ag_1(\lam) \stackrel{\rm def}{=} \Ag_0(\lam)S_1 B_1(\lam)^{-1}S_1 
\Ag_0(\lam). 
\eqn

\subsection{The case $H$ has singularities of the first kind at zero} 
\bgth \lbth(R-first) 
Suppose $\ax^{\c}V\in L^1 (\R^2)$ for a $\c>4$ and $H$ has 
singularities of the first kind at zero. Then, $\W_{\low,2a}$ 
is a good operator.
\edth

\refth(R-first) together with \refth(high) 
proves that $W_{\pm}$ are good operators if 
$V \in \ax^{-2}L^{\frac43}(\R^2) \cap \ax^{-\c}V\in L^1 (\R^2)$ for a $\c>4$  
and if $H$ has singularities of the first kind at zero. This slightly improves 
Theorem 1.1 (i) of \cite{EGG} which  assumes 
$|V(x)|\leq C \ax^{-6-\ep}$, $\ep>0$.

\bgpf 
We prove $v\Ag_1(\lam)v$ is a good producer for small $\lam>0$. 
Then, \reflm(A+S-inv) and \refeq(1-s) imply the same for $v M(\lam)^{-1} v$ 
and \refth(R-first) follows.  We have $\rank S_1=1$ and 
$S_1 = \z \otimes \z$ for a normalized $\z \in S_1\HL$, 
$T_1 = c_3 \z\otimes \z $ with $c_3 = \|P T_0 \z \|^2>0$ 
and 
\bqn 
B_1(\lam)^{-1}= h_2(\lam)^{-1} (\z\otimes \z), \ \  %\\  
h_2(\lam)= -c_3 h_1(\lam) (1+ \Og_2(g(\lam)^2\lam^2))
\lbeq(B11)
\eqn 
from \refeq(B1-s). Then, \refeq(A+S1) and \refeq(gPT0) imply that 
modulo $\Og_2(g^2(\lam)\lam^2)$, 
\[
\Ag_1(\lam) \equiv N(\lam) S_1 B_1(\lam)^{-1} S_1 N(\lam)
= -c_3^{-1} h_1(\lam)^{-1} N(\lam)(\z\otimes \z) + \Og_2(g(\lam)^3 \lam^2). 
\]
Since $QD_0 Q\z = \z$, this simplifies to    
\[ 
\Ag_1(\lam) \equiv -c_3^{-1} (h_1(\lam)^{-1}\z\otimes \z + 
(L_1\z\otimes \z + \z\otimes L_1\z) + h_1(\lam)(L_1\z)\otimes (L_1\z))).
%\lbeq(A-1)
\]
Recall that $v(x)\z(x)\in L^1(\R^2)$ and $\int_{\R^2} v(x) \z(x) dx =0$. 
Thus, $v \Ag_1(\lam) v$ is a good producer by virtue of \refprop(funda) (7). 
\edpf

\subsection{The case $H$ has singularities of the second kind at zero} 

\bgth \lbth(p-wave) Suppose $\ax^{\c}V \in L^1(\R^2)$ for a $\c>8$ 
and $H$ has singularities of the second kind at zero. 
Then, for sufficently small $a>0$, 
${\W}_{{\rm low},2a}$ is bounded in $L^p(\R^2)$ for $1<p\leq 2$ but 
is unbounded in $L^p(\R^2)$ for $2<p<\infty$. 
\edth 
{\it We assume in the rest of the paper 
$\ax^{\c} V \in L^1(\R^2)$ for a $\c>8$}.  
The proof of \refth(p-wave) is long and is given by a series of lemmas. 
It is well known that $W_{\pm}$ are isometries of $\HL$ 
and we assume $p\not=2$.  Recall that \refeq(G1-inner) implies 
$1 \leq \rank T_2=\rank S_2\leq 2$ and $T_2$ is negative.
In what follows we assume $\rank S_2=2$.  Modification for the case 
$\rank S_2=1$ is obvious. Then, $n=\rank S_1$ is two or three 
depending on the absence or the presence of $s$-wave resonances.  
We take the orthonormal basis $\{\z_1,\dots, \z_n\}$ of $S_1\HL$ such that 
$\z_1,\z_2\in S_2\HL$ are (real) 
eigenfunctions of $T_2$: $T_2 \z_j= -\k_j^2 \z_j$ for $\k_j>0$, $j=1,2$.

\subsubsection{\bf Threshold analysis 2} 
We study $M(\lam)^{-1}$ as $\lam \to 0$ when $T_1$ is singular in $S_1\HL$ 
but $T_2=S_2(vG_1 v)S_2\vert_{S_2\HL}$ is non-singular. 
We consider $S_1$ and $S_2$ are projections also 
in $L^2(\R^2)$ as previously. Recall that we are omitting 
``for small $\lam>0$''.

\bglm  \lblm(S2-effect) The projection $S_2$ annihilates $T_0$ and $L_1$:    
\bqn \lbeq(T0S2)
T_0 S_2 = S_2 T_0 =0, \quad S_2 L_1= L_1 S_2=0.
\eqn 
\edlm
\bgpf 
Since $QT_0 QS_1=0$, $P T_0 QS_1= T_0 QS_1$,    
$T_1 = (T_0 Q S_1)^\ast (T_0 Q S_1)$ 
and $\Ker_{S_1\HL}T_1 = \Ker_{S_1\HL}T_0 Q S_1$.
Thus, $T_0Q S_1 S_2=0$ or $T_0 S_2=0$ and $S_2 T_0=(T_0 S_2)^\ast=0$. 
We have $PS_2 = PQS_2=0$ and likewise $S_2 P=0$;   
$D_0 S_2 = D_0 S_1 S_2 = S_2$ and likewise $S_2 D_0 = S_2$.  
It follows   
\[
S_2L_1= S_2(P - PT_0 QD_0 - D_0 QT_0 P + D_0 QT_0 P T_0 Q D_0)=0 
\]
and $L_1S_2 =(S_2 L)^\ast=0$.  
\edpf

Recall that $B_1(\lam)= -h_1(\lam)(T_1-\lam^2 X(\lam))$. 
We define for simplicity 
\bqn \lbeq(A1-def)
A_1(\lam)= T_1 - \lam^2 X(\lam) \ \ \mbox{so that} \  \ 
B_1(\lam)= - h_1(\lam) A_1(\lam)\,.
\eqn  
We study $A_1(\lam)^{-1}$ via \reflm(JN). Define  
\bqn 
\tRg_1(\lam) = S_2 \Rg_1(\lam) S_2= S_2 v(G_1+ g(\lam)^{-1} G_2)v S_2.
\lbeq(tRg-def)
\eqn 
and let $C(\lam)= (c_{jk}(\lam))$ be the representation matrix  of 
$\tRg_1(\lam)$  with respect to the basis $\{\z_1, \z_2\}$. We have 
\bqn 
c_{jk}(\lam)= -\k_j^2 \d_{jk} 
+ g(\lam)^{-1} \la G_2 v \z_j, v\z_k\ra \, \quad j,k=1,2. \lbeq(R1basis)
\eqn 
$C(\lam)$ is clearly invertible and we write $D(\lam)=(d_{jk}(\lam))$ 
for $C(\lam)^{-1}$:   
\bqn \lbeq(tRg-basis) 
S_2 \tRg_1(\lam)^{-1}S_2 = \sum_{j,k=1}^2 d_{jk}(\lam)\z_j \otimes \z_k. 
\eqn 
Evidently $d_{jk}(\lam)\chi_{\leq 2a}(\lam) \in \Mg(\R^2)$, 
$j,k=1,2$ for small $a>0$ and $T_1 + S_2$ is invertible in $S_1\HL$. 

\bglm  \lblm(B2-inv)
{\rm (1)} $A_1(\lam)+S_2$  on 
$S_1\HL$ is invertible and 
\bqn \lbeq(A1+S2) 
(A_1(\lam)+S_2)^{-1}=
(T_1 + S_2)^{-1}
+ \lam^2 (T_1 + S_2)^{-1}X(\lam) (T_1 + S_2)^{-1} + \Og_2 (g(\lam)^4 \lam^4).
\eqn 
With $\lam$-independent operators $F^{(k)}_j$ on $S_1\HL$ 
and modulo $\Og_2(g(\lam)^3\lam^4)$  
\begin{gather} 
(A_1(\lam)+ S_2)^{-1}S_2  \equiv S_2 + \lam^2 h_1(\lam)^{-1}g(\lam) 
(F_0^{(1)}+ h_1(\lam)F_1^{(1)}), \lbeq(A1-S2-commute) \\
S_2 (A_1(\lam)+ S_2)^{-1} \equiv   S_2 +  \lam^2 h_1(\lam)^{-1}g(\lam) 
(F_0^{(2)}+ h_1(\lam)F_1^{(2)}),  \lbeq(S2-A1-commute) \\
S_2 (A_1(\lam)+ S_2)^{-1}S_2 \equiv 
S_2 - \lam^2 h_1(\lam)^{-1}g(\lam) \tRg_1(\lam)\,.
\lbeq(S2-A1-S_2commute)
\end{gather}
{\rm (2)} Define $B_2(\lam)=S_2 - S_2 (A_1(\lam)+ S_2)^{-1} S_2$ 
on $S_2\HL$. Then, 
$B_2(\lam)$ is invertible and, as $\lam \to 0$, 
\bqn 
B_2(\lam)^{-1}=\lam^{-2}h_1(\lam)g(\lam)^{-1}\tRg_1(\lam)^{-1} 
+ S_2 \Og_2(1)S_2  + S_2 \Og_2(g(\lam)^2\lam^2) S_2. 
\lbeq(B2-inv)
\eqn 
\edlm  
\bgpf (1) Since  
$A_1(\lam)+ S_2= (T_1 + S_2)(1 -\lam^2 (T+S_2)^{-1}X(\lam))$ and 
$X(\lam)\in \Og_2(g(\lam)^2)$, 
$(A_1(\lam)+S_2)^{-1}$ exists and 
\refeq(A1+S2) is satisfied. 
Since $S_2 (T_1 + S_2)^{-1}=(T_1 + S_2)^{-1} S_2 = S_2$ 
and $S_2 L_1= L_1 S_2=0$, the second half of (1) follows from 
\refeq(A1+S2), \refeq(T0S2) and \refeq(X-def). 

(2) By virtue of \refeq(S2-A1-S_2commute), we have as $\lam \to 0$, 
\begin{align} 
B_2(\lam) & = \lam^2 h_1(\lam)^{-1}g(\lam) 
\{\tRg_1(\lam) + S_2 \Og_2 (g(\lam)^2\lam^2) S_2\} \notag \\
& =\lam^2 h_1(\lam)^{-1}g(\lam) 
(1+ S_2 \Og_2 (g(\lam)^2 \lam^2) S_2 \tRg_1(\lam)^{-1}) \tRg_1(\lam)\,.
\lbeq(B2)
\end{align} 
Thus, $B_2(\lam)$ is invertible and \refeq(X+1) 
implies \refeq(B2-inv). 
\edpf 

\bgrm \lbrm(B2-inv-a)
\ben 
\item[{\rm (a)}] $\tRg_1(\lam)^{-1}=\Og_2(1)$ is not assumed in statement 
{\rm (1)}. %\\ 
\item[{\rm (b)}] If $\tRg_1(\lam)^{-1}=\Og_2(g(\lam))$ as in \reflm(eigen-lm), 
\refeq(B2-inv) remains to holds if the last two terms 
are replaced by $S_2(\Og_2(g(\lam)^2)+ \Og_2(g(\lam)^5\lam^2) S_2$.  
\een
\edrm 

\bgprop \lbprop(impo-1-2) Modulo a good producer we have that   
\bqn \lbeq(ag3def)
vM(\lam)^{-1}v\equiv - g(\lam)^{-1} 
\lam^{-2}vS_2 \tRg_1(\lam)^{-1} S_2 v\,. 
\eqn  
\edprop 
\bgpf  Since $v\Ag_0(\lam)v$ is a good producer, it suffices to show 
\refeq(ag3def) with $\Ag_1(\lam)$ replacing $M(\lam)^{-1}$. 
Recall $\Ag_1(\lam){=} \Ag_0(\lam)S_1 B_1(\lam)^{-1}S_1 \Ag_0(\lam)$. 
\reflms(JN,B2-inv) imply  
\[%bqn \lbeq(A1-inv) 
A_1(\lam)^{-1}=(A_1(\lam)+ S_2)^{-1} + (A_1(\lam)+ S_2)^{-1} 
B_2(\lam)^{-1}(A_1(\lam)+ S_2)^{-1}. 
\]%eqn  
and substituting $- h_1(\lam)^{-1}A_1(\lam)^{-1}$ 
for $B_1(\lam)^{-1}$ yields  
$\Ag_1(\lam)= -(\Ag_{11}(\lam)+ \Ag_{12}(\lam))$, where omitting the 
variable $\lam$ of $B_2^{-1}(\lam)$,  
\begin{align}  
\Ag_{11}(\lam)& =h_1(\lam)^{-1}\Ag_0(\lam) 
S_1(A_1(\lam)+ S_2)^{-1}S_1 \Ag_0(\lam), \lbeq(Ag11) \\ 
\Ag_{12}(\lam)& =
h_1(\lam)^{-1}\Ag_0(\lam)S_1 (A_1(\lam)+ S_2)^{-1}  
S_2 B_2^{-1}S_2(A_1(\lam)+ S_2)^{-1}S_1\Ag_0(\lam) . \lbeq(Ag12)
\end{align}  
We first show that $v\Ag_{11}(\lam)v$ is a good producer. 
The following lemma also does not assume $\tRg_1(\lam)^{-1}=\Og_2(1)$. 
Recall $N(\lam)= h_1(\lam) L_1 + Q D_0 Q$ (see \refeq(gPT0)). 

\begin{lemma} \lblm(first-line)
{\rm (1)} We have that  
\bqn \lbeq(A-S1-commute) 
\left\{ \br{rl} 
\Ag_0(\lam)S_1 & = S_1 + h_1(\lam) L_1 S_1 + \Og_2(\lam^2 g(\lam) )S_1, \\ 
S_1\Ag_0(\lam) & = S_1 + h_1(\lam) S_1 L_1 + S_1 \Og_2(\lam^2 g(\lam))
\er \right. 
\eqn 
and 
\bqn 
\lbeq(A+S-S2)
\left\{ \br{rl} 
\Ag_0(\lam)S_2 & 
= S_2 + g(\lam)\lam^2 N(\lam)\Rg_1(\lam) S_2 +\Og_2(\lam^4 g(\lam)^2), \\
S_2 \Ag_0(\lam) & = S_2 + 
g(\lam)\lam^2 S_2 \Rg_1(\lam)N(\lam)+\Og_2(\lam^4 g(\lam)^2).
\er
\right. 
\eqn 
{\rm (2)} $v\Ag_{11}(\lam)v$ is a good producer.  
\edlm 

\bgpf (1) If we multiply \refeq(A+S1) by $S_1$ from the right or the left 
and use that $S_1 QD_0 Q= QD_0 Q S_1= S_1$, \refeq(A-S1-commute) follows. 
If we multiply \refeq(A+S1-a) by $S_2$ from the right or from the left and 
apply \refeq(T0S2), we obtain \refeq(A+S-S2).  \\
(2) Since $X(\lam) = \Og_2(g(\lam)^2)$, we have from \refeq(A1+S2) that 
\[
S_1 (A_1(\lam)+ S_2)^{-1}S_1 = S_1(T_1+S_2)^{-1} S_1 + \Og_2(g(\lam)^2 \lam^2).
\]
Then, \refeq(Ag11) and \refeq(A-S1-commute) imply 
that modulo $\Og_2(g(\lam)^3 \lam^2)$ 
\begin{align}
\Ag_{11}(\lam)& \equiv h_1(\lam)^{-1}(S_1 + h_1(\lam) L_1 S_1)(T_1+ S_2)^{-1}
( S_1 + h_1(\lam) S_1 L_1)   \notag \\
& = h_1(\lam)^{-1} S_1(T_1+ S_2)^{-1}S_1 
+ F_3 + h_1(\lam) F_4 \lbeq(ag11)
\end{align} 
with $\lam$-independent finite rank operators $F_3$ and $F_4$ 
and, \refprop(funda) implies 
$v \Ag_{11} v$ is a good producer. Note that  
$S_1(T_1+ S_2)^{-1}S_1= \sum_{j,k=1}^n a_{jk}(\z_j \otimes \z_k)$ 
with $z_1, \dots, z_n \in S_1\HL$, hence  
$v \z_j\in L^1(\R^2)$ , $\la \z_j, v \ra =0$ for $j=1, \dots, n$. 
\edpf 

Plugging \refeq(B2-inv) and \refeq(Ag12) implies   
$\Ag_{12}(\lam)\equiv \Ag_{12}^{(1)}(\lam) + \Ag_{12}^{(2)}(\lam)$ 
modulo $\Og_2(g(\lam)^3 \lam^2)$: 
\begin{align} 
\Ag_{12}^{(1)}(\lam)& = 
g(\lam)^{-1} \lam^{-2}\Ag_0(\lam)S_1 (A_1(\lam)+ S_2)^{-1}S_2 
\lbeq(impo-10) \\
& \times 
\tRg_1(\lam)^{-1}S_2(A_1(\lam)+ S_2)^{-1}S_1\Ag_0(\lam) , \notag 
\\
\Ag_{12}^{(2)}(\lam) 
& = \Ag_0(\lam)S_1 (A_1(\lam)+ S_2)^{-1}S_2 \Og_2(g(\lam))S_2 
\lbeq(impo-11) \\
& \times 
(A_1(\lam)+ S_2)^{-1}S_1\Ag_0(\lam) +\Og_2(g(\lam)^3 \lam^2).   \notag 
\end{align}
We first show 
\bglm \lblm(impo-1-1) 
$v\Ag_{12}^{(2)}(\lam)v$ is a good producer for small 
$\lam>0$. 
\edlm 
\bgpf  $\Ag_{12}^{(2)}(\lam)=
\Ag_0(\lam)S_2 \Og_2(g(\lam))S_2 \Ag_0(\lam) 
+\Og_2(\lam^2g(\lam)^3)$ 
by virtue of \refeq(A1-S2-commute) and  \refeq(S2-A1-commute) 
and, further applying \refeq(A-S1-commute) and \refeq(T0S2) yields  
that $\Ag_{12}^{(2)}(\lam)= S_2 \Og_2(g(\lam))S_2  + \Og_2(\lam^2g(\lam)^3)$. 
In the basis $\{\z_1, \z_2\}$ of $S_2\HL$, we have that 
\bqn \lbeq(s222s)
S_2 \Og_2(g(\lam))S_2 = g(\lam) \sum b_{jk}(\lam) \z_j \otimes \z_k ,
\quad b_{jk}(\lam)\in \Mg(\R^2).
\eqn 
\refprop(funda)(7) then implies that $v\Ag_{12}^{(2)}(\lam)v$ as previously. 
\edpf 

\bgrm \lbrm(A12-good)  
\reflm(impo-1-1) holds if $\tRg_1(\lam)=\Og_2(g(\lam))$ as in 
\reflm(eigen-lm) because \refeq(tRg-basis) remains to hold 
with $d_{jk}(\lam)$ being replaced by $g(\lam)\tilde{d}_{jk}(\lam)$,  
$\tilde{d}_{jk}(\lam) \in \Mg(\R^2)$ (see \refeq(dd) below) and, by 
virtue of \refrm(B2-inv-a), 
this  will only change  $\Og_2(g(\lam))$ to $\Og_2(g(\lam)^3)$ 
and $\Og_2(\lam^2g(\lam)^3)$ to $\Og_2(\lam^2g(\lam)^5)$ 
in the proof above and \refprop(funda) {\rm (7)} applies without any change. 
\edrm 

\paragraph{\bf Completion of the proof of \refpropb(impo-1-2)} 
We have only to show that $v\Ag_{12}^{(1)}(\lam)v \equiv 
- g(\lam)^{-1} \lam^{-2}vS_2 \tRg_1(\lam)^{-1} S_2 v$ 
modulo a good producer. Define  
\bqn \lbeq(Pg)
\Pg(\lam)= \lam^{-2}g(\lam)^{-1} (A_1(\lam)+ S_2)^{-1}S_2 
\tRg_1(\lam)^{-1}S_2(A_1(\lam)+ S_2)^{-1} 
\eqn 
so that 
\bqn \lbeq(Ag121)
\Ag_{12}^{(1)}(\lam) = 
\Ag_0(\lam)S_1 \Pg(\lam)S_1\Ag_0(\lam) .
\eqn 
We substitute $(A_1(\lam)+ S_2)^{-1}S_2$ and $S_2(A_1(\lam)+ S_2)^{-1}$ 
by right sides of \refeq(A1-S2-commute) and \refeq(S2-A1-commute) 
respectively and expand the sum. Then,    
\begin{gather*}
\Pg(\lam)=\lam^{-2}g(\lam)^{-1}S_2\tRg_1(\lam)^{-1}S_2  
+ h_1(\lam)^{-1} L_3(\lam) + \Og_2(g(\lam)^3\lam^2), %\lbeq(center) 
\\  
L_3(\lam)= 
 (F_{0}^{(1)}+h_1(\lam)F_{1}^{(1)})\tRg_1(\lam)^{-1}  + 
\tRg_1(\lam)^{-1}(F_{0}^{(2)}+ h_1(\lam) F_{k}^{(2)}), %\lbeq(L3)
\end{gather*} 
which we insert for $\Pg(\lam)$ of \refeq(Ag121). This 
produces three operators for $\Ag_{12}^{(1)}(\lam)$.  After sandwitched by $v$, 
the one produced by $\Og_2(g(\lam)^3\lam^2)$ is a good producer 
by \refprop(funda) (5). The definition of $L_3(\lam)$ 
and \refeq(A-S1-commute) imply  
\begin{multline}  \lbeq(middle)
\Ag_0(\lam)S_1(h^{-1} L_3(\lam))S_1\Ag_0(\lam) \\
= h^{-1} L_3(\lam)+ L_1S_1 L_3(\lam)+ L_3(\lam)S_1 L_1  + \Og_2(g(\lam)^2\lam^2).
\end{multline}
Again $v\Og_2(g(\lam)^2\lam^2)v$ is a good producer 
and the same is true for first three terms by virtue of \refprop(funda). 
This is because \refeq(tRg-basis) implies 
$L_3(\lam)= \sum_{j,k=1}^2 b_{jk}(\lam) \z_j \otimes \z_k$ 
with good multipliers $b_{jk}(\lam)$ and they are the sum of the 
operators which appear in \refprop(funda) (7). 
By virtue of \refeq(A+S-S2) the one produced by 
$g(\lam)^{-1} \lam^{-2}S_2 \tRg_1(\lam)^{-1}S_2 $ may be 
expressed in the form 
\begin{multline} \lbeq(Df)
g(\lam)^{-1} \lam^{-2} S_2 \tRg_1(\lam) S_2
+ N(\lam)\Rg_1(\lam)S_2 \tRg_1(\lam)^{-1}S_2 \\
+ S_2 \tRg_1(\lam)^{-1}S_2\Rg_1(\lam) N(\lam) 
+ \Og_2(g(\lam)\lam^2).
\end{multline}
Here the last three terms are good producers by the same reason 
as above and $vM(\lam)v \equiv 
 g(\lam)^{-1} \lam^{-2} v S_2 \tRg_1(\lam)^{-1}S_2v$. 
This proves the proposition 
\edpf 

\bgrm \lbrm(prop-517) \refprop(impo-1-2) holds if 
$\tRg_1(\lam)=\Og_2(g(\lam))$ as in \reflm(eigen-lm) 
because it only changes $d_{jk}(\lam)$ in \refeq(tRg-basis) 
by $g(\lam)\tilde{d}_{jk}(\lam)$ with possibly another 
$\tilde{d}_{jk}(\lam) \in \Mg(\R^2)$ and because 
and the power $k$ of $g(\lam)^k$ is arbitrary in \refprop(funda) (7). 
\edrm 

\subsubsection{\bf Good and bad parts.} 
By virtue of \refprop(impo-1-2) the proof of \refth(p-wave) 
is reduced to proving that 
\bqn \lbeq(tWlow-a-1) 
\tilde{\W}_{{\rm low},2a}u=
- \int_0^\infty g(\lam)^{-1}\lam^{-2} G_0(-\lam) vS_2 \tRg_1(\lam)^{-1}
S_2v \Pi(\lam)u {\lam}\chi_{\leq 2a}(\lam) d\lam
\eqn   
defined by \refeq(tWlow-a) 
is bounded in $L^p(\R^2)$ for $1<p\leq 2$ and unbounded for $2<p<\infty$. 
Substituting \refeq(tRg-basis) for $S_2 \tRg_1(\lam)^{-1}S_2$ shows that 
$\tilde{\W}_{{\rm low},2a}u$ is the sum over $j,k=1,2$ of 
\bqn  
\tilde{\W}_{{\rm low},2a}^{(j,k)}u= 
- \int_0^\infty g(\lam)^{-1}\lam^{-2} d_{jk}(\lam)
G_0(-\lam)|v\z_j\ra \la \z_k v, \Pi(\lam)u\ra {\lam}\chi_{\leq 2a}(\lam) d\lam\,.
\lbeq(tWlow-1)
\eqn  
Recall that for $u \in \Dg_\ast$,
$\Pi(\lam) u=0$ for $\lam \not\in (\a,\b)$ for $0<\a<\b<\infty$ 

Since $\la v, \z_j\ra =0$, $j=1,2$, we may replace $\Pi(\lam)u(z)$  
by $\Pi(\lam)u(z)-\Pi(\lam)u(0)$ which we decompose into the sum of 
the good part \refeq(good-def-intro) and the bad part \refeq(bad-def-intro):   
\bqn \lbeq(g+b)
\Pi(\lam)u(z)-\Pi(\lam)u(0)= \tilde{g}(\lam,z)+ \tilde{b} (\lam,z).
\eqn 
Define ${\tilde{\W}}_{(g)}^{(j,k)}$ and 
${\tilde{\W}}_{(b)}^{(j,k)}$ 
by \refeq(tWlow-1) by replacing $\Pi(\lam)u(z)$ by   
$\tilde{g}(\lam,z)$ and $\tilde{b} (\lam,z)$ respectively and 
\bqn 
\tilde{\W}_{(g)}u(x)
{=} \sum_{j,k=1}^2 {\tilde{\W}}_{(g)}^{(j,k)}u(x), \quad  
\tilde{\W}_{(b)}u(x)= \sum_{j,k=1}^2  
{\tilde{\W}}_{(b)}^{(j,k)}u(x)
\eqn 
so that $\tilde{\W}_{{\rm low},2a}u(x)= \tilde{\W}_{(g)}u(x)
+ \tilde{\W}_{(b)}u(x)$. We recall 
\begin{align}
& \tilde{g}(\lam,z)= - \sum_{l,m=1}^2 z_l z_m 
\lam^2 \int_0^1 (1-\th) \left(\frac1{2\pi}
\int_{{\mathbb S}^1} \Fg(\tau_{-\th{z}}R_l R_m u)(\lam\w)d\w \right) 
d\th. \notag \\% \lb \\
& \tilde{b}(\lam,z)= 
\frac{i\lam}{2\pi}\int_{{\mathbb S}^1} 
({z}\w)\hat{u}(\lam{\w})d\w
=\sum_{l=1}^2 \frac{i\lam{z_l} }{2\pi} \int_{{\mathbb S}^1} 
(\Fg R_l{u})(\lam{\w})d\w\,.  \lbeq(bad-def)
\end{align} 

\bgprop \lbprop(good) Te good part 
${\tilde{\W}}_{(g)}u(x)$ is a good operator. 
\edprop

\bgpf Define $\mu_{jk}(\lam) = d_{jk}(\lam)g(\lam)^{-1}\chi_{\leq 2a}(\lam)$, 
$j,k=1,2$. $\mu_{jk}(\lam)$ is a good multiplier. Then, after 
changing the order of integrations, we can express 
${\tilde{\W}}_{(g)}^{(j,k)}u(x)$ as 
\begin{multline*} %\lbeq(good-1)
{\tilde{\W}}_{(g)}^{(j,k)}u(x)=\sum_{l,m=1}^2 \int_0^1 (1-\th) d\th
\int_{\R^4}dy dz\, v(y)\z_j(y) v(z)\z_k(z)z_l z_m \\
\times \left\{ \int_0^\infty  \m_{jk}(\lam) \Gg_{-\lam}(x-y)  
\left( \frac{1}{2\pi}\int_{{\mathbb S}^1}
\Fg(\t_{-\th{z}}R_l R_m u)(\lam\w)d\w \right)\lam d\lam\right\} .
\end{multline*}
The identity \refeq(mult) and the definition \refeq(DFEK) of $K$ 
imply that the second line is equal to  
\begin{align}
& \tau_y \int_0^\infty \Gg_{-\lam}(x)  
\left( \frac{1}{2\pi}\int_{{\mathbb S}^1}
\Fg(\m_{jk}(|D|)\t_{-\th{z}}R_l R_m u)(\lam\w)d\w \right)\lam d\lam \notag \\
& \qquad = (\t_y K \m_{jk}(|D|)\t_{-\th{z}}R_l R_m u)(x). %\lbeq(good-1a)
\end{align}
Since $\|\t_y K \m_{jk}(|D|)\t_{-\th{z}}R_l R_m u \|_p \leq C\|u\|_p$ with 
$C$ independent of $y,z\in \R^2$, $0<\th<1$ and $u\in \Dg_\ast$, 
Minkowki's inequality implies 
\bqn \lbeq(good-2)
\|{\tilde{\W}}_{(g)}^{(j,k)}u\|_p  \leq  C_p \left(\sum  
\|z_l z_m v\z_k\|_1 \|v\z_j\|_1 \right)\|u\|_p . 
\eqn 
Since $p$-wave resonances $\z_j$ satisfies $\ax^{-\d}\z_j \in \HL$ 
for any $0<\d<1$, Schwarz inequality implies 
$\|\la z \ra^2 \z_j v\|_1^2 \leq C\|\ax^{4+\d} V \|_1$. 
Thus, ${\tilde{\W}}_{(g)}$ is a good operator. 
\edpf

\subsubsection{\bf Estimate of bad part 1, Positive result} 

The following proposition completes the proof of the positive part of 
\refth(p-wave). 
\bgprop \lbprop(low-b)
${\tilde{\W}}_{(b)}$ is bounded in $L^p(\R^2)$ for $1<p\leq 2$.
\edprop 
\bgpf  Let 
$\m_{jk}(\lam)= g(\lam)^{-1}d_{jk}(\lam)\chi_{\leq 2a}(\lam)\in \Mg(\R^2)$ 
for $j,k=1,2$ as previously. We have 
\bqn 
{\tilde{\W}}_{(b)}^{(j,k)}u(x) 
= - \int_0^\infty 
\lam^{-1} \mu_{jk}(\lam) G_0(-\lam) (v\z_j)(x)
\la \z_k, \tilde{b}(\lam)\ra  d\lam.  \lbeq(tWlow-b)
\eqn 
Using $\chi_{\leq 4a}(\lam) + \chi_{>4a}(\lam)=1$, we decompose 
${\tilde{\W}}_{(b)}^{(j,k)}u(x)$ as follows:  
\[
{\tilde{\W}}_{(b)}^{(j,k)}u 
= \chi_{>4a}(|D|){\tilde{\W}}_{(b)}^{(j,k)}u 
+\chi_{\leq 4a}(|D|){\tilde{\W}}_{(b)}^{(j,k)}u .
\]
Then, the following two lemmas prove the proposition.  
\edpf 

\bglm \lblm(low-high-positive)  For $1<p\leq 2$,   
$\chi_{>4a}(|D|){\tilde{\W}}_{(b)}\in \Bb(L^p(\R^2))$. 
\edlm 
\bgpf For $j,k,l=1,2$, define $X_{jkl}u(x)$ by 
\bqn%egin{multline} 
\int_0^\infty 
\chi_{>4a}(|D|)G_0(-\lam)(v \z_j)(x)\left(
\int_{{\mathbb S}^1}(\Fg \m_{jk}(|D|)R_l u)(\lam\w)d\w\right)
d\lam. \lbeq(6-0)
\eqn
Substituting \refeq(bad-def) for $\tilde{b}(\lam,z)$ 
in \refeq(tWlow-b), we obtain 
\bqn \lbeq(6-0-)
\chi_{>4a}(|D|){\tilde{\W}}_{(b)}^{(j,k)}u(x)
=- \sum_{l=1}^2 \frac{i}{2\pi}\la z_l v, \z_k\ra X_{jkl}u(x) 
\eqn 
and we prove $X_{jkl}u(x)$, $j,k,l=1,2$ are bounded in $L^p(\R^2)$ 
for $1<p\leq 2$. Let $\m(\xi)= \chi_{>4a}(\xi)|\xi|^{-2}$ as in the 
proof of \reflm(swave-r).  Then, we have (see \refeq(6-2))
\bqn \lbeq(6-1)
\chi_{>4a}(|D|)G_0(-\lam)v \z_j(x)=
\frac{1}{2\pi}\hat{\m}\ast (v\z_j)(x)   
+ \lam^2 \int_{\R^2}\mu(|D|)\Gg_{-\lam}(x-y)(v \z_j)(y) dy.  
\eqn 
Plugging \refeq(6-1) and \refeq(6-0) yields   
\bqn \lbeq(6-01)
X_{jkl}u(x) = X_{jkl}^{(b)}u(x) + X_{jkl}^{(g)}u(x), 
\eqn 
where $X_{jkl}^{(b)}u(x)$ and $X_{jkl}^{(g)}u(x)$ 
are the functions produced by 
the first and the second terms on the right of of \refeq(6-1) respectively. 
Denote ${\r}_{jk}(\lam) = \lam \mu_{jk}(\lam)$. Then, 
by integrating with respect to $d\lam$ first and by recalling \refeq(DFEK), 
we obtain that $X_{jkl}^{(g)}u(x)$ is equal to   
\begin{align*}
& \int_{\R^2} (v \z_j)(y) 
\left\{ 
\int_0^\infty 
\mu(|D|) \t_y \Gg_{-\lam}(x) 
\left(\int_{{\mathbb S}^1}
(\Fg \r_{jk}(|D|)R_lu)(\lam)d\w\right)\lam d\lam \right\}dy \\
& = 2\pi\int_{\R^2} (v \z_j)(y) 
(\mu(|D|)K {\r}_{jk}(|D|)R_lu)(x-y)dy. 
\end{align*}
Minkowski's inequality, \refeq(K-est) and the multiplier theory 
then imply 
\bqn \lbeq(Xgood)
\|X_{jkl}^{(g)}u\|_p \leq C_p \|v \z_j\|_1 \|u\|_p , \quad 1<p<\infty.
\eqn 
and  $X_{jkl}^{(g)}u(x)$ is a good operator, $1\leq j,k,l\leq 2$. 

By observing that $\hat{\m}\ast (v\z_j)(x)$ is independent of $\lam$ 
and by using the polar coordinates $\xi=\lam\w$, 
we represent 
\begin{align} 
X_{jkl}^{(b)}u(x)& = \frac{1}{2\pi} (\hat{\m} \ast v\z_j)(x) 
\int_0^\infty \left(\int_{{\mathbb S}^1}(\Fg \m_{jk}(|D|)R_l u)(\lam)d\w\right)
d\lam \notag \\
& =\frac{1}{2\pi} (\hat{\m} \ast v\z_j)(x) 
\int_{\R^2} (\Fg \m_{jk}(|D|)R_l u)(\xi)|\xi|^{-1}d\xi .  
\lbeq(Xbad)
\end{align}
Here $\hat{\m} \ast v\z_j\in L^p(\R^2)$ for any $1<p<\infty$ since 
so is $\hat{\m}$  and $v\z_j \in L^1(\R^2)$. Since 
$\m_{jk}(\lam)|\xi|^{-1}\in L^p(\R^2)$ for $1\leq p \leq 2$, 
Hausdorff-Young's inequality implies that  
\bqn \lbeq(Xbad-a)
\int_{\R^2} (\Fg \m_{jk}(|D|)R_l u)(\xi)|\xi|^{-1}d\xi
= \int_{\R^2} u(x)
\Fg(\m_{jk}(|\xi|)\xi_l|\xi|^{-2})(x)dx
\eqn 
are bounded linear functionals on $L^p(\R^2)$ for $1\leq p\leq 2$. 
Thus, $X_{jkl}^{(b)}$ are bounded in $L^p(\R^2)$ for $1< p \leq 2$. 
This together with \refeq(Xgood) proves the lemma. 
\edpf 

\bglm \lblm(low-low) Let $1<p<2$. Then,  
$\chi_{\leq 4a}(|D|){\tilde{\W}}_{(b)}$ is bounded 
from $L^p(\R^2)$ to itself and, hence the same holds for 
${\tilde{\W}}_{(b)}$.
\edlm 
\bgpf The proof is a slight modification of that of 
\reflm(swave-l) and we shall be a little sketchy in some places. 
As in the proof of 
\reflm(low-high-positive) it suffices to show the lemma for   
$\tilde{X}_{jkl}u(x)$, $j,k,l=1,2$ defined by replacing 
$\chi_{>4a}(|D|)$ by $\chi_{\leq 4a}(|D|)$ in \refeq(6-0): 
\bqn \lbeq(Xjkl)
\tilde{X}_{jkl}u(x)=
\int_0^\infty 
\chi_{\leq 4a}(|D|)G_0(-\lam)v \z_j(x) 
 \left(\int_{{\mathbb S}^1}(\Fg \m_{jk}(|D|)R_l u)(\lam\w)d\w\right)
d\lam \,.
\eqn  
We denote $u_{jkl}= \m_{jk}(|D|)R_l u$. We use \refeq(6-20-a) for 
$\p= v \z_j$ that 
\bqn \lbeq(6-20)
\chi_{\leq 4a}(|D|)G_0(-\lam)v \z_j(x)
 = \frac1{2\pi} \int_{\R^2} 
\frac{e^{ix\xi}\chi_{\leq 4a}(|\xi|)\Fg(v\z_j)(\xi)}
{\xi^2 -(-\lam +i0)^2}d\xi, 
\eqn 
which, by virtue of \refeq(exp-Fp) adapted for $v\z_j(x)$, is equal to   
\[%begin{multline} 
\frac{-i}{2\pi}\sum_{m=1}^2 \int_0^1 
\int_{\R^2}z_m(v\z_j)(z) \tau_{{\th}z}
\left(\frac1{2\pi}\int_{\R^2}e^{ix\xi}\frac{\xi_m\chi_{\leq 4a}(|\xi|)}
{\xi^2 -(-\lam +i0)^2} d\xi \right) dz d\th 
\]
(see \refeq(6-21-a)). The inner integral is computed in \refeq(rmchi-a) 
and is equal to 
\bqn \lbeq(rmchi)
2\pi R_m \lam \chi_{\leq 4a}(|D|)\Gg_{-\lam}(x) + R_m\left(\frac1{2\pi} 
\int_{\R^2} e^{ix\xi} \frac{\chi_{\leq 4a}(|\xi|)}{|\xi|+\lam}d\xi\right)\,.
\eqn
The contribution of the first term for $\tilde{X}_{jkl}u(x)$ is  
given as in \refeq(40) by 
\[
-2\pi{i}\sum_{m=1}^2 \int_0^1 d\th\int_{\R^2}z_m(v\z_j)(z) 
(R_m \chi_{\leq 4a}(|D|)K{u_{jkl}})(x-\th{z})dz d\th %\lbeq(6-22)
\]
and is a good operator being bounded by 
\bqn \lbeq(40-good)
C \int_{\R^2}|z_m(v\z_j)(z)| \|u_{jkl}\|_p \leq C 
\|z_m v(z)\z_j(z)\|_1 \|u\|_p\,, \quad 1<p<\infty.
\eqn 

The contribution of the second term of \refeq(rmchi) for 
$\tilde{X}_{jkl}u(x)$ is given after changing the order of integrations 
by the $\sum_{l=m}^2 \int_0^1 d\th$ of 
\bqn
\frac{-iR_m}{(2\pi)^2}
\int_{\R^2}z_m(v\z_j)(z) \tau_{{\th}z} \left(\int_{\R^2} 
\frac{e^{ix\xi}\chi_{\leq 4a}(|\xi|)}{(|\xi|+|\eta|)|\eta|}
\widehat{u_{jkl}}(\eta)d\xi d\eta \right)dz.  \lbeq(6-23)
\eqn 
If we define the integral operator $L_1$ by the integral kernel 
\[
L_1 (x,y)= \int_{\R^4} 
\frac{e^{ix\xi-iy\eta}\chi_{\leq 4a}(|\xi|)\chi_{\leq 2a}(|\eta|)}
{(|\xi|+|\eta|)|\eta|g(|\eta|)}d\xi d\eta 
\]
then, the inner integral of \refeq(6-23) becomes 
$L_1 d_{jk}(|D|)R_l u(x)$ (recall that $u_{jkl}= \m_{jk}(|D|)R_l u$)
and  
\[
\refeq(6-23)= \frac{-iR_m}{(2\pi)^2}
\int_{\R^2}z_m(v\z_j)(z) (L d_{jk}(|D|)R_m u)(x-\th{z}) dz \,. 
\] 
We have shown in Lemma 4.7 of \cite{Ya-point} that $L_1$ is a 
bounded operator of $L^p(\R^2)$ for $1<p\leq 2$. Thus, Minkowski's inequality 
implies 
\[
\|\refeq(6-23)\|_p \leq C\|z_m(v\z_j)(z)\|_1 \|u\|_p, \quad 1<p \leq 2.
\]
This to gether with \refeq(40-good) proves \reflm(low-low). \edpf 

\subsubsection{\bf Estimate of bad part 2, Negative result} 

\bgprop \lbprop(low-high-negative)  
$\chi_{>4a}(|D|){\tilde{\W}}_{(b)}$ is unbounded in $L^p(\R^2)$ 
for $2<p<\infty$ if $a>0$ is sufficiently small.  
\edprop
Since $\chi_{>4a}(|D|)$ is a good operator, the proposition implies 
${\tilde{\W}}_{(b)}$ is unbounded in $L^p(\R^2)$ 
for $2<p<\infty$ for small $a>0$. 
It follow, since ${\tilde{\W}}_{(g)}$ is a good operator for any $a>0$ 
by virtue of \refprop(good), $W_{+}$ is unbounded in $L^p(\R^2)$ 
for $2<p<\infty$ if $H$ has a singularity of the second kind at zero. 

\bgpf By virtue of \refeq(6-0-), \refeq(6-01) and \refeq(Xgood), 
it suffices to show that the operator 
$\tilde{\W}_{\rm low, \ast}u(x)\stackrel{\rm def}{=} 
\sum_{j,k,l=1}^2 i\la z_l v | \z_k\ra X_{jkl}^{(b)} u (x)$ 
is unbounded  in $L^p(\R^2)$ if $2<p<\infty$ and if $a>0$ is small enough. 
By virtue of \refeqs(Xbad,Xbad-a) 
\bqn 
\tilde{\W}_{\rm low, \ast}u(x) 
=\sum_{j=1}^2 (\hat{\m} \ast(v\z_j))(x)\ell_j(u)
\eqn 
where $\ell_1(u)$ and $\ell_2(u)$ are linear functionals defined by 
\bqn
\ell_j(u) = \sum_{k,l=1}^2  i\la z_l v | \z_k\ra \int_{\R^2} u(x)
\Fg(\m_{jk}(|\xi|)\xi_l|\xi|^{-2})(x)dx\, \ j=1,2\,. \lbeq(6-79)
\eqn 
It is obvious that 
$\hat{\m} \ast(v\z_1)$ and $\hat{\m} \ast(v\z_2)\in L^p(\R^2)$  
for $2<p<\infty$ and, for small $a>0$,
they are linearly independent 
in $L^p(\R^2)$ as will be shown in the next subsection where 
we prove a more general result. Then, if $\tilde{\W}_{\rm low\ast}$ 
were bounded in $L^p(\R^2)$ for $2<p<\infty$, it must be that 
both $\ell_1$ and $\ell_2$ were continuous functionals on $L^p(\R^2)$ 
by the Hahn-Banach theorem, that for $q=p/(p-1)$, $1<q<2$
\[
\sum_{k,l=1}^2 \la z_l v | \z_k\ra 
\Fg(\m_{jk}(|\xi|)\xi_l|\xi|^{-2}) \in L^q(\R^2) , \quad j=1,2
\]
by the Riesz representation theorem and that  
\bqn \lbeq(Hau)
\sum_{k=1}^2  d_{jk}(|\xi|) 
\sum_{l=1}^2 \la z_l v | \z_k\ra \chi_{\leq 2a}(\xi)\xi_l 
g(|\xi|)^{-1}|\xi|^{-2} \in L^p(\R^2), 
\eqn 
by Haudorff-Young's inequality where we have restored 
$g(\lam)^{-1}d_{jk}(\lam)\chi_{\leq 2a}(\lam)$ for 
$\m_{jk}(\lam)$. 
\refeq(Hau) means in the matrix notation that 
\bqn  \lbeq(6-100) 
\frac{\chi_{\leq 2a}(|\xi|)}{g(|\xi|)|\xi|^{2}}
D(|\xi|)\begin{pmatrix}
\la z_1 v | \z_1\ra \xi_1+ \la z_2 v | \z_1\ra \xi_2 \\ 
\la z_1 v | \z_2 \ra \xi_1+ \la z_2 v | \z_2\ra 
\xi_2  \end{pmatrix} \in L^p(\R^2, \C^2), 
\eqn 
Then, since $D(\lam)=C(\lam)^{-1}$, \refeq(R1basis) and \refeq(6-100) 
imply that   
\bqn \lbeq(6-100a)
\frac{\chi_{\leq 2a}(|\xi|)}{g(|\xi|)|\xi|^{2}}
(\la z_1 v | \z_k\ra \xi_1+ \la z_2 v | \z_k\ra \xi_2) \in L^p(\R^2), 
\ \ k=1,2\,.
\eqn 
But for $p>2$ this can happen only when 
$\la z_j v | \z_k\ra= 0$ for $j,k=1,2$. 
However, $-\k_j^2 =\la v G_1 v\z_k, \z_k \ra 
= -\frac12 \sum_{j=1}^2 \left|\la z_j v, \z_k \ra \right|^2<0$, $k=1,2$
and this is a contradiction. Thus, 
$\tilde{\W}_{\rm low\ast}$ must be unbounded in $L^p(\R^2)$ for 
any $2<p<\infty$ if $a>0$ is small enough. 
\edpf 

\subsection{The case $H$ has sigularities of the third kind at zero}
In this case $T_2$ is singular in $S_2 \HL$;   
$S_3$ is the projection in $S_2\HL$ onto ${\rm Ker}\, T_2$ 
and $T_3= S_3 v G_2 v S_3$ is non-singular in $S_3 \HL$. 
We first assume that $S_2 \ominus S_3 \not=0$. As previously 
we shall omit ``for small $\lam>0$''.  

\bgth \lbth(third-1) 
Suppose that $H$ has singularities of the third kind at zero 
and $S_3 \subsetneq S_2$. Then $W_{+}$ is bounded in $L^p(\R^2)$ 
for $1<p\leq 2$ and is unbounded in $L^p(\R^2)$ for $2<p<\infty$. 
\edth 

\subsubsection{\bf Threshold anaysis 3} 
For shortening formulas, define 
\[
T= S_2 v G_1 v S_2 (=T_2), \ \ \tT= S_2 v G_2 v S_2, \ 
\]
so that ${\tRg}_1(\lam) = T + g(\lam)^{-1}\tT$ and $T_3= S_3 \tT S_3$. 
We define $\Xg_2 = (S_2 \ominus S_3)\HL$ and $\Xg_3 = S_3 \HL$ and 
denote by $P_2$ and $P_3$ the projections in $S_2\HL$ onto $\Xg_2$ 
and $\Xg_3$ respectively. We express $\tRg_1(\lam)$ in $S_2\HL$ 
as the operator matrix in the decomposition $S_2\HL= \Xg_2\oplus \Xg_3$: 
\bqn \lbeq(FS-Rg)
\tRg_1(\lam) = 
\begin{pmatrix} T_{22} + g(\lam)^{-1}\tT_{22} &  g(\lam)^{-1}\tT_{23} \\
g(\lam)^{-1}\tT_{32} & g(\lam)^{-1}\tT_{33} 
\end{pmatrix}
\eqn 
where $T_{jk}= P_j T P_k$ and etc. and we have used that 
$T_{23}=0, T_{32}=0$ and $T_{33}=0$. 
Note that $T_{22}$ and $\tilde{T}_{33}=T_3$ are invertible 
in $\Xg_2$ and in $\Xg_3$ respectively. It follows that   
$T_{22} + g(\lam)^{-1}(\tilde{T}_{22} - \tT_{23}\tT_{33}^{-1}\tT_{32})$ 
is invertible in $\Xg_2$ for small $\lam>0$ and we define 
\bqn 
\tilde{d}(\lam) \stackrel{\rm def}{=} T_{22}^{-1}
(1_{\Xg_2}+g(\lam)^{-1}(\tT_{22} - 
\tT_{23}\tT_{33}^{-1}\tT_{32})T_{22}^{-1} )^{-1}.  \lbeq(dd)
\eqn 
We have $\tilde{d}(\lam)= T_{22}^{-1}+ \Og_2(g(\lam)^{-1})$ 
for small $\lam>0$ and entries of $\tilde{d}(\lam)$ are good 
multipliers. 
The following lemma follows immediately from Feshbach \reflm(FS).

\bglm \lblm(eigen-lm) $\tRg_1(\lam)$ is invertible in $S_2\HL$ and 
\begin{gather}\lbeq(R1-inv-matrix)
\tRg_1(\lam)^{-1}
= g(\lam) S_3 T_3^{-1} S_3 + L_4(\lam), \\
\lbeq(L4-def)
L_4(\lam)= \begin{pmatrix}  \tilde{d}(\lam) & -\tilde{d}(\lam) \tT_{23} \tT_{33}^{-1} \\
-\tT_{33}^{-1} \tT_{32} \tilde{d}(\lam)  
&  \tT_{33}^{-1}\tT_{32}\tilde{d}(\lam) \tT_{23} \tT_{33}^{-1}
\end{pmatrix}.
\end{gather} 
The entries of $L_4(\lam)$ are good multipliers.
\edlm 

By virtue of \reflm(eigen-lm) the proof of \refprop(impo-1-2) 
goes through with slight modifications 
(see \refrmss(B2-inv-a,A12-good,prop-517)) and the result 
of the proposition that 
$ vM(\lam)^{-1}v \equiv 
- g(\lam)^{-1} \lam^{-2}vS_2 \tRg_1(\lam)^{-1} S_2v$ 
remains to hold if $H$ has singularities 
of the third kind at zero. Thus, we have only to show that 
${\tilde{\W}}_{{\rm low},2a}$ defined by \refeq(tWlow-a-1) 
satisfies the statement of \refth(p-wave). 

We take the basis $\{\z_1, \dots, \z_m\}$ of normalized vectors 
of $S_2\HL$ such that 
\[
T_2 \z_j = -\k_j^2 \z_j, \ \ \k_j >0, \ j=1,2
\]
and $\{\z_{l+1}, \dots, \z_m\}$ is the basis of $\Xg_3$ such that 
\[
T_3 \z_k =  \a_k \z_k, \quad k=l+1, \dots, m. 
\]
Since $\la \z_j, \z_k\ra= -\k_j^{-2}\la T_2\z_j, \z_k\ra= 
-\k_j^{-2}\la \z_j, T_2 \z_k\ra= 0$ also for
for $=1,2$ and $k=l+1,\dots, m$,  $\{\z_1, \dots, \z_m\}$ is orthonormal. 
Since $T_2\z_k=0$, we have  the extra cancellation property:
\bqn \lbeq(cancellation-2)
\int_{\R^2} x_1 \z_k(x)v(x)dx  = \int_{\R^2} x_2 \z_k(x)v(x)dx  = 0. 
\quad k=l+1, \dots, m.
\eqn 

The following lemma proves the positive part of \refth(third-1). 

\bglm \lblm(third-positive) For any $a>0$, 
${\tilde{\W}}_{{\rm low},2a}$ is bounded in $L^p(\R^2)$ for $1<p\leq 2$. 
\edlm 
\bgpf Let $\tRg_1^{-1}=g(\lam) S_3 T_3^{-1} S_3 + L_4(\lam)$ 
be represented by matrix $D(\lam)=(d_{jk}(\lam))$ 
with respect to the basis $\{\z_{1}, \dots, \z_m\}$. 
$D(\lam)$ is the sum of $\tilde{D}(\lam)$ for 
$g(\lam) S_3 T_3^{-1} S_3$ and $B(\lam)=(\b_{jk}(\lam))$ for $L_4(\lam)$. 
By virtue of \refeq(L4-def) and 
\refeq(dd), $\b_{jk}(\lam)\chi_{\leq 2a}(\lam)$ are good multipliers. 
Recall that $\la \z_j, v \ra=0$ $for j=1, \dots, m$. 
Let $\tilde{\W}^B$ and $\tilde{\W}^D$ be defined by \refeq(tWlow-a-1) 
respectively with $L_4(\lam)$ and $g(\lam) S_3 T_3^{-1} S_3$ replacing 
$S_2 \tRg_1(\lam)^{-1}S_2$ so that 
$\tilde{\W}_{{\rm low},2a}= \tilde{\W}^B + \tilde{\W}^D$. Then, 
$\tilde{\W}^B u$ is the sum over $1\leq j,k \leq m$ of 
\bqn \lbeq(tWlow-B-jk)
\tilde{\W}^B_{jk} u=\int_0^\infty g(\lam)^{-1} \lam^{-2}\b_{jk}(\lam)
|G_0(-\lam)v\z_j \ra \la \z_k v, 
\Pi({\lam})u \ra \chi_{\leq 2a}(\lam)\lam d\lam. 
\eqn 
This is the same as \refeq(tWlow-1) with $\b_{jk}(\lam)$ replacing   
by $d_{jk}(\lam)$ and the former function can play the role of the latter 
in the proof of the positive part of \refth(p-wave). Thus, 
$\tilde{\W}^B$ is bounded in $L^p(\R^2)$ for $1<p<2$. 

Since $\tilde{D}(\lam)=g(\lam){\rm diag}(\a_1, \dots, \a_m)$ is 
diagonal and $\a_1=\a_2=0$,   
\[
\tilde{\W}^Du(x)=
\sum_{j=3}^m \a_j \int_0^\infty \lam^{-1}(G_0(-\lam)\z_j v)(x) 
\la \z_j v, \Pi(\lam) u\ra \chi_{\leq 2a}(\lam)d\lam.
\]
This is of the form \refeq(v-w) but with the 
singular factor $\lam^{-1}$ in place of $\lam$. However,  
for $j=3, \dots, m$ the extra cancellation property 
\refeq(cancellation-2) implies that 
$\la \z_j v, \tilde{b}(\lam,\cdot)\ra=0$ 
and $\la \z_j v, \Pi(\lam)u \ra = \la \z_j v, \tilde{g}(\lam,\cdot)\ra$.
It follows that $\tilde{\W}^Du(x)$ is equal to  
$\sum_{l,k=1}^2 \int_0^1 (1-\th) d\th$ of 
\[
- \int_{\R^4} z_l z_k (\z_j v)(y) (\z_j v)(z)
(\t_y K \chi_{\leq 2a}(|D|)R_k R_l 
\t_{-{\th}z}u)(x)dy dz . 
\] 
Then Mikowski's inequality, \refeq(K-est) 
and the multiplier theorem imply that for any $\c>4$  
\[
\|\tilde{\W}^Du \|_p \leq C\sum_{j=3}^m 
\|z^2 \z_j v\|_1 \|\z_j v\|_1 \|u\|_p 
\leq  C \|\az^{\c} V\|_1 \|u\|_p, \quad 1<p<\infty 
\]
and $\tilde{\W}^D$ is a good operator. This proves the lemma.   
\edpf 

If $T_2=0$, then $\tRg_1(\lam)= g(\lam)S_3 T_3^{-1} S_3 $ 
does not contain 
$L_4(\lam)$ and the proof of \reflm(third-positive) above implies 
the following theorem. 

\bgth \lbth(third-2) 
Suppose that $H$ has singularities of the third kind at zero 
and $T_2=0$, viz. $S_2=S_3$. Then $W_{+}$ is a good operator. 
\edth

The following lemma completes the proof of \refth(third-1)

\bglm 
Suppose that $H$ has singularities of the third kind at zero 
and $S_3 \subsetneq S_2$. Then $W_{+}$ is unbounded in $L^p(\R^2)$ 
if $2<p<\infty$. 
\edlm 
\bgpf 
In view of the proof of \reflm(third-positive), it suffices 
to prove that 
$\tilde{\W}^B \stackrel{\rm def}{=} \sum_{j,k=1}^m \tilde{\W}^B_{jk}$  
is unbounded in $L^p(\R^2)$ if $2<p<\infty$.  
Since $\la \z_j, v\ra=0$, we may replace 
$\Pi(\lam) u(x)$ of \refeq(tWlow-B-jk) by 
$\Pi(\lam) u(x)-\Pi(\lam) u(0)= \tilde{g}(\lam,x)+ \tilde{b} (\lam,x)$ 
as previously, which produces  
$\tilde{\W}^Bu = {\tilde{\W}}^{B}_{g}u+ {\tilde{\W}}^B_{b}u$, 
where the definition of 
${\tilde{\W}}^B_{g}$ and ${\tilde{\W}}^B_{b}$ should be obvious. 
Then, the factor $\lam^2$ in 
$\tilde{g}(\lam,x)$ cancels the singularity in the integrand of 
\refeq(tWlow-B-jk) and ${\tilde{\W}}^{B}_{g}$ becomes a good operator 
as in \refprop(good). Then, since $\chi_{>4a}(|D|)$ is a good operator, 
it suffices to show that $\chi_{>4a}(|D|){\tilde{\W}}^B_{b}$ 
is unbounded in $L^p(\R^2)$ if $2<p<\infty $. If we further 
express $\chi_{>4a}(|D|)G_0(-\lam)(v\z_j)$ as the sum as in \refeq(6-1), 
then the proof of \reflm(low-high-positive) shows that the seccond term on 
the right of \refeq(6-1) produces a good operator and 
we have only to prove that  
\bqn  \lbeq(redu-1)
Z_a u(x)= \sum_{j,k=1}^m (\hat{\m}_a\ast (v\z_j)(x) 
\int_0^\infty   \b_{jk}(\lam)
\la\z_k{v}, \tilde{b}(\lam,\cdot)\ra \r_a(\lam)d\lam
\eqn  
is unbounded in $L^p(\R^2)$ for $2<p<\infty$. Here 
$\m_a(\xi)=\chi_{>4a}(|\xi|)|\xi|^{-2}$ satisfies 
$\hat{\mu}_a (x)\in L^p(\R^2)$ for $1<p<\infty$ 
and we have introduced the notation 
$\r_a(\lam)= g(\lam)^{-1}\lam^{-1}\chi_{2a}(\lam)$. 
To make the argument transparent, we introduce 
the vector notation 
\[
\Zg = \begin{pmatrix} \Zg_1, \\ \Zg_2   \end{pmatrix}, \quad 
\Zg_1 = \begin{pmatrix} \z_1, \\ \z_2 \end{pmatrix} \quad 
\Zg_2 = \begin{pmatrix} \z_3 \\ \vdots \\ \z_m  \end{pmatrix} \quad 
\Zg^\ast = (\Zg_1^\ast, \Zg_2^\ast)
\]
and express $Z_a u(x)$ in the form 
\bqn \lbeq(redu-2)
Z_a u(x) = \Big\la \hat{\m}_a \ast v(x)\Zg(x),  
\int_0^\infty  B(\lam)
\la \Zg v, \tilde{b}(\lam,\cdot)\ra \r_a(\lam)d\lam \Big\ra_{\C^2}.
\eqn 
Here the extra moment condition \refeq(cancellation-2) and 
\refeq(bad-def):  
\[
\tilde{b}(\lam, z)= \sum_{l=1}^2 z_l b_l(\lam), \quad 
b_l(\lam)= \frac{i}{2\pi}\int_{{\mathbb S}^1}
\w_l \hat{u}(\lam\w) d\w 
\]
imply that $\la \Zg_2 v, \tilde{b}(\lam,\cdot)\ra=0$. It follows that    
\[
B(\lam) \la \Zg v, \tilde{b}(\lam,\cdot)\ra = 
\begin{pmatrix} {\bf 1}_2 \\ A \end{pmatrix} \tilde{d}(\lam) 
\la v\Zg_1, \tilde{b}(\lam, \cdot) \ra\, , 
\]
where ${\bf 1}_2$ for the $2\times 2$ identity matrix,  
$A=(a_{jk})_{3\leq j \leq m, 1\leq k \leq 2}$ is the reprentation matrix for 
$\tilde{T}^{-1}_{33}\tilde{T}_{32}$ and $\tilde{d}(\lam)$ is defined in 
\refeq(dd). Since the matrix $\begin{pmatrix} {\bf 1}_2 \\ A \end{pmatrix}$ 
does not depend on $\lam$, we can transpose it in front of  
$\hat{\m}_a \ast v(x)\Zg(x)$ in \refeq(redu-2) and express $Z_a u(x)$ 
in the form  
\bqn 
Z_a u(x) =  f_{1}(x)\tilde{\ell}_1(u) + f_{2}(x)\tilde{\ell}_2(u)
\eqn 
where $f_j(x)$ and $\tilde{\ell}_j(u)$, $j=1,2$ are defined by 
\[
\begin{pmatrix} f_1(x) \\ f_2(x) \end{pmatrix}= 
\hat{\m}_a \ast \begin{pmatrix} w_1(x) \\ w_2(x) \end{pmatrix}, 
\quad 
\begin{pmatrix} w_1(x) \\ w_2(x) \end{pmatrix} 
= v(x)\Zg_1(x) + A^\ast v(x) \Zg_2(x) 
\]
and 
\[
\begin{pmatrix} \tilde{\ell}_1(u) \\ \tilde{\ell}_2(u)
\end{pmatrix}=\int_0^\infty  
\tilde{d}(\lam)\la \Zg_1 v, \tilde{b}(\lam,\cdot)\ra \r_a(\lam)d\lam
\]
and both depend on $a>0$. 
It is obvious that $f_1, f_2 \in L^p(\R^2)$ for $1<p<\infty$ 
and we show that they are linearly independent if $a>0$ is 
sufficiently small. Indeed if otherwise, there exists for any $a>0$ 
a point $\begin{pmatrix} C_{1a} \\ C_{2a}\end{pmatrix} \in {\mathbb S}^1$ 
such that $C_{1a}f_1(x)+ C_{2a}f_2(x)=0$ and via Fourier transform 
\bqn \lbeq(suff)
\Fg(C_{1a}w_1 + C_{2a} w_2)(\xi)= 0 \quad \mbox{for}\ |\xi|\geq 4a
\eqn 
Then the subset $\Cg_a \subset {\mathbb S}^1$ of points 
$\begin{pmatrix} C_{1a}\\ C_{2a}\end{pmatrix}$ 
which satisfy \refeq(suff) is compact and non-empty and 
$\Cg_a\subset \Cg_b$ if $a<b$. 
Thus, $\Cg=\cap_{a>0}\Cg_a\not= \emptyset$ and, for a  
$\begin{pmatrix} C_{1}\\ C_{2}\end{pmatrix}\in \Cg$,  
we must have 
\[ 
\Fg(C_{1} w_1 + C_{2}w_2)(\xi) = 0, \quad \xi \in \R^2
\]
and hence $C_{1} w_1(x) + C_{2}w_2(x)=0$ or  
\[
v(C_{1} \z_1 + C_{2}\z_2 + 
(C_{1}a_{13} + C_{2}a_{23}) \z_3 + \cdots+ 
(C_{1}a_{1m} + C_{2}a_{2m}) \z_m )=0 
\]
However, 
$-\k_j^2 \z_j(x) = v(x)(G_1 v \z_j)(x)$, $j=1, 2$  
and $\a_j \z_j(x) = v(x)(G_2 v \z_j)(x)$  for $j=3, \dots, m$. 
Hence $\z_j(x) = 0$, $j=1, \dots, m$, if $v(x)=0$.    
It follows that 
\[
C_{1} \z_1 + C_{2}\z_2 + 
(C_{1}a_{13} + C_{2}a_{23}) \z_3 + \cdots+ 
(C_{1}a_{1m} + C_{2}a_{2m}) \z_m =0 
\]
and it must be that $C_{1}=C_2=0$ which is a contradiction. Thus, 
$f_1$ andd $f_2$ must be linearly independent for small $a>0$. 
Incidentaly, this simultaneously proves the corresponding statement in the proof 
of \refprop(low-high-negative). 

The rest of the proof is the repetition of that of 
\refprop(low-high-negative). If $Z_a u$ were  
bounded in $L^p(\R^2)$ for $2<p<\infty$, then 
$\ell_1(u)$ and $\ell_2(u)$ must be bounded functionals in $L^p(\R^2)$ 
by the Hahn-Banach theorem and we must have \refeq(6-100) with 
$\tilde{d}(|\xi|)$ in place of $D(|\xi|)$, viz 
\bqn \lbeq(6-101)
\frac{\chi_{\leq 2a}(|\xi|)}{g(|\xi|)|\xi|^{2}}
\tilde{d}(|\xi|)\Ng(\xi)\in L^p(\R^2), \quad 
\Ng(\xi)=\begin{pmatrix} \la z_1 v, \z_1\ra \xi_1+ \la z_2 v, \z_1\ra\xi_2 \\ 
\la z_1 v, \z_2\ra \xi_1+ \la z_2 v, \z_2\ra \xi_2  \end{pmatrix}\,.
\eqn 
Since $\tilde{d}(|\xi|)^{-1}$ is bounded (see \refeq(dd)), \refeq(6-101) 
would lead to $\la z_1v, \z_k\ra = \la z_2v, \z_{k}\ra=0$ for $k=1,2$ 
which contradicts to \refeq(G1-inner). The lemma follows. 
\edpf

\end{document}